\pdfoutput=1
\documentclass[11pt, oneside]{article}   	
\usepackage{geometry}
\usepackage{graphicx}			
\usepackage{amsmath,bm}				
\usepackage{amssymb}
\usepackage{xcolor}
\usepackage{rotating}
\usepackage{caption}
\usepackage{subcaption}

\usepackage{framed}

\begin{document}
\baselineskip=2pc

\vspace{.5in}

\begin{center}

{\large\bf
Fast sparse grid simulations of fifth order WENO scheme
for high dimensional hyperbolic PDEs
\footnote{Research was partially supported by NSF grant DMS-1620108.}
}

\end{center}

\vspace{.1in}

\centerline{Xiaozhi Zhu\footnote{Department of Applied and Computational Mathematics and Statistics,
University of Notre Dame, Notre Dame, IN 46556, USA. E-mail: xzhu4@nd.edu
},
Yong-Tao Zhang\footnote{Corresponding author. Department of Applied and Computational Mathematics and Statistics,
University of Notre Dame, Notre Dame, IN 46556, USA. E-mail: yzhang10@nd.edu
}
}

\vspace{.6in}

\abstract{The weighted essentially non-oscillatory (WENO) schemes, especially the fifth order WENO schemes,
are a popular class of high order accurate numerical methods for solving hyperbolic
partial differential equations (PDEs).
However when the spatial dimensions are high, the number of spatial grid points increases significantly. It leads to
large amount of operations and computational costs in the numerical simulations by using nonlinear high order accuracy WENO schemes such as a fifth order WENO scheme.
How to achieve fast simulations by high order WENO methods for high spatial dimension hyperbolic PDEs is a challenging and important question.
In the literature, sparse-grid technique has been developed as a very efficient approximation tool for high dimensional problems. In a recent work [Lu, Chen and Zhang, Pure and Applied Mathematics Quarterly, 14 (2018) 57-86], a third order finite difference WENO method with sparse-grid
combination technique was designed to solve multidimensional hyperbolic equations including both linear advection equations and nonlinear Burgers' equations.
Numerical experiments showed that WENO computations on sparse grids achieved comparable third order accuracy in smooth regions of the solutions and nonlinear stability as that for computations on regular single grids. In application problems, higher than third order WENO schemes are often preferred in order to efficiently resolve the complex solution structures. In this paper, we extend the approach to higher order WENO simulations specifically the fifth order WENO scheme. A fifth order WENO interpolation is applied in the prolongation part of the sparse-grid
combination technique to deal with discontinuous solutions. Benchmark problems
are first solved to show that significant CPU times are saved while
both fifth order accuracy and stability of the WENO scheme are preserved for simulations on sparse grids.
The fifth order sparse grid WENO method is then applied to kinetic problems modeled by high dimensional Vlasov based PDEs to further demonstrate
large savings of computational costs by comparing with simulations on regular single grids.}

\vfill


\pagenumbering{arabic}

\newpage

\section{Introduction}

A popular class of high order accuracy numerical methods for solving hyperbolic PDEs is the class of weighted essentially non-oscillatory (WENO) schemes. They have been applied extensively in computational fluid dynamics and other scientific problems. High order WENO schemes are especially efficient for solving problems containing both singularities and complicated smooth solution structures \cite{SZS},
for example, the Rayleigh-Taylor instability problems \cite{ZSSZ,ZSZ}, the shock vortex interactions \cite{ZZS3}, biological population dynamics \cite{ZZS}, etc.
A third order finite volume WENO scheme was firstly constructed
in \cite{LOC}. In \cite{JS}, Jiang and Shu designed
arbitrary order accurate finite
difference WENO schemes for efficiently computing multidimensional problems. A general
framework for the design of the smoothness indicators and nonlinear
weights was provided. In WENO schemes, a weighted combination
of several local reconstructions based on different stencils which are called ``small stencils'' is used as the final WENO
reconstruction. The combination coefficients are called
nonlinear weights. Smoothness indicators measure the smoothness of
reconstructed functions in the relevant small stencils.
They are incorporated in nonlinear weights to achieve both nonlinear stability in non-smooth regions and
high order accuracy in smooth regions of the solution for WENO schemes.
In classical WENO schemes \cite{JS}, nonlinear weights are also designed to
increase the order of accuracy over that on each small stencil.
WENO schemes have been studied extensively in the literature.
For example, to deal with complex domain geometries, WENO schemes on unstructured meshes were developed in e.g. \cite{HS,ZS,LNSZ,DK2,ZS2,LZ,ZhuS}. To efficiently solve steady state problems of hyperbolic PDEs by high order WENO schemes, fast sweeping WENO schemes and homotopy WENO schemes were developed in \cite{ZZQ,ZZC,XZZS,HHSSXZ,WZZS}. High order Krylov implicit integration factor methods \cite{CZ} were applied to WENO schemes in \cite{JiangZhang, JiangZhang2, LuZhang2} for solving stiff convection-diffusion-reaction PDEs. Efforts have been made to simplify or improve the accuracy in high order WENO schemes.
Strategies include modifying the linear or nonlinear weights, modifying the smoothness indicators, etc, see e.g. \cite{LPR, HAP, YC, CCD, AMM, ZhuQ, ZhuS}. For an overview on WENO schemes, see e.g. \cite{Shu2, ZSNA}.

Because WENO schemes deal with problems with both complicated solution structures and discontinuities / sharp gradient regions,
their sophisticated nonlinear properties and high order accuracy require more operations than that in
many other schemes. The computational cost increases significantly when the number of grid points is large or the spatial dimensions of the PDEs are high. Especially for long time simulations, how to achieve fast computations by high order WENO methods is a challenging and important question.

In the literature, sparse-grid techniques have been developed as an efficient approximation tool for solving high-dimensional problems in scientific and engineering applications. Discretizations on sparse grids involve
$O(N\cdot (\log N)^{d-1})$ degrees of freedom only, where $d$ denotes the dimensionality of the underling problems and
$N$ is the number of grid points in one coordinate direction.
See \cite{BG} for a detailed review on sparse-grid techniques.
Sparse-grid techniques were introduced by Zenger \cite{Zg} in 1991 to reduce the number of degrees of freedom in finite element method. The sparse-grid combination technique, which was introduced in 1992 by Griebel et al. \cite{GSZ}, is one of the approaches about practical implementation of sparse-grid techniques.
In the sparse-grid combination technique, the final solution is a linear combination of solutions on semi-coarsened grids, where the coefficients of the combination are chosen such that there is a canceling in leading-order error terms and the accuracy order can be kept to be the same as that on a single full grid \cite{LKV1,LKV2,GSZ}.
Recently in \cite{LuZhang1}, the sparse-grid combination technique has been used in Krylov implicit integration factor methods
to efficiently solve high spatial dimension convection-diffusion equations.
Our next goal is to apply sparse-grid techniques in high order WENO schemes for solving multidimensional hyperbolic PDEs which may develop singular solutions and obtain much faster computations than that in their regular performance.
The challenge is how to implement WENO computations on sparse grids such that comparable high order accuracy of WENO schemes in smooth regions and
essentially non-oscillatory stability in non-smooth regions of the solutions can be
 preserved as that for computations on regular single grids. This is not straightforward due to the high nonlinearity of high order WENO schemes.
The first effort we have made is in \cite{LuChenZhang3}, where a third order finite difference WENO method with sparse-grid
combination technique was designed to solve multidimensional hyperbolic equations including both linear advection equations and nonlinear Burgers' equations. Numerical experiments showed that WENO computations on sparse grids achieved comparable third order accuracy in smooth regions of the solutions as that for computations on regular single grids, and nonlinear stability of the scheme was maintained. In application problems, higher than third order WENO schemes are often preferred in order to efficiently resolve the complex solution structures. Especially the fifth order WENO (WENO5) schemes are very popular and have been used broadly. It is an open problem whether fast sparse grid simulations for higher order WENO schemes such as a fifth order WENO scheme can be performed as that for lower order WENO schemes, especially whether the accuracy order and nonlinear stability of higher order WENO schemes can be preserved in sparse grid simulations. This is not a trivial question since theoretical analysis of higher order
WENO schemes is difficult due to their high nonlinearity and numerical experiment is one of the major tools to study them.

In this paper, we apply the sparse-grid combination technique to a fifth order WENO finite difference scheme for solving hyperbolic PDEs defined on high spatial dimension domains.
To deal with discontinuity / sharp gradient in solutions of hyperbolic PDEs, we apply a fifth order WENO interpolation for the prolongation part in sparse-grid
combination technique. Two dimensional (2D), three dimensional (3D) and four dimensional (4D) numerical examples with smooth or non-smooth solutions are presented to show that significant computational times are saved, while both accuracy and stability of the fifth order WENO scheme are maintained for simulations on sparse grids. The fifth order sparse grid WENO method is then applied to kinetic problems modeled by high dimensional Vlasov based PDEs to further demonstrate
large savings of computational costs by comparing with simulations on regular single grids.
The rest of the paper is organized as following. In Section 2, we describe the algorithm how to apply the sparse-grid combination technique to a WENO5 scheme, with a novel WENO5 prolongation. In Section 3, various numerical experiments including solving high dimensional kinetic problems modeled by Vlasov based PDEs
are presented to test the sparse grid WENO5 method and show significant savings in
computational costs by comparisons with single-grid computations. Conclusions and discussions are given in Section 4.

\section{Description of the numerical method}
We study efficient and high order accuracy numerical methods for solving multidimensional hyperbolic equations
\begin{equation}
\label{eq1}
{u}_t+\nabla\cdot\vec{f}({u})=0.
\end{equation}
Here ${u}(\vec x,t)$ is the unknown function, and $\vec{f}=({f}_1, \cdots, {f}_d)^T$ is the vector of
flux functions. $d$ is the dimension of the spatial domain on which the PDE is defined.
We apply
the method of lines (MOL) approach, namely,
the equation (\ref{eq1}) is first discretized in the spatial directions to obtain a system of ordinary differential equations (ODEs), then the ODE system is marched by a temporal scheme. In this paper, we use the classical fifth order finite difference WENO scheme
with Lax-Friedrichs flux splitting \cite{JS} for the spatial discretizations, and show that very efficient computations are achieved by performing such high order WENO simulations on sparse grids. However, we like to point out that such
approach can be easily applied to other recently developed high order WENO schemes, for example \cite{ZhuS0, ZhuQ, CCD}, etc.
In the following,
the classical fifth order finite difference WENO spatial discretization is reviewed at first, then we give a detailed
description about how to implement it on sparse grids by
the sparse-grid combination techniques with a fifth order WENO prolongation. A complete algorithm is summarized at last.

\subsection{The fifth order WENO discretization}
The conservative finite difference scheme is used to discretize the hyperbolic equation (\ref{eq1}). The point values of unknown functions in the PDEs are approximated at a
uniform (or smoothly varying) grid. One of big advantages of the finite difference schemes is that they
approximate multi-dimensional derivatives in a dimension by dimension way, hence they are very efficient for solving multi-dimensional problems defined on regular domains \cite{SO}. In the following
we will just describe the discretization of derivatives for
one spatial direction. Similar
procedures are used to discretize derivatives of other spatial directions.
Without the loss of generality,
we take the $x$-direction derivative $f(u)_x$ as the example. Its value at a grid point with $x$-coordinate $x_i$ on a uniform grid with $x$-direction grid size $\Delta x$
is approximated by a conservative flux difference
\begin{equation}
\label{eq2.9}
f(u)_x|_{x=x_i}\approx \frac{1}{\Delta x}(\hat f_{i+1/2}-\hat f_{i-1/2}).
\end{equation}
Here $\hat f_{i+1/2}$ is the numerical flux at the point $x_{i+1/2}$, where $x_{i+1/2}=(x_i+x_{i+1})/2$.
For the classical fifth order WENO scheme, if the wind is positive (namely, when $f'(u)\geq 0$ for a scalar
equation, or when the corresponding eigenvalue is positive for a
system of equations with a local characteristic decomposition), the numerical flux $\hat
f_{i+1/2}$ depends on a five-point stencil (called the big stencil) with numerical values $f(u_{l})$,
$l=i-2,i-1,i,i+1,i+2$.
Note that for the purpose of
the simplicity of notations, here we use $u_l$ to denote the value of
the numerical solution $u$ at the grid point $x=x_l$ along the grid lines of other spatial directions, e.g.,
$y=y_j, z=z_k,$ etc, with the understanding that the value could be
different for different coordinates of other spatial directions. In
the classical fifth order WENO scheme,
the numerical flux $\hat f_{i+1/2}$ is computed using a convex combination
of three third order approximations of the numerical flux. These third order approximations are based on three different substencils (i.e., small stencils) of three grid points each. The union of these three substencils is the original big stencil.
The coefficients of the convex combination, called ``nonlinear weights'',
depend on ``smoothness indicators'' which measure the regularity of
the numerical solution in each substencil. The detailed formulae are
\begin{equation}
\label{eq2.7}
\hat f_{i+1/2}=w_0\hat f^{(0)}_{i+1/2}+w_1\hat f^{(1)}_{i+1/2}+w_2\hat f^{(2)}_{i+1/2},
\end{equation}
where
\begin{eqnarray}
\label{eq2.7.1}
\hat f^{(0)}_{i+1/2} & = & \frac{1}{3}f(u_{i-2})-\frac{7}{6}f(u_{i-1})+\frac{11}{6}f(u_{i}),
\nonumber \\
\hat f^{(1)}_{i+1/2} & = & -\frac{1}{6}f(u_{i-1})+\frac{5}{6}f(u_{i})+\frac{1}{3}f(u_{i+1}),
\\
\hat f^{(2)}_{i+1/2} & = & \frac{1}{3}f(u_{i})+\frac{5}{6}f(u_{i+1})-\frac{1}{6}f(u_{i+2}).
\nonumber
\end{eqnarray}
Here
\begin{equation}
w_r=\frac{\alpha_r}{\alpha_1+\alpha_2+\alpha_3}, \qquad
\alpha_r=\frac{d_r}{(\epsilon+\beta_r)^2}, \qquad r=0, 1, 2.
\label{eq2.6}
\end{equation}
$d_0=0.1, d_1=0.6, d_2=0.3$ are called the ``linear weights'', and
$\beta_0, \beta_1, \beta_2$ are
 the smoothness indicators. They are actually quadratic functions of numerical values $f(u_{l})$ on the substencils and have the following explicit formulae
\begin{eqnarray}
\label{eq2.7.2}
\beta_0 & = & \frac{13}{12}(f_{i-2}-2f_{i-1}+f_{i})^2+\frac{1}{4}(f_{i-2}-4f_{i-1}+3f_{i})^2, \nonumber \\
\beta_1 & = & \frac{13}{12}(f_{i-1}-2f_{i}+f_{i+1})^2+\frac{1}{4}(f_{i-1}-f_{i+1})^2,  \\
\beta_2 & = & \frac{13}{12}(f_{i}-2f_{i+1}+f_{i+2})^2+\frac{1}{4}(3f_{i}-4f_{i+1}+f_{i+2})^2,
\nonumber
\end{eqnarray}
where $f_{l}$ denotes $f(u_{l})$.
$\epsilon$ is a small
positive number chosen to avoid the denominator becoming to $0$.

For the case that the wind is negative (namely, when $f'(u)< 0$), the right-biased
big stencil with numerical values $f(u_{i-1}), f(u_{i}), f(u_{i+1}), f(u_{i+2})$ and
$f(u_{i+3})$ is used to reconstruct a fifth order WENO
approximation to the numerical flux $\hat f_{i+1/2}$. The formulae
 for the negative and the positive wind cases are symmetric with respect to the
point $x_{i+1/2}$. For the general case of the flux $f(u)$, a flux splitting is performed to
separate the positive and the negative wind parts.
Here we use the
``Lax-Friedrichs flux splitting"
\begin{equation}
f^+(u)=\frac{1}{2}(f(u)+\alpha u), \qquad
f^-(u)=\frac{1}{2}(f(u)-\alpha u),
\label{eqLFflux}
\end{equation}
where $\alpha=\max_u|f'(u)|$. $f^+(u)$ is the positive wind part,
and $f^-(u)$ is the negative wind part. Then the corresponding WENO
approximations are applied in constructing numerical fluxes $\hat
f^+_{i+1/2}$ and $\hat f^-_{i+1/2}$ respectively. The final numerical flux
$\hat{f}_{i+1/2}=\hat f^+_{i+1/2}+\hat f^-_{i+1/2}$.
More details can be found in the review articles e.g. \cite{Shu,ZSNA}.

\subsection{The sparse grid WENO scheme}

In order to improve the efficiency of using high order WENO schemes to solve multidimensional hyperbolic equations (\ref{eq1}), we study the techniques about how to implement WENO schemes on sparse grids. One effective method is to use sparse-grid combination approach. In this paper, we focus on the classical fifth order finite difference WENO scheme described in the last section.

In the sparse-grid combination technique, a computational domain is partitioned into a group of sparse grids which have different grid sizes. Among them the most refined mesh is corresponding to the usual single full grid in our regular single-grid computations. The PDEs are solved on certain semi-coarsened grids rather than the single full grid. Since the sum of the numbers of grid points of the semi-coarsened grids on which the PDEs are solved is much smaller than that of the single full grid, computational costs are saved a lot. At last the solutions on these semi-coarsened grids are
combined to obtain the final solution on the most refined mesh.
The final solution obtained by a good sparse-grid combination technique is expected to
have a comparable accuracy to that resulted from solving the PDEs directly on a single full grid. For example see \cite{GSZ,LKV1,LKV2,LuZhang1,LuChenZhang3} for some early work on this approach.

\begin{figure}
    \centering
    \includegraphics[width=4in]{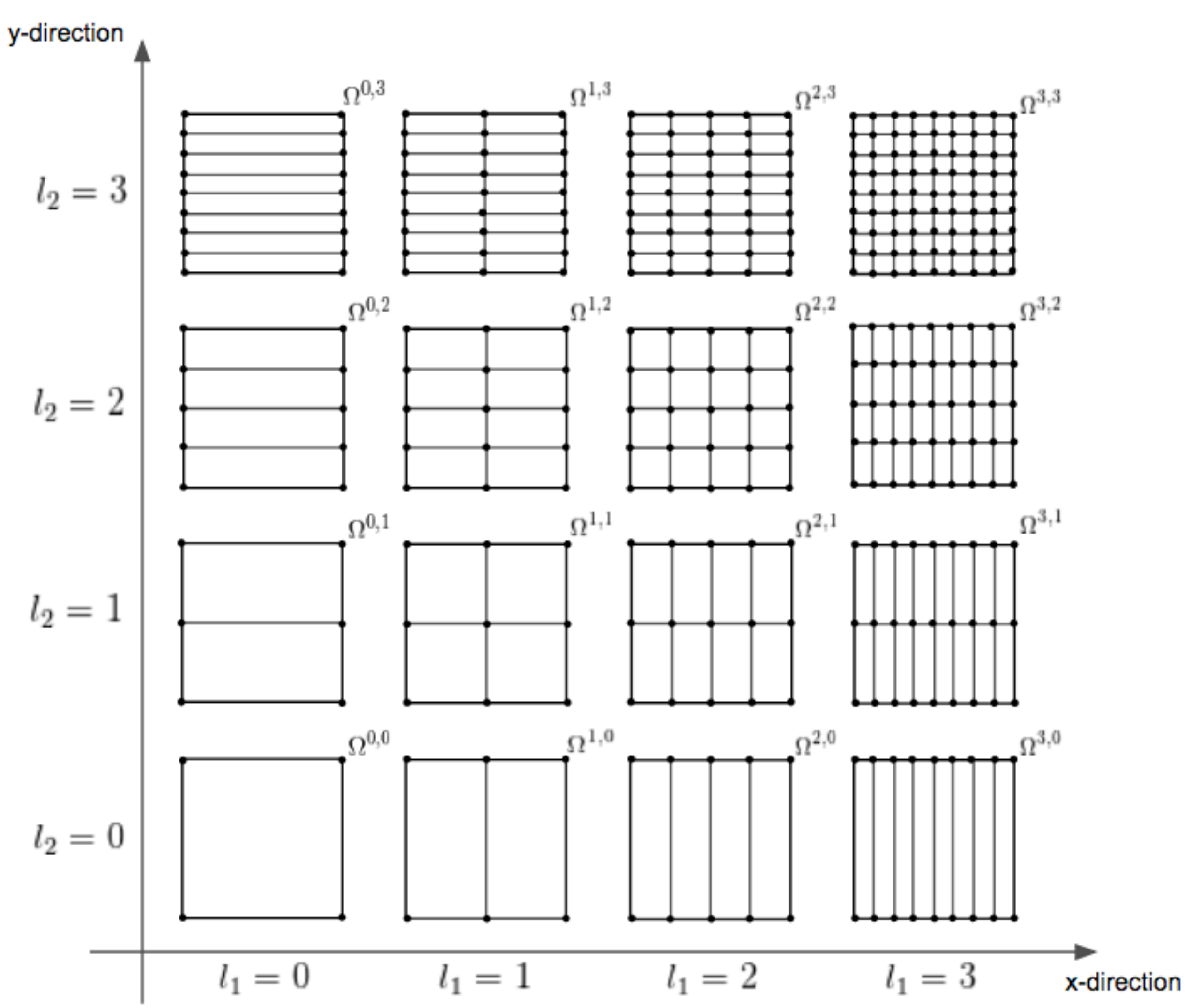}
    \caption{Illustration of two dimensional sparse grids $\{\Omega^{l_1,l_2}\}$ for one cell of a root grid. Here the cell indicated by $\Omega^{0,0}$ at the level $l_1=0, l_2=0$ is one cell of the whole root grid $\Omega^{0,0}$, and the side length of the cell is $H$. The finest level $N_L=3$.}
\label{sparse_grids}
\end{figure}

Without the loss of generality, we use the two dimensional (2D) case to illustrate the idea. Algorithm procedures for the higher dimensional cases are similar. Consider a two dimensional computational domain $[a, b]^2$. Note that for the simplicity of the description, here we use a square domain. However, it is not necessary that the domain needs to be a square one. The procedures presented here can be applied to any rectangular domain straightforwardly.
The semi-coarsened grids are obtained as the following. First the domain is partitioned into the coarsest mesh with $N_r$ cells in each spatial direction. It is called a root grid and denoted by $\Omega^{0,0}$. The grid size of the root grid is $H=\frac{b-a}{N_r}$. Then a multi-level refinement on the root grid is carried out to obtain a family of semi-coarsened grids \{$\Omega^{l_1,l_2}$\}. These semi-coarsened grids \{$\Omega^{l_1,l_2}$\} have mesh sizes $h_{l_1}=2^{-l_1}H$ in the $x$ direction and $h_{l_2}=2^{-l_2}H$ in the $y$ direction, where $l_1=0,1,\cdots,N_L$, $l_2=0,1,\cdots,N_L$. The superscripts $l_1, l_2$ represent the levels of refinement relative to the root grid $\Omega^{0,0}$ in the $x$ and the $y$ directions respectively, and $N_L$ denotes the finest level. Hence, the finest grid here is $\Omega^{N_L,N_L}$ with the grid size $h=2^{-N_L}H$ for both $x$ and $y$ directions. In Figure \ref{sparse_grids}, we show a family of two dimensional semi-coarsened sparse grids \{$\Omega^{l_1,l_2}$\} for one cell of a root grid, with the finest level $N_L=3$.

We apply the spare-grid combination techniques to the fifth order WENO scheme described in the last section for
solving multidimensional hyperbolic equations (\ref{eq1}).
 For the time discretization, the third order total variation diminishing (TVD) Runge-Kutta (RK) scheme in \cite{SO,GST} is used. The PDE (\ref{eq1}) is not solved on a single full grid, but on the following $(2N_L+1)$ sparse grids $\{\Omega^{l_1,l_2}\}_I$:
\[
\Big\{
\Omega^{0,N_L}, \Omega^{1,N_L-1}, \cdots, \Omega^{N_L-1,1}, \Omega^{N_L,0}
\Big\}\quad
\text{and} \quad
\Big\{
\Omega^{0,N_L-1}, \Omega^{1,N_L-2}, \cdots, \Omega^{N_L-2,1}, \Omega^{N_L-1,0}
\Big\}.
\]
Here the notation $I$ represents the index set
\[
I=\big\{(l_1,l_2)| l_1+l_2=N_L \quad \text{or} \quad l_1+l_2=N_L-1\big\}.
\]
By performing time marching of numerical solutions of the PDE (\ref{eq1}) using the TVD RK scheme with the fifth order WENO spatial discretization on these $(2N_L+1)$ sparse grids, we have $(2N_L+1)$ sets of numerical solutions $\{U^{l_1,l_2}\}_I$, where each set of numerical solutions is corresponding to each sparse grid of $\{\Omega^{l_1,l_2}\}_I$.
However the essential point in sparse-grid computations is that the PDE is never solved directly on the full grid $\Omega^{N_L,N_L}$ in order to save computational costs.
Hence the next key task is to combine solutions on sparse grids to obtain the final solution on the finest grid $\Omega^{N_L,N_L}$.  This is carried out in two steps. The first step is to extend numerical solutions $\{U^{l_1,l_2}\}_I$ on these sparse grids to obtain $(2N_L+1)$ solutions on the finest grid $\Omega^{N_L,N_L}$. This is called ``prolongation''. The second step is to combine all $(2N_L+1)$ sets of numerical solutions on $\Omega^{N_L,N_L}$ to form the final solution $\hat{U}^{N_L,N_L}$.
In the following we first describe the prolongation technique, then provide a summary of the complete algorithm.

\subsubsection{Prolongation and WENO interpolation}

The prolongation procedure is denoted by the operator $P^{N_L,N_L}$.
A prolongation operator $P^{N_L,N_L}$ maps numerical solutions $\{U^{l_1,l_2}\}_I$ on sparse grids onto the finest grid $\Omega^{N_L,N_L}$. For example, with the numerical solution $U^{l_1,l_2}$ on $\Omega^{l_1,l_2}$, the prolongation $P^{N_L,N_L}U^{l_1,l_2}$ generates numerical values for all grid points on the most refined mesh $\Omega^{N_L,N_L}$.
Implementation of prolongation operators is usually done via interpolation procedure. If solutions are smooth, regular Lagrange interpolations are applied directly for prolongations. However, since hyperbolic equations may develop discontinuities in their solutions, we use more robust WENO interpolations in prolongations rather than Lagrange interpolations for a general case. Studies in \cite{GSZ,LKV1,LKV2} for linear schemes and in \cite{LuZhang1,LuChenZhang3} for nonlinear schemes show that the final solution $\hat{U}^{N_L,N_L}$ resulted from the spare-grid combination techniques can achieve the similar accuracy orders as the numerical schemes, as long as the accuracy order of interpolations in the prolongations is not less than the accuracy order of the numerical schemes used to solve PDEs on sparse grids. Hence we use fifth order interpolations here for prolongations.
The interpolations are performed in the dimension by dimension way. We first describe a fifth order WENO interpolation procedure for the one dimensional case.

Given numerical values $u_{i-2}$, $u_{i-1}$, $u_i$, $u_{i+1}$ and $u_{i+2}$ at the
grid points $x_{i-2}$, $x_{i-1}$, $x_i$, $x_{i+1}$ and $x_{i+2}$, we find a fifth order WENO interpolation $u_{WENO}(x)$ for any point $x \in [x_{i-1/2}, x_{i+1/2})$. Here $x_{i-1/2}=(x_{i-1}+x_i)/2$ and $x_{i+1/2}=(x_{i}+x_{i+1})/2$. The grid is uniform with a grid size $h$. The positive linear weights derived as a function of the interpolation location $x$ in \cite{LiuSZ} are used here in the fifth order WENO interpolation.
Here we have the big stencil $S=\{x_{i-2},~x_{i-1},~x_{i},~x_{i+1},~x_{i+2}\}$, and three substencils
$S_0=\{x_{i-2},~x_{i-1},~x_{i}\}$, $S_1=\{x_{i-1},~x_{i},~x_{i+1}\}$, $S_2=\{x_{i},~x_{i+1},~x_{i+2}\}$.
Three polynomials $P_k(x),~k=0,~1,~2,$ of degree at most $2$ are constructed to interpolate $u$ on the three
substencils $S_0,~S_1,~S_2$ respectively. It is found in \cite{LiuSZ} that the combination given by
\begin{equation}
    u_{Lagr}(x)~=~\sum_{k=0}^{2}C_k(x)P_k(x)
\label{lwlag}
\end{equation}
reproduces the Lagrange interpolating polynomial of degree at most $4$ which interpolates $u$ on the big stencil $S$ and is a $5$th order approximation to $u(x)$ if it is smooth on $S$. Here
\begin{align*}
    C_0(x)&=\frac{(x-x_{i+1})(x-x_{i+2})}{12h^2},\\
    C_1(x)&=\frac{(x-x_{i-2})(x-x_{i+2})}{-6h^2},\\
    C_2(x)&=\frac{(x-x_{i-2})(x-x_{i-1})}{12h^2}
\end{align*}
 and they are positive for any $x\in[x_{i-{1}/{2}},x_{i+{1}/{2}})$ that we consider. To find the WENO interpolation $u_{WENO}(x)$, we need to replace the linear weights $C_k(x)$ in (\ref{lwlag}) by the nonlinear weights $w_k(x)$. The nonlinear weights $w_k(x)$ are defined as
\begin{equation}
w_k(x)=\frac{\Tilde{C}_k(x)}{\Tilde{C}_0(x)+\Tilde{C}_1(x)+\Tilde{C}_2(x)}, \qquad
\Tilde{C}_k(x)=\frac{C_k(x)}{(\epsilon+\beta_k)^2}, \qquad k=1, 2, 3,
\label{eqnlnwe}
\end{equation}
where $\epsilon$ is a small
positive number chosen to avoid the denominator becoming $0$.
$\beta_k$ is the the smoothness
indicator for the interpolating polynomial $P_k(x)$ and it is chosen as in the classical fifth order WENO scheme \cite{JS},
\begin{equation}
\beta_k=\sum_{l=1}^2 h^{2l-1} \int_{x_{i-1/2}}^{x_{i+1/2}} \left (\frac{d^l}{dx^l} P_k(x) \right)^2 dx,
\label{si}
\end{equation}
which is actually a quadratic function of numerical values $u_{l}$ on the substencil $S_k$. The final fifth order
WENO interpolation for the point $x$ is
\begin{equation}
    u_{WENO}(x)~=~\sum_{k=0}^{2}w_k(x)P_k(x).
\label{nlweno}
\end{equation}

The interpolations for multi-dimensional cases use the dimension by dimension approach. For example,
in the two dimensional case,
first in every grid line of the $x$ direction with a fixed $y$-coordinate on the sparse grid $\Omega^{l_1,l_2}$, we construct the fifth order WENO interpolations $u_{WENO}(x)$ in (\ref{nlweno}), for $2^{l_1}N_r+1$ different small intervals. Each WENO interpolation uses five adjacent grid points.
Then we evaluate $u_{WENO}(x)$ on the grid points of $\Omega^{N_L,l_2}$, which is the most refined mesh in the $x$ direction. Next, in every grid line of the $y$ direction with a fixed $x$-coordinate on the grid $\Omega^{N_L,l_2}$, we construct the fifth order WENO interpolations $u_{WENO}(y)$ for $2^{l_2}N_r+1$ different small intervals, and evaluate them on the grid points of $\Omega^{N_L,N_L}$. In this way we obtain the fifth order WENO prolongation $P^{N_L,N_L}U^{l_1,l_2}$ on the finest grid $\Omega^{N_L,N_L}$. If in this procedure we replace the fifth order
WENO interpolations by the fifth order Lagrange interpolations as in (\ref{lwlag}), then the fifth order Lagrange prolongation is gotten.

Both the fifth order Lagrange prolongation and WENO prolongation are tested in our numerical experiments in the next section.

\subsubsection{Summary of the algorithm}

We summarize the algorithm of the fifth order sparse grid WENO scheme as following.

%\begin{framed}
%\bigskip
%\bigskip
\newpage
\noindent\textbf{Algorithm: sparse grid WENO5 scheme}
\begin{itemize}
\item Step 1: Restrict the initial condition $u(x,y,0)$ of the equation to these $(2N_L+1)$ sparse grids $\{\Omega^{l_1,l_2}\}_I$ which are aforementioned. Here ``Restrict'' means that functions are
evaluated at grid points;
\item Step 2: On each sparse grid $\Omega^{l_1,l_2}$ of the set $\{\Omega^{l_1,l_2}\}_I$, solve the equation (\ref{eq1}) by the fifth order WENO scheme with the TVD RK time marching to reach the final time $T$. Then we obtain $(2N_L+1)$ sets of solutions $\{U^{l_1,l_2}\}_I$;
\item Step 3: At the final time $T$,
    \begin{itemize}
    \item on each grid $\Omega^{l_1,l_2}$, apply the prolongation operator $P^{N_L,N_L}$ on $U^{l_1,l_2}$. Then we find $P^{N_L,N_L}U^{l_1,l_2}$ on the most refined mesh $\Omega^{N_L,N_L}$. For smooth solutions, the fifth order Lagrange prolongation can be used directly. In general, we use the fifth order WENO prolongation for more robust computations;
    \item carry out the combination to obtain the final solution
    \begin{equation}
    \hat{U}^{N_L,N_L}=\sum_{l_1+l_2=N_L}P^{N_L,N_L}U^{l_1,l_2}-\sum_{l_1+l_2=N_L-1}P^{N_L,N_L}U^{l_1,l_2}.
    \end{equation}
    \end{itemize}
\end{itemize}
%\end{framed}

In general, for higher dimensional problems, the algorithm is similar to the 2D case although prolongation operations are performed in additional spatial directions. The sparse-grid combination formula for a $d$  dimensional problem has the following form (\cite{GSZ}):
\begin{equation}
    \hat{U}^{N_L,\cdots,N_L}~=~\sum_{m=N_L}^{N_L+d-1}(-1)^{d+N_L-(m+1)}\left( \begin{array}{c}{d-1} \\ {m-N_L}\end{array}\right)\sum_{|I_d|=m-(d-1)}P^{N_L,\cdots,N_L}U^{l_1,\cdots,l_d}.
\label{genformusgc}
\end{equation}
Here $N_L$ is the finest level of the sparse grids. $I_d=(l_1,l_2,\cdots,l_d)$ denotes the index of the levels
  of sparse grid $\Omega^{l_1,l_2,\cdots,l_d}$, and $|I_d|=l_1 + l_2 + \cdots + l_d$. $U^{l_1,\cdots,l_d}$ is the numerical solution on
 the sparse grid $\Omega^{l_1,\cdots,l_d}$, and $P^{N_L,\cdots,N_L}$ is the prolongation operator onto the finest  grid $\Omega^{N_L,\cdots,N_L}$. $\hat{U}^{N_L,\cdots,N_L}$ is the final solution on $\Omega^{N_L,\cdots,N_L}$ after the sparse-grid combination.

Specifically in the numerical experiments of this paper, we use the three dimensional (3D) formula
\begin{equation}
\begin{aligned}
&\hat{U}^{N_L,N_L,N_L}=\sum_{l_1+l_2+l_3=N_L}P^{N_L,N_L,N_L}U^{l_1,l_2,l_3}-2\sum_{l_1+l_2+l_3=N_L-1}P^{N_L,N_L,N_L}U^{l_1,l_2,l_3}\\
&+\sum_{l_1+l_2+l_3=N_L-2}P^{N_L,N_L,N_L}U^{l_1,l_2,l_3}\\
\end{aligned}
\end{equation}
and the four dimensional (4D) formula
\begin{equation}
\begin{aligned}
&\hat{U}^{N_L,N_L,N_L,N_L}=\sum_{l_1+l_2+l_3+l_4=N_L}P^{N_L,N_L,N_L,N_L}U^{l_1,l_2,l_3,l_4}\\
&-3\sum_{l_1+l_2+l_3+l_4=N_L-1}P^{N_L,N_L,N_L,N_L}U^{l_1,l_2,l_3,l_4}
+3\sum_{l_1+l_2+l_3+l_4=N_L-2}P^{N_L,N_L,N_L,N_L}U^{l_1,l_2,l_3,l_4}\\
&-\sum_{l_1+l_2+l_3+l_4=N_L-3}P^{N_L,N_L,N_L,N_L}U^{l_1,l_2,l_3,l_4}.\\
\end{aligned}
\end{equation}

%\newpage
\section{Numerical Experiments}
In this section, we use various numerical examples to show
the high computational efficiency of the presented sparse grid fifth order WENO scheme by comparing to the simulations on the corresponding regular grids.
2D, 3D and 4D numerical examples with either smooth or non-smooth solutions are tested.

Analysis on linear schemes for solving linear problems \cite{GSZ, LKV1} shows that the canceling in leading-order error terms of numerical solutions on semi-coarsened grids results in
that the accuracy order of the final solution in the sparse-grid combination technique is kept to be the same as that on a single full grid. By replacing nonlinear weights $w_0$, $w_1$ and $w_2$ in the fifth order WENO approximation (\ref{eq2.7}) with linear weights $d_0$, $d_1$ and $d_2$, the fifth order linear scheme is obtained.
Linear schemes serve as the important base schemes for high order nonlinear WENO schemes.
In the following numerical experiments, we first use the fifth order linear scheme to solve some linear problems with smooth solutions
and verify the linear analysis results in the literature. Due to the high nonlinearity of the fifth order WENO scheme,
it is difficult to perform theoretical analysis to show the accuracy of the sparse grid computations. Hence we use
numerical simulations on nonlinear problems to show that the fifth order sparse grid WENO scheme achieve similar
numerical accuracy as that on the corresponding single grid, which is consistent with the linear cases.
For all examples, we record and compare CPU times for simulations on sparse grids and the corresponding single grids.
In this section, we use $N_h$ to
denote the number of computational cells in one spatial direction of the most refined mesh in sparse grids or the corresponding single grid.

In the numerical examples, the mesh refinement study is performed and numerical convergence rates of the schemes are computed on
successively refined meshes. To refine meshes for simulations on sparse grids, we
refine the root grid (e.g., $\Omega^{0,0,0}$ for the 3D case), but
keep the number of semi-coarsened sparse-grid levels (total $N_L+1$ levels) unchanged. For example, 3D sparse grids with a $10 \times 10\times 10$ root grid and
$N_L=3$ have the finest mesh $80 \times 80\times 80$. When the root grid is refined once to be $20 \times 20\times 20$ and $N_L=3$ is kept unchanged,
the finest mesh $160 \times 160\times 160$ is obtained then.

As that discovered in \cite{LuChenZhang3}, time step sizes used to march the PDEs on all semi-coarsened sparse grids have to be determined by the spatial grid size $\triangle x$ of the most refined grid, i.e., $\Omega^{N_L,N_L}$ in the 2D problems,
or $\Omega^{N_L,N_L,N_L}$ in the 3D problems. Note that $\triangle x$ is the minimum one among the grid sizes of all spatial directions in the 2D grid $\Omega^{N_L,N_L}$, or the 3D grid $\Omega^{N_L,N_L,N_L}$, and so on.
At each time step, the same time step size are used for each individual time evolution on different semi-coarsened sparse grids. It is determined by the spatial grid size of the most refined grid (e.g., the $\Omega^{N_L,N_L}$ in the 2D problems), the chosen CFL number and the current wave speed of the PDEs.
Note that the time step sizes could be different at different time steps due to the change of wave speeds in the nonlinear problems.
The desired numerical accuracy in sparse-grid computations is obtained if the time step sizes are taken in this way
for a general problem.
We follow this method to choose time step sizes in this paper.
The third order TVD-RK scheme is used for the time evolution, so
in the examples with a smooth solution which are solved for accuracy order tests of the schemes, we follow the common practice in the literature and take the time step size $\triangle t = (\triangle x)^{\frac{5}{3}}$ in order to observe the fifth order convergence rate, while the usual CFL condition is applied for other examples (i.e., $\triangle t \sim O(\triangle x)$). All of the numerical simulations in this paper are performed on a 2.3 GHz, 16GB RAM Linux workstation.

\subsection{Linear cases}

High order nonlinear WENO schemes are built on the corresponding high order linear schemes. We first test the fifth order linear scheme by solving linear advection problems on sparse grids and show consistent results with the theoretical linear analysis in the literature, e.g. \cite{GSZ, LKV1}.

\begin{figure}[p]
    \centering
    \begin{subfigure}[b]{0.48\textwidth}
        \includegraphics[width=0.9\textwidth]{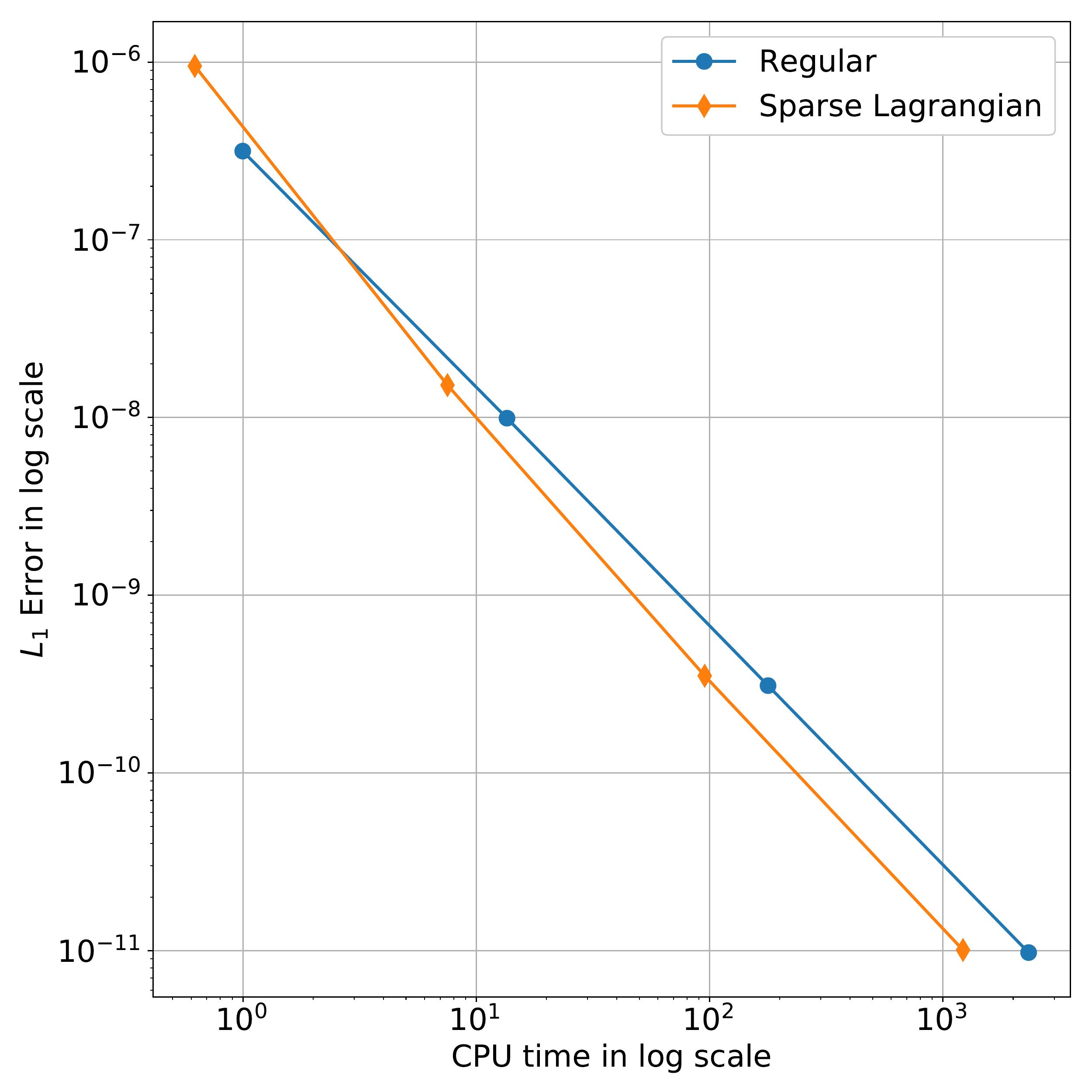}
        \caption{$L^1$ error and CPU time dependence}
        \label{2d_ln_l1}
    \end{subfigure}
        \begin{subfigure}[b]{0.48\textwidth}
        \includegraphics[width=0.9\textwidth]{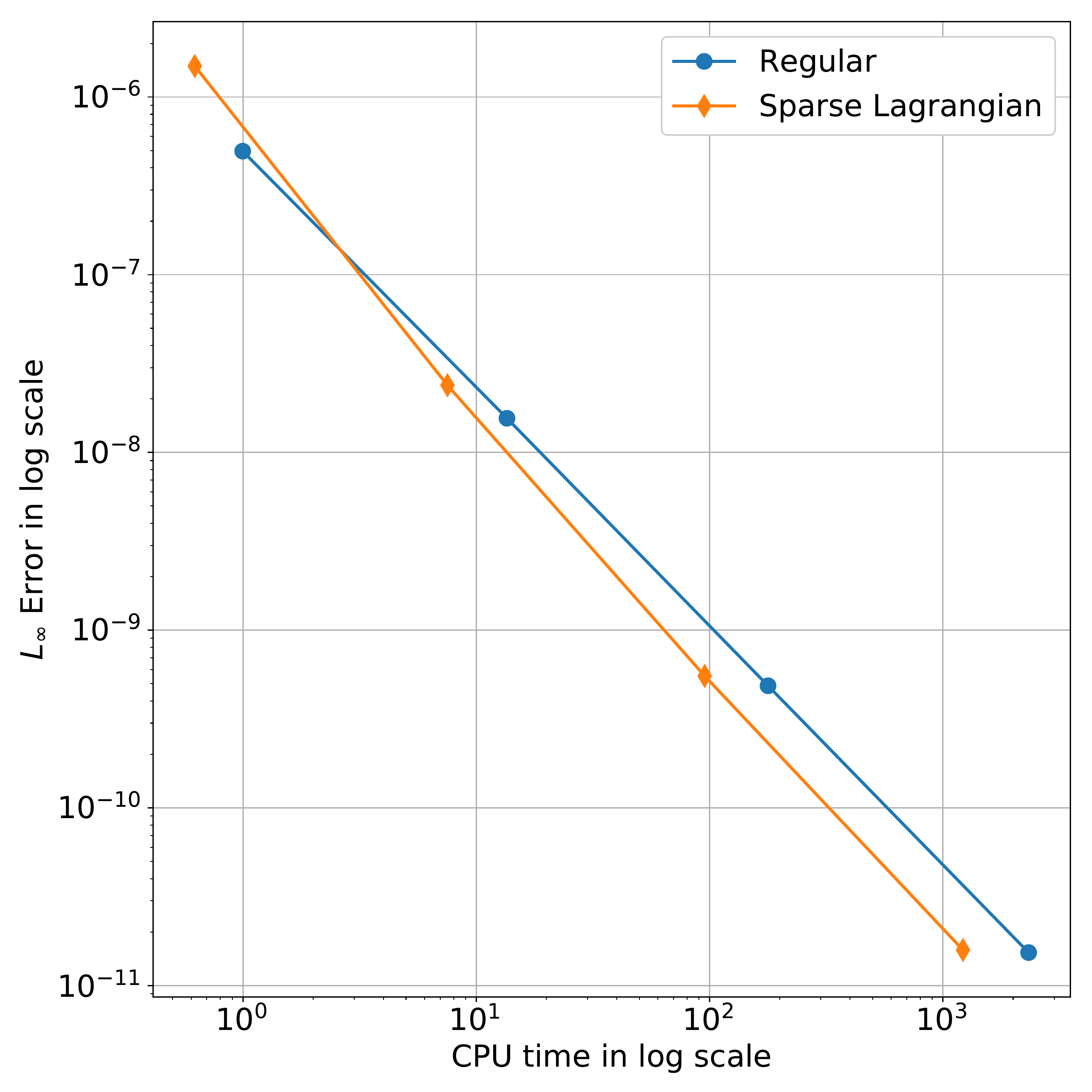}
        \caption{$L^\infty$ error and CPU time dependence}
        \label{2d_ln_linf}
    \end{subfigure}
    \caption{Example 1, a 2D linear advection equation is solved by the linear scheme. log-log plots of the dependence of numerical errors and CPU times for computations on single-grid and sparse-grid. Orange lines with diamonds: sparse-grid computations; blue lines with solid circles: single-grid computations.}
    \label{2dplot_linear}
\end{figure}

\paragraph{Example 1 (A 2D linear advection equation).} In this example, we solve a 2D linear advection equation
\begin{equation*}
    u_t~+~u_x~+~u_y~=~0,
\end{equation*}
with the initial condition
\begin{equation*}
    u_0(x,~y)~=~0.3+0.7\sin(\frac{\pi}{2}(x+y))
\end{equation*}
and periodic boundary conditions. The computational domain is $[0,~4]\times[0,~4]$.
The fifth order linear scheme on both single grids and sparse grids are used to compute the numerical solution of the problem at the final time $T=0.5$. The time step size $\triangle t = (\triangle x)^{\frac{5}{3}}$.
The efficiency of single-grid computation and sparse-grid computation is compared.
We report the $L^1$ errors, $L^\infty$ errors, the corresponding numerical accuracy orders, and CPU times in Table \ref{tab:ex1}. It is observed that both single-grid computations and sparse-grid computations achieve the fifth order accuracy, and their numerical errors are comparable, especially on refined meshes. However, computations on sparse grids are more efficient than those on single grids. From the recorded CPU times in Table \ref{tab:ex1}, on relatively refined $N_h\times N_h$ meshes  we see that close to
$50\%$ computation time is saved for computations on sparse grids
to reach the similar error levels as that on single grids. This is further indicated in Figure \ref{2dplot_linear}, in which log-log plots of the numerical errors as functions of CPU times are shown for computations on single grids and sparse grids. It can be clearly seen that less CPU time is needed for sparse-grid computation than
for single-grid computation to reach a small numerical error.

\begin{table}[htbp]\footnotesize
\centering
\begin{tabular}{c c c c c c c c}
\hline
\multicolumn{8}{c}{ Single-grid}\\ \hline
& & $N_h\times N_h$		&$L^1$ error		&Order  &$L^\infty$ error    &Order  &CPU(s)			\\
\hline
&&$80\times 80$         &3.1556e-07 &        - &  4.9572e-07 &          - &     0.998    	\\
&&$160\times 160$       &9.9122e-09 &    4.993 &  1.5571e-08 &      4.993 &    13.544	\\
&&$320\times 320$       &3.0987e-10 &    4.999 &  4.8676e-10 &          5.000 &   178.068	\\
&&$640\times 640$       &9.7480e-12 &     4.990 &  1.5350e-11 &      4.987 &  2329.770	\\
\hline
\multicolumn{8}{c}{ Sparse-grid}\\ \hline
$N_r$   &$N_L$  &$N_h\times N_h$	&$L^1$ error   &Order	    &$L^\infty$ error    &Order  &CPU(s)	     \\
\hline
10  &3  &$80\times 80$         &9.5053e-07 &        - &  1.4942e-06 &          - &     0.622	\\
20  &3  &$160\times 160$       &1.5212e-08 &    5.965 &  2.3898e-08 &      5.966 &     7.528 	\\
40  &3  &$320\times 320$       &3.5197e-10 &    5.434 &  5.5290e-10 &      5.434 &    95.318  	\\
80  &3  &$640\times 640$       &1.0078e-11 &    5.126 &  1.5873e-11 &      5.122 &  1219.498 	\\
\hline
\end{tabular}
\caption{\footnotesize{Example 1, a 2D linear advection equation is solved by the linear scheme.
Comparison of numerical errors and CPU times for computations on single-grid and sparse-grid. Lagrange interpolation for prolongation is used in sparse-grid computations.
Final time $T=0.5$.
$N_r$: number of cells in each spatial direction of a root grid.
$N_L$: the finest level in a sparse-grid computation.
CPU: CPU time for a complete simulation. CPU time unit: seconds.}}
\label{tab:ex1}
\end{table}

\begin{figure}[p]
    \centering
    \begin{subfigure}[b]{0.48\textwidth}
        \includegraphics[width=0.9\textwidth]{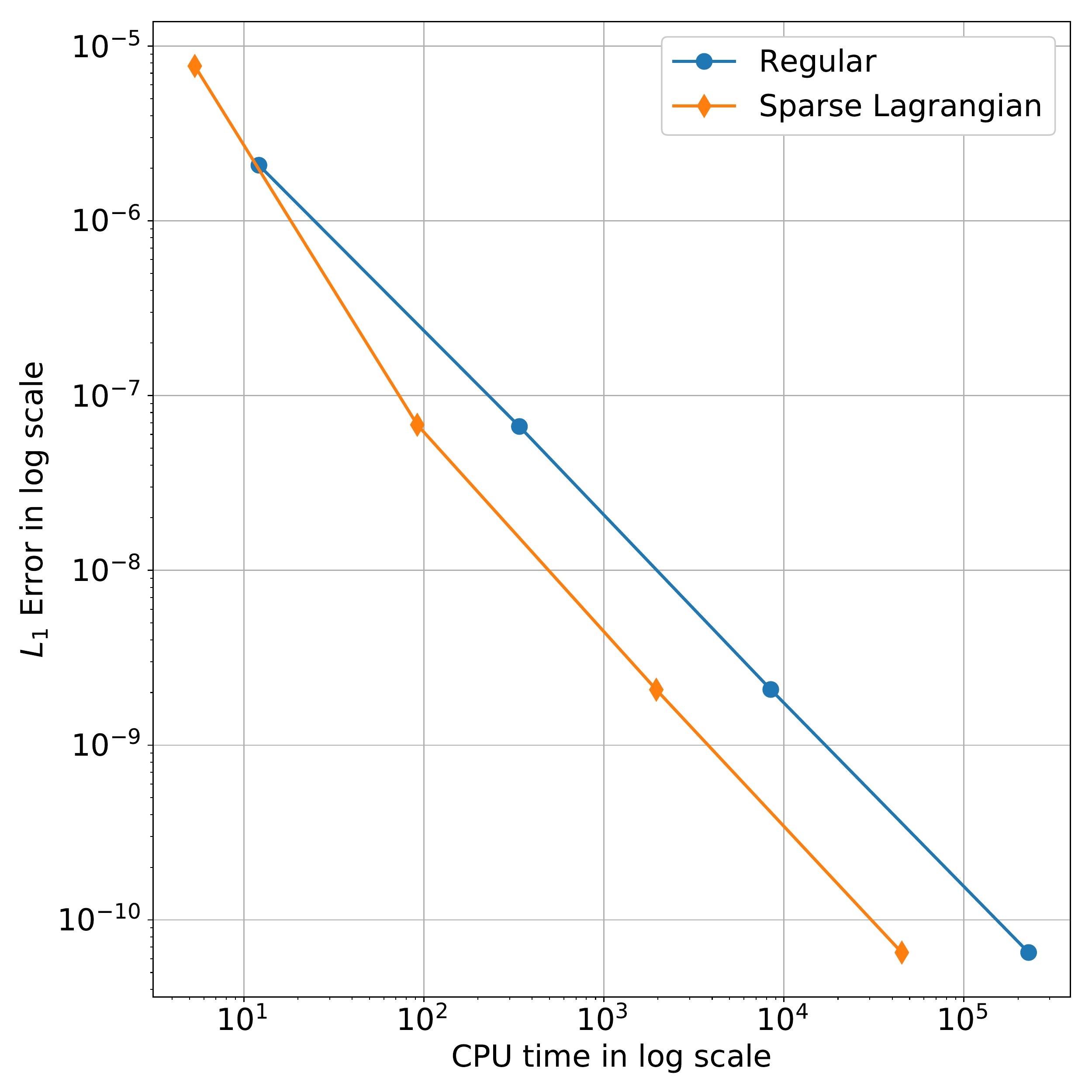}
        \caption{$L^1$ error and CPU time dependence}
        \label{3d_ln_l1}
    \end{subfigure}
        \begin{subfigure}[b]{0.48\textwidth}
        \includegraphics[width=0.9\textwidth]{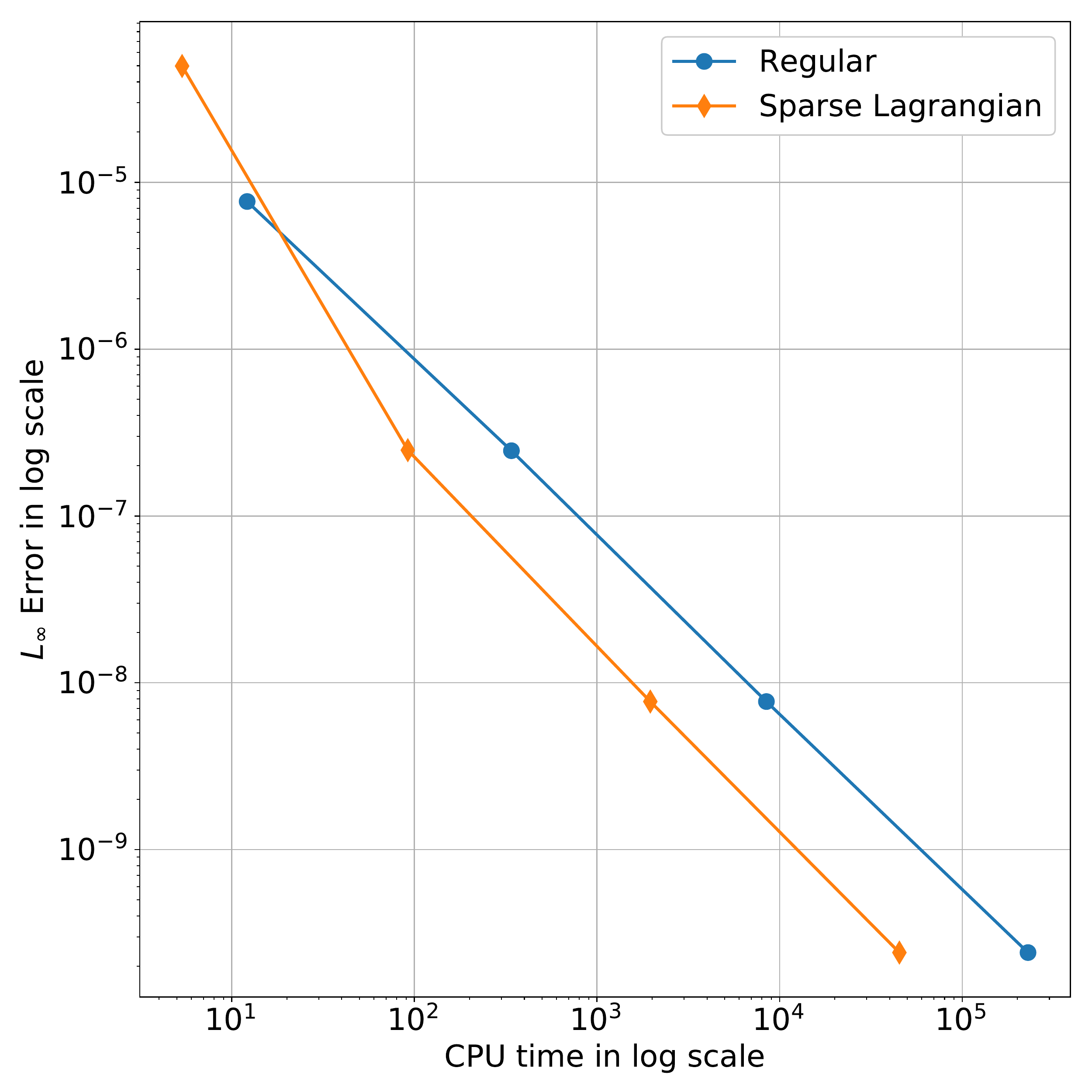}
        \caption{$L^\infty$ error and CPU time dependence}
        \label{3d_ln_linf}
    \end{subfigure}
    \caption{Example 2, a 3D linear advection equation is solved by the linear scheme. log-log plots of the dependence of numerical errors and CPU times for computations on single-grid and sparse-grid. Orange lines with diamonds: sparse-grid computations; blue lines with solid circles: single-grid computations.}
    \label{3dplot_linear}
\end{figure}

\paragraph{Example 2 (A 3D linear advection equation).} In this example, we go one dimension higher than the last example and solve a 3D linear advection equation
\begin{equation*}
    u_t~+~u_x~+~u_y~+~u_z~=~0,
\end{equation*}
with the initial condition
\begin{equation*}
    u_0(x,~y)~=~1.0~+~0.5\sin(x~+~y~+~z)
\end{equation*}
and periodic boundary conditions. The computational domain is $[0,~2\pi]\times[0,~2\pi]\times[0,~2\pi]$.
We use the fifth order linear scheme on both single grids and sparse grids to compute the numerical solution of the problem at the final time $T=0.5$, and compare the efficiency of single-grid computation and sparse-grid computation. The time step size $\triangle t = (\triangle x)^{\frac{5}{3}}$. In Table \ref{tab:ex2}, the $L^1$ errors, $L^\infty$ errors, the corresponding numerical accuracy orders, and CPU times are reported. Similar as the last example, both single-grid computations and sparse-grid computations achieve the fifth order accuracy, and their numerical errors are comparable on refined meshes, except that on the relatively coarse mesh the numerical errors on sparse grids are larger than that on the corresponding single grid. Overall the computations on sparse grids are much more efficient than those on single grids. From the recorded CPU times in Table \ref{tab:ex2}, on relatively refined $N_h\times N_h\times N_h$ meshes we see that
$70\% \sim 80\%$ computation time is saved for computations on sparse grids
to reach the similar error levels as that on the corresponding single grids. More CPU time savings are achieved in this 3D problem than the last 2D example. In Figure \ref{3dplot_linear}, log-log plots of the numerical errors as functions of CPU times are presented for computations on single grids and sparse grids, which shows that the sparse-grid computation is much more efficient than the single-grid computation to reach a small numerical error.

\begin{table}[htbp]\footnotesize
\centering
\begin{tabular}{c c c c c c c c}
\hline
\multicolumn{8}{c}{ Single-grid}\\ \hline
& & $N_h\times N_h\times N_h$		&$L^1$ error		&Order  &$L^\infty$ error    &Order  &CPU(s)			 \\
\hline
&&80 $\times$ 80 $\times$ 80& 2.0810e-06 &        - &  7.6743e-06 &          - &      12.132 \\
&&160 $\times$ 160 $\times$ 160& 6.6581e-08 &    4.966 &  2.4607e-07 &      4.963 &     339.748 \\
&&320 $\times$ 320 $\times$ 320& 2.0825e-09 &    4.999 &  7.7131e-09 &      4.996 &    8461.060 \\
&&640 $\times$ 640 $\times$ 640& 6.5097e-11 &    5.000 &  2.4123e-10 &      4.999 &  231386.000 \\
\hline
\multicolumn{8}{c}{ Sparse-grid}\\ \hline
$N_r$   &$N_L$  &$N_h\times N_h\times N_h$	&$L^1$ error   &Order	    &$L^\infty$ error    &Order  &CPU(s)	     \\
\hline
10 & 3 &   80 $\times$ 80 $\times$ 80&  7.6908e-06 &        - &  4.9805e-05 &          - &      5.342 \\
20 & 3 & 160 $\times$ 160 $\times$ 160& 6.8082e-08 &     6.820 &  2.4792e-07 &       7.650 &     91.994 \\
40 & 3 & 320 $\times$ 320 $\times$ 320&  2.0761e-09 &    5.035 &  7.7136e-09 &      5.006 &   1957.354 \\
80 & 3 & 640 $\times$ 640 $\times$ 640&  6.5039e-11 &    4.996 &  2.4124e-10 &      4.999 &  45248.280 \\
\hline
\end{tabular}
\caption{\footnotesize{Example 2, a 3D linear advection equation is solved by the linear scheme.
Comparison of numerical errors and CPU times for computations on single-grid and sparse-grid. Lagrange interpolation for prolongation is used in sparse-grid computations.
Final time $T=0.5$.
$N_r$: number of cells in each spatial direction of a root grid.
$N_L$: the finest level in a sparse-grid computation.
CPU: CPU time for a complete simulation. CPU time unit: seconds.}}
\label{tab:ex2}
\end{table}

\subsection{Nonlinear cases with smooth solutions}

Now we consider the nonlinear cases, and apply the sparse grid WENO method to nonlinear problems. First the problems with smooth solutions
are tested. We use 2D and 3D nonlinear burgers' equations as examples.

\bigskip
\noindent{\bf Example 3 (Burgers' equations with smooth solutions).}

\begin{figure}[p]
    \centering
    \begin{subfigure}[b]{0.48\textwidth}
        \includegraphics[width=0.9\textwidth]{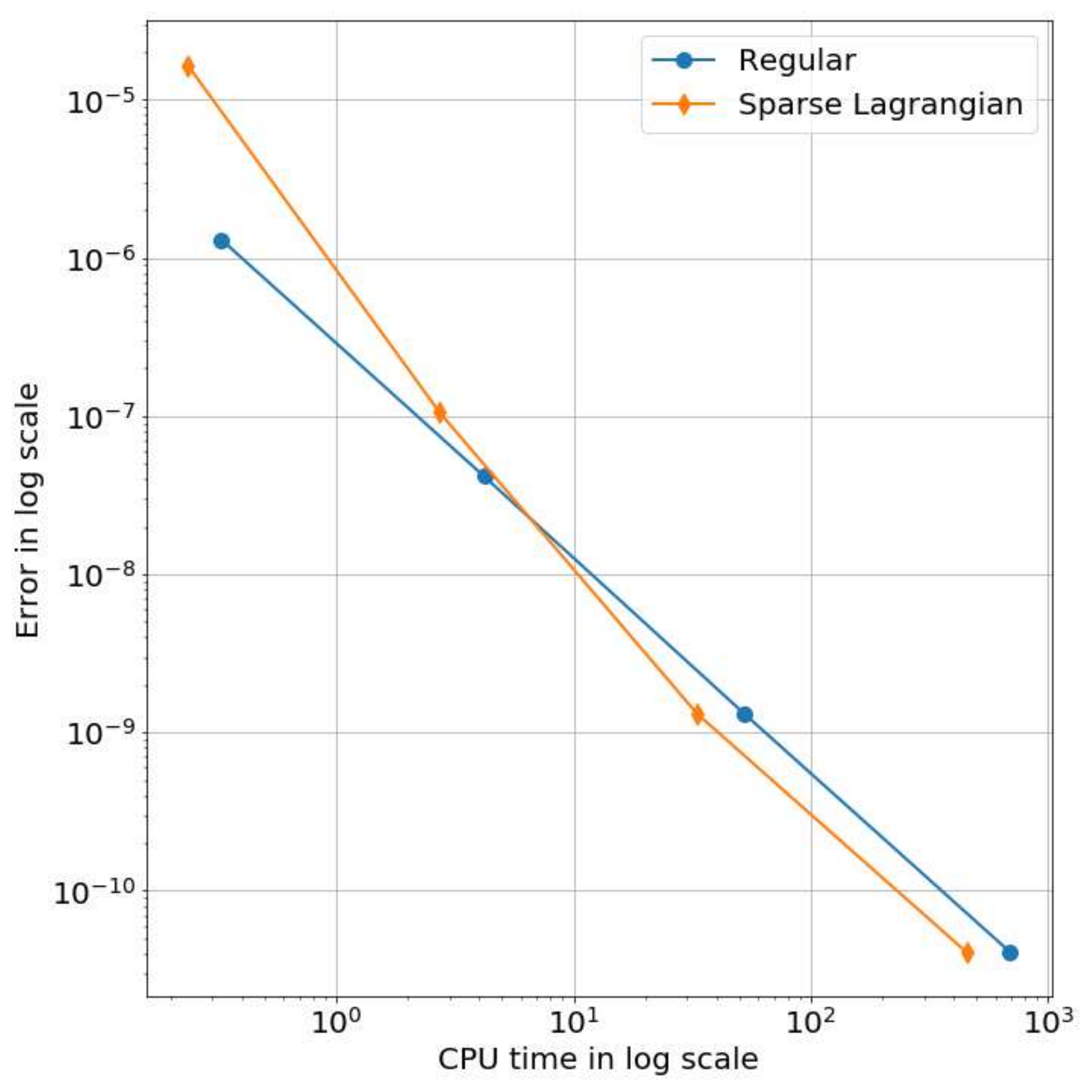}
        \caption{\footnotesize{Linear scheme, $L^1$ error vs. CPU time}}
        \label{2d_burg_ln_l1}
    \end{subfigure}
        \begin{subfigure}[b]{0.48\textwidth}
        \includegraphics[width=0.9\textwidth]{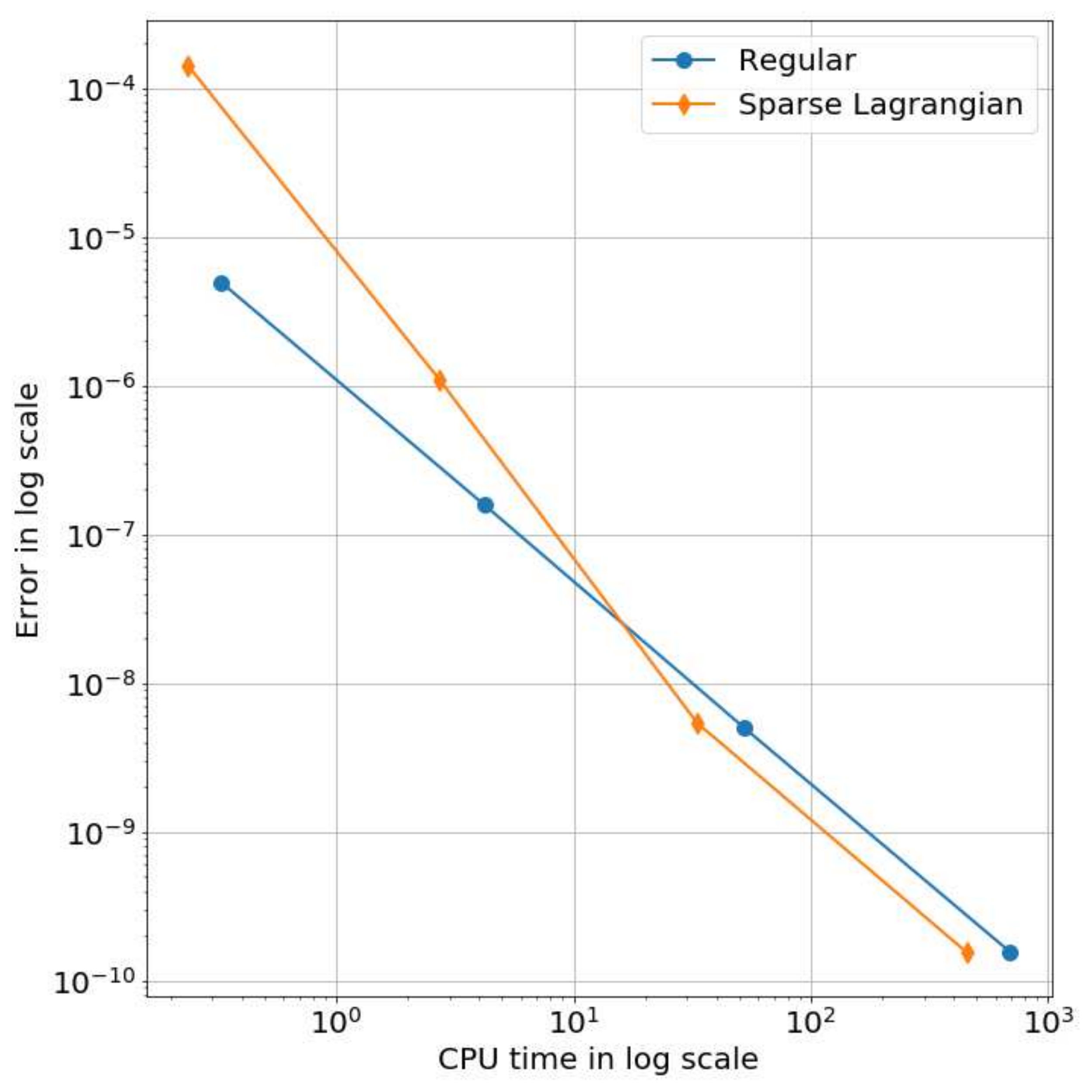}
        \caption{\footnotesize{Linear scheme, $L^\infty$ error vs. CPU time}}
        \label{2d_burg_ln_linf}
    \end{subfigure}
    \begin{subfigure}[b]{0.48\textwidth}
        \includegraphics[width=0.9\textwidth]{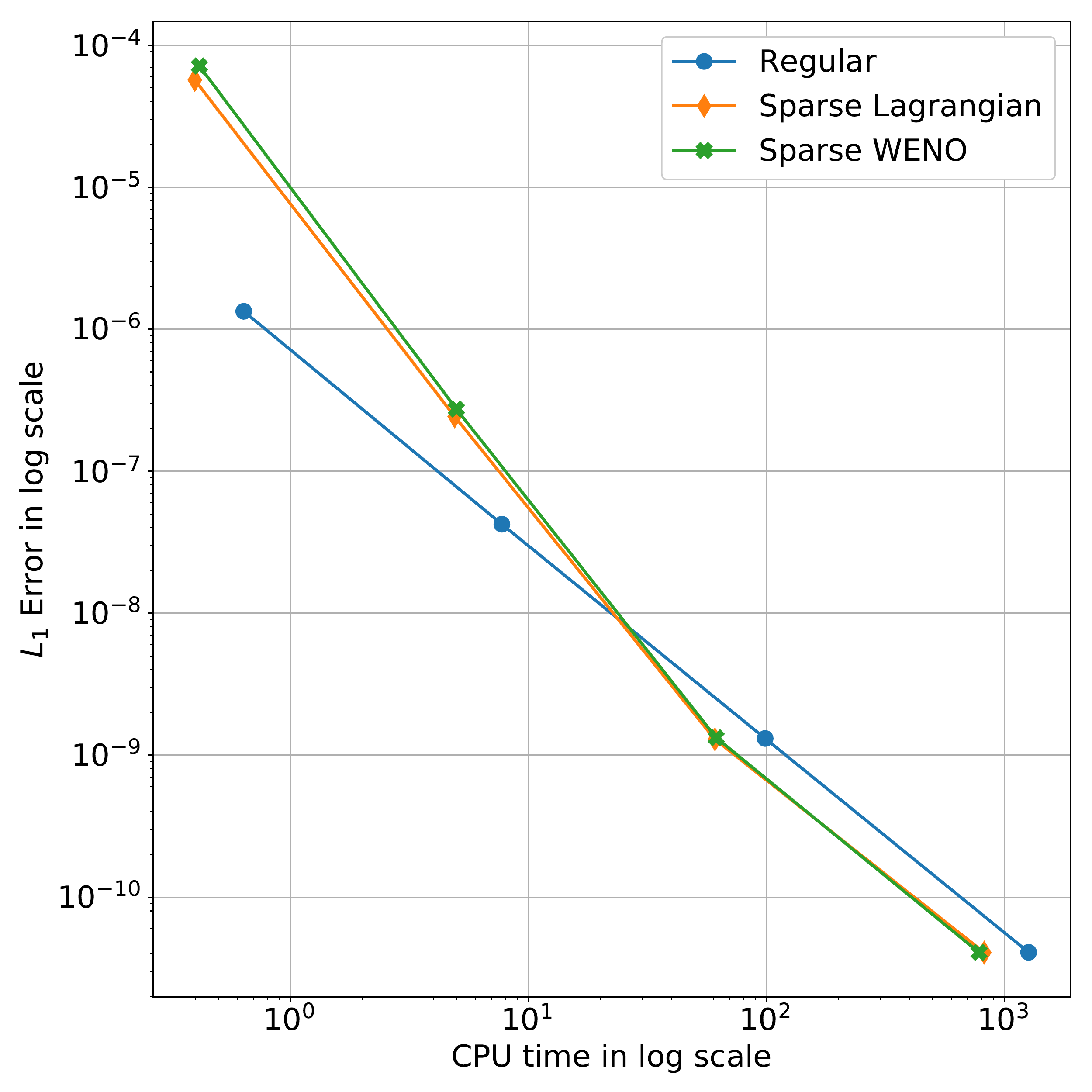}
        \caption{\footnotesize{WENO scheme, $L^1$ error vs. CPU time}}
        \label{2d_burg_weno_l1}
    \end{subfigure}
        \begin{subfigure}[b]{0.48\textwidth}
        \includegraphics[width=0.9\textwidth]{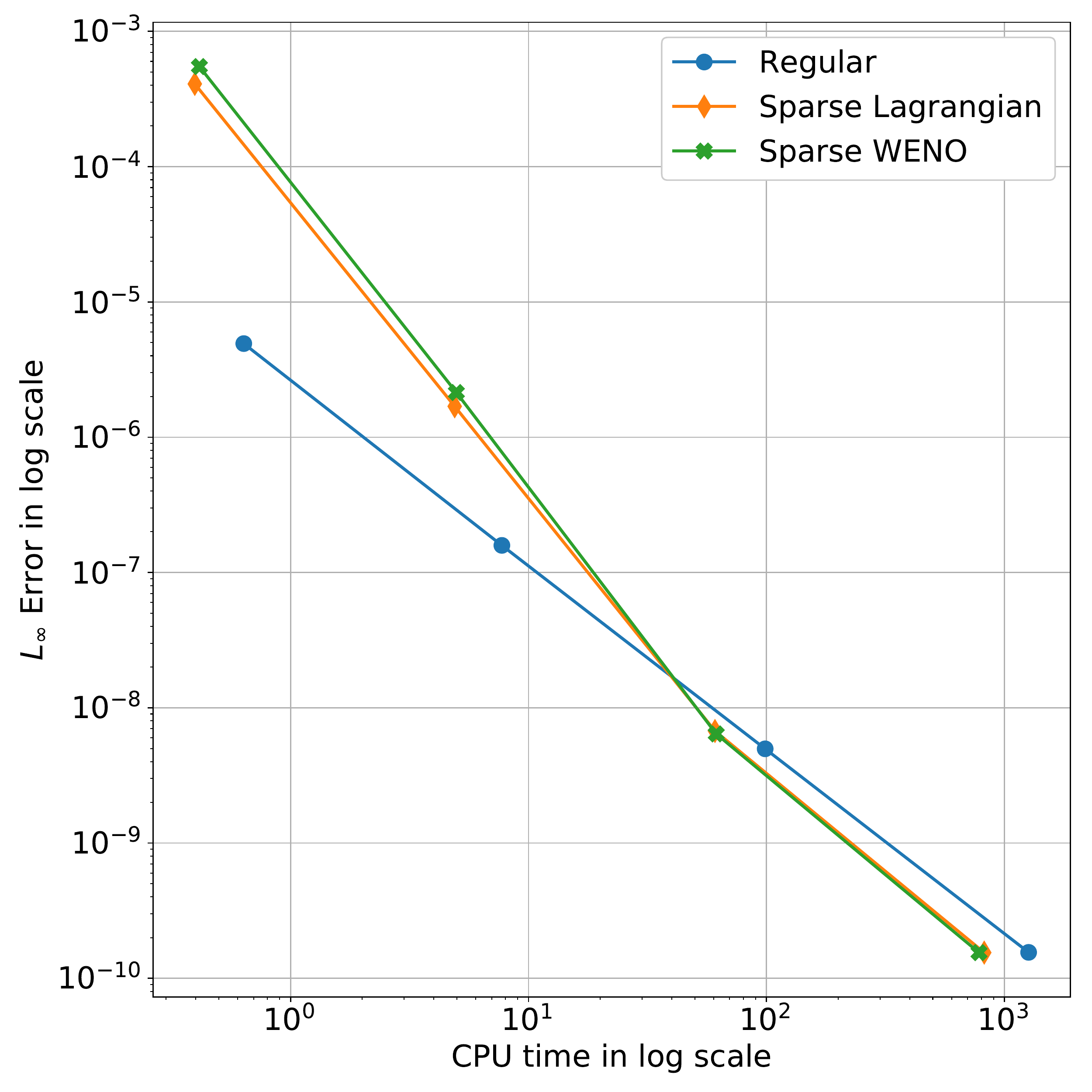}
        \caption{\footnotesize{WENO scheme, $L^\infty$ error vs. CPU time}}
        \label{2d_burg_weno_linf}
    \end{subfigure}
    \caption{\footnotesize{Example 3(a), a 2D Burgers' equation is solved by the linear scheme and the WENO scheme. log-log plots of the dependence of numerical errors and CPU times for computations on single-grid and sparse-grid. Green lines with crosses: sparse-grid computations with WENO prolongation; orange lines with diamonds: sparse-grid computations with Lagrange prolongation; blue lines with solid circles: single-grid computations.}}
    \label{2dplot_burgeqn}
\end{figure}

\noindent{\bf(a)} The equation of 2D case is
\begin{equation*}
    u_t~+~\Big(\frac{u^2}{2}\Big)_x~+~\Big(\frac{u^2}{2}\Big)_y~=~0,
\end{equation*}
with the initial condition
\begin{equation*}
    u_0(x, y)~=~1~+~\frac{1}{2}\sin(x+y)
\end{equation*}
and periodic boundary conditions.
The computational domain is $[0,~2\pi]~\times~[0,~2\pi]$.
We use both the fifth order linear scheme and WENO scheme on sparse grids and the corresponding single grids to compute the numerical solution of the problem at the final time $T=0.3$, when the solution is still smooth. Both the Lagrange interpolation and the WENO interpolation for prolongation are respectively used in sparse-grid computations of the WENO scheme.
The time step size $\triangle t = (\triangle x)^{\frac{5}{3}}$. We report the $L^1$ errors, $L^\infty$ errors, the corresponding numerical accuracy orders, and CPU times in Table \ref{tab:ex3alin} for the linear scheme and in Table \ref{tab:ex3aweno} for the WENO scheme. The fifth order accuracy is obtained for all cases along with the mesh refinement. The numerical errors and computational efficiency of the sparse-grid WENO computations by using the Lagrange interpolation and the WENO interpolation for prolongation are comparable.
However, it is more obvious than the previous linear examples that on the relatively coarse meshes the numerical errors on sparse grids are larger than that on the corresponding single grids for these nonlinear cases, including the case that the linear scheme is applied to the nonlinear equation. With more refined meshes, sparse-grid computations show a behavior of superconvergence and achieve comparable numerical errors as single-grid computations. Hence for this nonlinear problem, the sparse-grid computations are more efficient than the corresponding single-grid computations
on refined meshes to reach similar numerical errors. From the recorded CPU times, we see that more than $30\%$ CPU time is saved by computations on sparse grids for this 2D problem. This is also illustrated in Figure \ref{2dplot_burgeqn}.

\begin{figure}[p]
    \centering
    \begin{subfigure}[b]{0.48\textwidth}
        \includegraphics[width=0.9\textwidth]{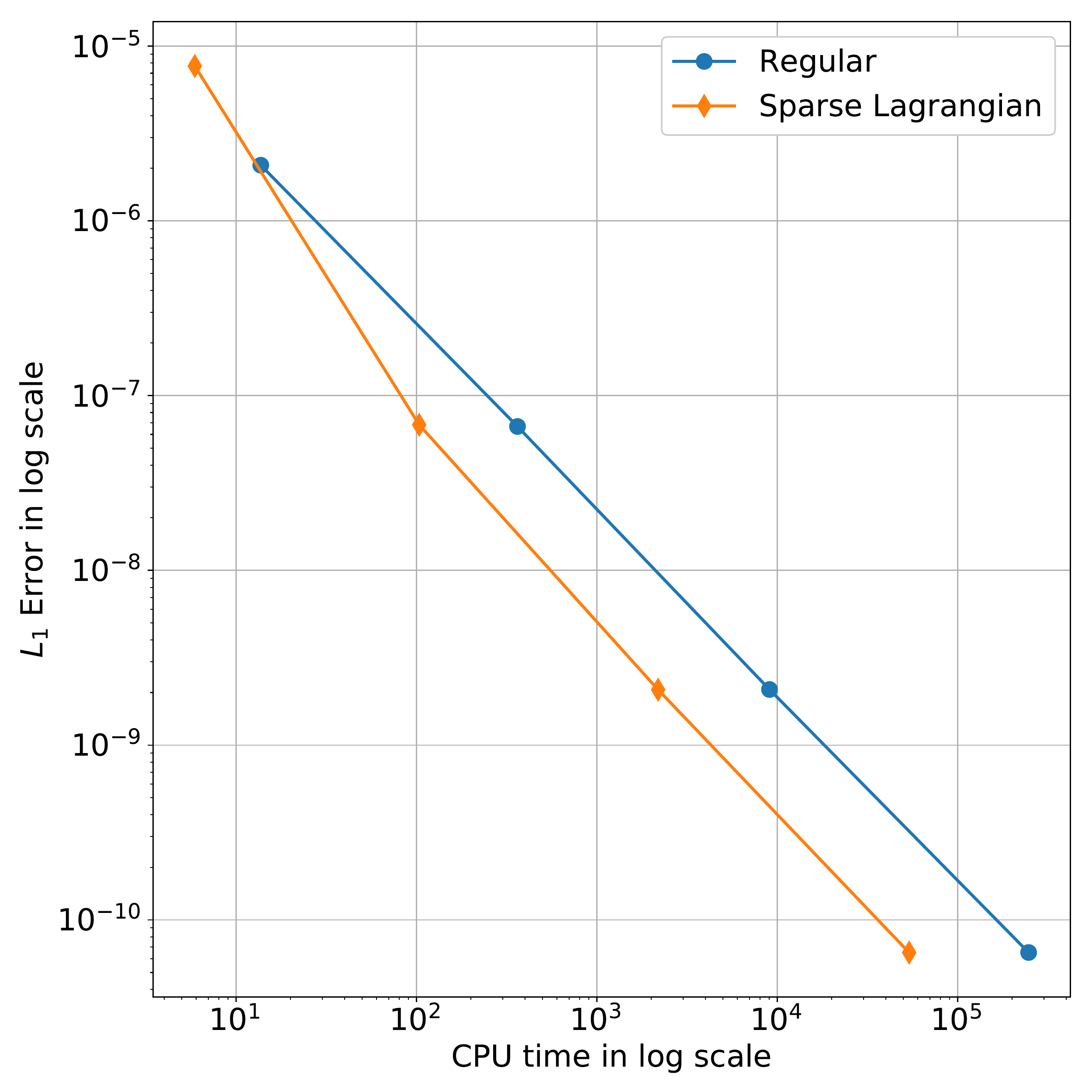}
        \caption{\footnotesize{Linear scheme, $L^1$ error vs. CPU time}}
        \label{3d_burg_ln_l1}
    \end{subfigure}
        \begin{subfigure}[b]{0.48\textwidth}
        \includegraphics[width=0.9\textwidth]{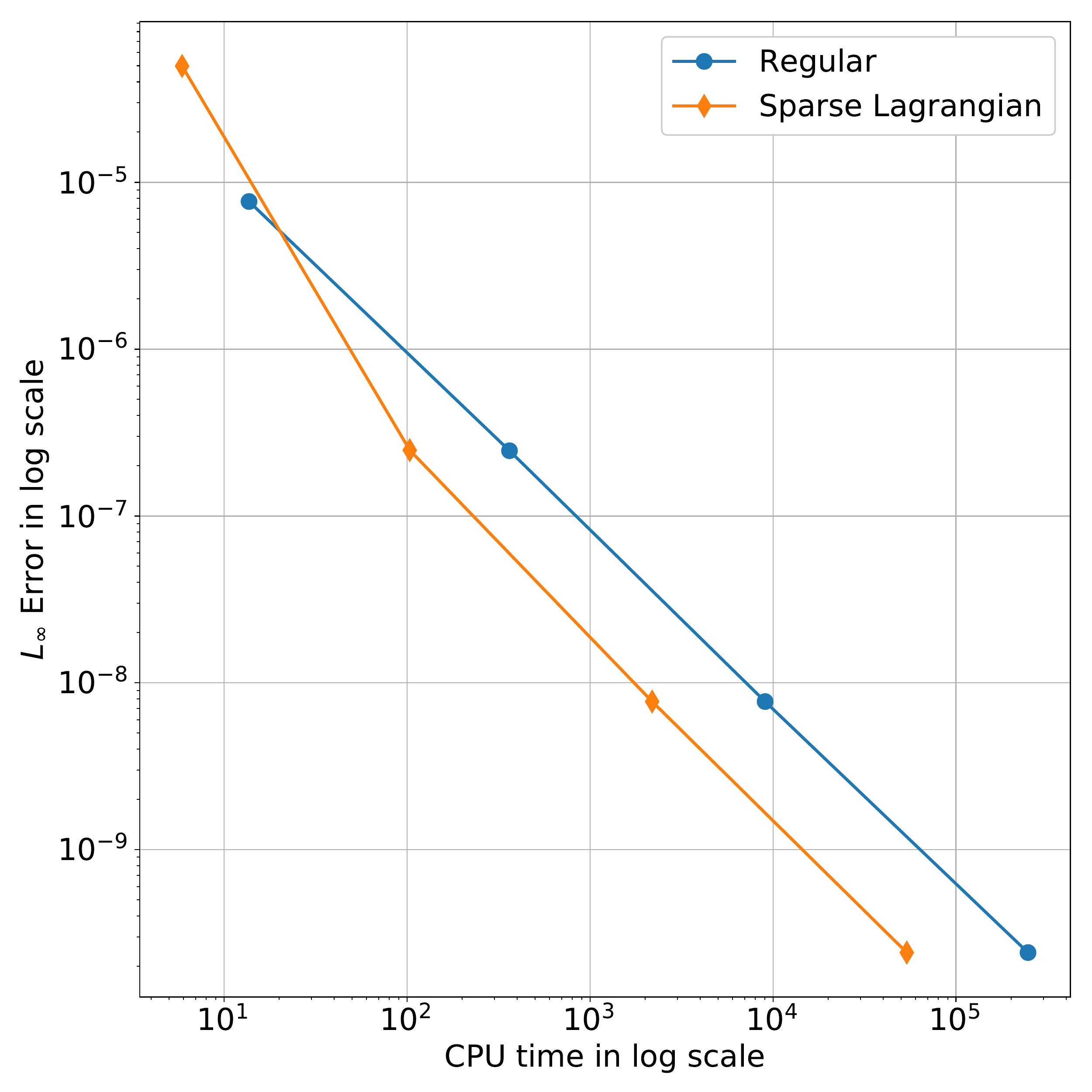}
        \caption{\footnotesize{Linear scheme, $L^\infty$ error vs. CPU time}}
        \label{3d_burg_ln_linf}
    \end{subfigure}
    \begin{subfigure}[b]{0.48\textwidth}
        \includegraphics[width=0.9\textwidth]{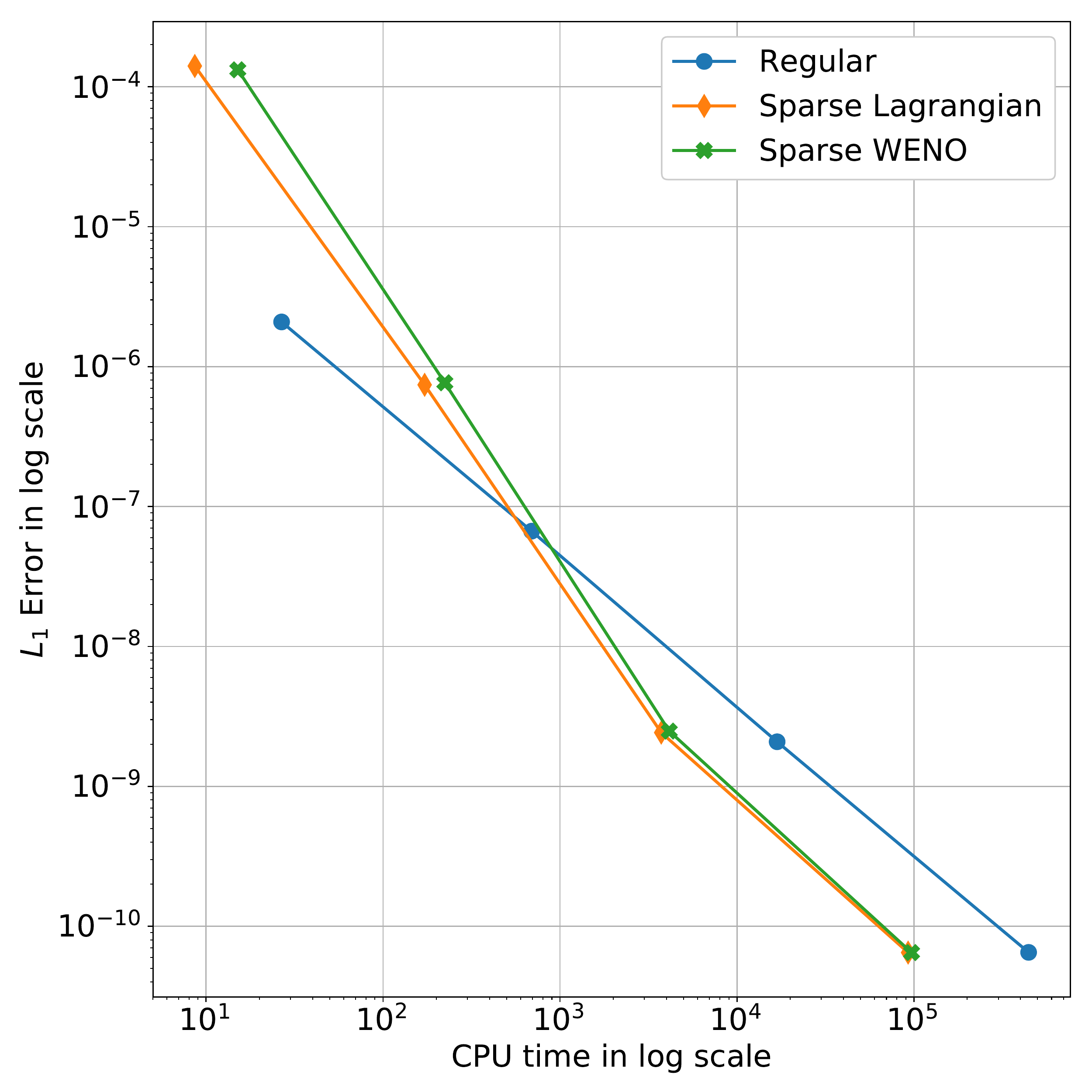}
        \caption{\footnotesize{WENO scheme, $L^1$ error vs. CPU time}}
        \label{3d_burg_weno_l1}
    \end{subfigure}
        \begin{subfigure}[b]{0.48\textwidth}
        \includegraphics[width=0.9\textwidth]{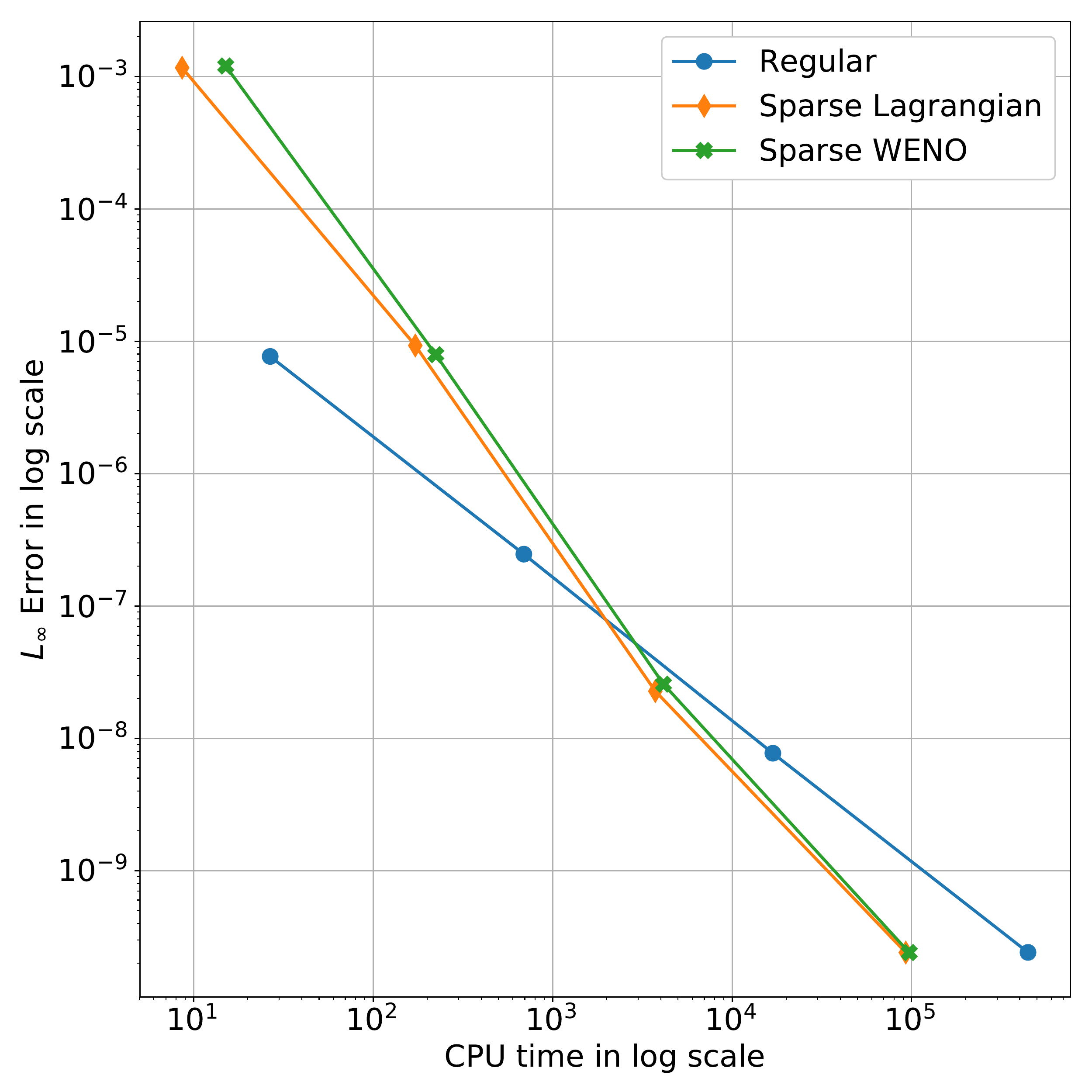}
        \caption{\footnotesize{WENO scheme, $L^\infty$ error vs. CPU time}}
        \label{3d_burg_weno_linf}
    \end{subfigure}
    \caption{\footnotesize{Example 3(b), a 3D Burgers' equation is solved by the linear scheme and the WENO scheme. log-log plots of the dependence of numerical errors and CPU times for computations on single-grid and sparse-grid. Green lines with crosses: sparse-grid computations with WENO prolongation; orange lines with diamonds: sparse-grid computations with Lagrange prolongation; blue lines with solid circles: single-grid computations.}}
    \label{3dplot_burgeqn}
\end{figure}

\bigskip
\noindent{\bf(b)}
The equation of 3D case is
\begin{equation*}
    u_t~+~\Big(\frac{u^2}{2}\Big)_x~+~\Big(\frac{u^2}{2}\Big)_y~+~\Big(\frac{u^2}{2}\Big)_z~=~0,
\end{equation*}
with the initial condition
\begin{equation*}
    u_0(x, y, z)~=~1~+~\frac{1}{2}\sin(x+y+z)
\end{equation*}
and periodic boundary conditions.
The computational domain is $[0,~2\pi]\times[0,~2\pi]\times[0,~2\pi]$.
Similar as the 2D case,
we use both the fifth order linear scheme and WENO scheme on sparse grids and the corresponding single grids to compute the numerical solution of the problem at the final time $T=0.1$, when the solution is still smooth. We use both the Lagrange interpolation and the WENO interpolation for prolongation in sparse-grid computations of the WENO scheme, and
compare their accuracy. Again the time step size is taken to be $\triangle t = (\triangle x)^{\frac{5}{3}}$ in this example for the accuracy order test. The $L^1$ errors, $L^\infty$ errors, the corresponding numerical accuracy orders, and CPU times are reported in Table \ref{tab:ex3-balin} for the linear scheme  and in Table \ref{tab:ex3-baweno} for the WENO scheme. As the 2D case, we see that the desired accuracy orders are obtained for all cases when the meshes are refined. Comparing the results of using the Lagrange interpolation and the WENO interpolation for prolongation in the sparse-grid WENO computations in Table \ref{tab:ex3-baweno}, we observe comparable numerical errors while the WENO interpolation for prolongation takes slightly more CPU time than the Lagrange interpolation. Also similar as the 2D case, on the relatively coarse meshes the numerical errors on sparse grids are larger than that on the corresponding single grids for these nonlinear cases. However, this issue is resolved on refined meshes and sparse-grid computations achieve comparable numerical errors as their corresponding single-grid computations. Especially much more CPU times are
saved in the 3D sparse-grid computations than the 2D case. From the recorded CPU times, we see that $70\% \sim 80\%$ CPU time is saved by computations on sparse grids for this 3D problem. This is also shown in Figure \ref{3dplot_burgeqn}.

\begin{table}[htbp]\footnotesize
\centering
\begin{tabular}{c c c c c c c c}
\hline
\multicolumn{8}{c}{ Single-grid}\\ \hline
& & $N_h\times N_h$		&$L^1$ error		&Order  &$L^\infty$ error    &Order  &CPU(s)			 \\
\hline
&&80 $\times$ 80 & 1.3025e-06 &        - &  4.9233e-06 &          - &    0.308 \\
&&160 $\times$ 160 & 4.1639e-08 &    4.967 &  1.5854e-07 &      4.957 &    3.788 \\
&&320 $\times$ 320 & 1.3048e-09 &    4.996 &  4.9715e-09 &      4.995 &   51.144 \\
&&640 $\times$ 640 & 4.0818e-11 &    4.998 &  1.5559e-10 &      4.998 &  659.734 \\
\hline
\multicolumn{8}{c}{ Sparse-grid}\\ \hline
$N_r$   &$N_L$  &$N_h\times N_h$	&$L^1$ error   &Order	    &$L^\infty$ error    &Order  &CPU(s)	     \\
\hline
10 & 3 &   80 $\times$ 80 & 1.6557e-05 &        - &  1.4279e-04 &          - &    0.236 \\
20 & 3 & 160 $\times$ 160 & 1.0572e-07 &    7.291 &  1.0969e-06 &      7.024 &    2.720 \\
40 & 3 & 320 $\times$ 320 & 1.3135e-09 &    6.331 &  5.4073e-09 &      7.664 &   33.156 \\
80 & 3 & 640 $\times$ 640 & 4.0759e-11 &     5.010 &  1.5506e-10 &      5.124 &  454.298 \\
\hline
\end{tabular}
\caption{\footnotesize{Example 3(a), a 2D Burgers' equation is solved by the linear scheme.
Comparison of numerical errors and CPU times for computations on single-grid and sparse-grid. Lagrange interpolation for prolongation is used in sparse-grid computations.
Final time $T=0.3$.
$N_r$: number of cells in each spatial direction of a root grid.
$N_L$: the finest level in a sparse-grid computation.
CPU: CPU time for a complete simulation. CPU time unit: seconds.}}
\label{tab:ex3alin}
\end{table}

\begin{table}[htbp]\footnotesize
\centering
\begin{tabular}{c c c c c c c c}
\hline
\multicolumn{8}{c}{ Single-grid}\\ \hline
& & $N_h\times N_h$		&$L^1$ error		&Order  &$L^\infty$ error    &Order  &CPU(s)			 \\
\hline
&&80 $\times$ 80  &  1.3362e-06 &        - &  4.9201e-06 &          - &     0.572  \\
&&160 $\times$ 160 &  4.2306e-08 &    4.981 &  1.5860e-07 &      4.955 &     7.156  \\
&&320 $\times$ 320 &  1.3119e-09 &    5.011 &  4.9725e-09 &      4.995 &    92.854  \\
&&640 $\times$ 640 &  4.0883e-11 &    5.004 &  1.5560e-10 &      4.998 &  1149.946  \\
\hline
\multicolumn{8}{c}{ Sparse-grid, Lagrange interpolation for prolongation}\\ \hline
$N_r$   &$N_L$  &$N_h\times N_h$	&$L^1$ error   &Order	    &$L^\infty$ error    &Order  &CPU(s)	     \\
\hline
10 & 3 &   80 $\times$ 80 & 5.6865e-05 &        - &  4.0871e-04 &          - &    0.396 \\
20 & 3 & 160 $\times$ 160 & 2.4276e-07 &    7.872 &  1.6913e-06 &      7.917 &    4.896 \\
40 & 3 & 320 $\times$ 320 & 1.2912e-09 &    7.555 &  6.7338e-09 &      7.972 &   60.762 \\
80 & 3 & 640 $\times$ 640 & 4.0720e-11 &    4.987 &  1.5498e-10 &      5.441 &  822.254 \\
\hline
\multicolumn{8}{c}{ Sparse-grid, WENO interpolation for prolongation}\\ \hline
$N_r$   &$N_L$  &$N_h\times N_h$	&$L^1$ error   &Order	    &$L^\infty$ error    &Order  &CPU(s)	     \\
\hline
10 & 3 &   80 $\times$ 80 & 7.1354e-05 &        - &  5.4916e-04 &          - &    0.414 \\
20 & 3 & 160 $\times$ 160 & 2.7404e-07 &    8.024 &  2.1403e-06 &      8.003 &    4.976 \\
40 & 3 & 320 $\times$ 320 & 1.3265e-09 &    7.691 &  6.4093e-09 &      8.383 &   61.510 \\
80 & 3 & 640 $\times$ 640 & 4.0784e-11 &    5.024 &  1.5497e-10 &       5.370 &  781.990 \\
\hline
\end{tabular}
\caption{\footnotesize{Example 3(a), a 2D Burgers' equation is solved by the WENO scheme.
Comparison of numerical errors and CPU times for computations on single-grid and sparse-grid. Both Lagrange and WENO interpolations for prolongation are respectively used in sparse-grid computations.
Final time $T=0.3$.
$N_r$: number of cells in each spatial direction of a root grid.
$N_L$: the finest level in a sparse-grid computation.
CPU: CPU time for a complete simulation. CPU time unit: seconds.}}
\label{tab:ex3aweno}
\end{table}

\begin{table}[htbp]\footnotesize
\centering
\begin{tabular}{c c c c c c c c}
\hline
\multicolumn{8}{c}{ Single-grid}\\ \hline
& & $N_h\times N_h\times N_h$		&$L^1$ error		&Order  &$L^\infty$ error    &Order  &CPU(s)			 \\
\hline
&&80 $\times$ 80 $\times$ 80  &  2.0810e-06 &        - &  7.6743e-06 &          - &      13.692  \\
&&160 $\times$ 160 $\times$ 160 &  6.6581e-08 &    4.966 &  2.4607e-07 &      4.963 &     362.752  \\
&&320 $\times$ 320 $\times$ 320 &  2.0825e-09 &    4.999 &  7.7131e-09 &      4.996 &    9045.370  \\
&&640 $\times$ 640 $\times$ 640 &  6.5097e-11 &        5.000 &  2.4123e-10 &      4.999 &  247133.600  \\
\hline
\multicolumn{8}{c}{ Sparse-grid}\\ \hline
$N_r$   &$N_L$  &$N_h\times N_h\times N_h$	&$L^1$ error   &Order	    &$L^\infty$ error    &Order  &CPU(s)	     \\
\hline
10 & 3 &   80 $\times$ 80 $\times$ 80  & 7.6908e-06 & - & 4.9805e-05 & - & 5.898  \\
20 & 3 & 160 $\times$ 160 $\times$ 160 & 6.8082e-08 & 6.820 & 2.4792e-07 & 7.650 & 103.416  \\
40 & 3 & 320 $\times$ 320 $\times$ 320 & 2.0761e-09 & 5.035 & 7.7136e-09 & 5.006 & 2183.176  \\
80 & 3 & 640 $\times$ 640 $\times$ 640 & 6.5039e-11 & 4.996 & 2.4124e-10 & 4.999 & 53751.975  \\
\hline
\end{tabular}
\caption{\footnotesize{Example 3(b), a 3D Burgers' equation is solved by the linear scheme.
Comparison of numerical errors and CPU times for computations on single-grid and sparse-grid. Lagrange interpolation for prolongation is used in sparse-grid computations.
Final time $T=0.1$.
$N_r$: number of cells in each spatial direction of a root grid.
$N_L$: the finest level in a sparse-grid computation.
CPU: CPU time for a complete simulation. CPU time unit: seconds.}}
\label{tab:ex3-balin}
\end{table}

\begin{table}[htbp]\footnotesize
\centering
\begin{tabular}{c c c c c c c c}
\hline
\multicolumn{8}{c}{ Single-grid}\\ \hline
& & $N_h\times N_h\times N_h$		&$L^1$ error		&Order  &$L^\infty$ error    &Order  &CPU(s)			 \\
\hline
&&80 $\times$ 80 $\times$ 80  &  2.0866e-06 &        - &  7.6725e-06 &          - &      26.706  \\
&&160 $\times$ 160 $\times$ 160 &  6.6687e-08 &    4.968 &  2.4606e-07 &      4.963 &     691.410  \\
&&320 $\times$ 320 $\times$ 320 &  2.0836e-09 &        5.000 &  7.7132e-09 &      4.996 &   16868.220  \\
&&640 $\times$ 640 $\times$ 640 &  6.5107e-11 &        5.000 &  2.4124e-10 &      4.999 &  444972.000  \\
\hline
\multicolumn{8}{c}{ Sparse-grid, Lagrange interpolation for prolongation}\\ \hline
$N_r$   &$N_L$  &$N_h\times N_h\times N_h$	&$L^1$ error   &Order	    &$L^\infty$ error    &Order  &CPU(s)	     \\
\hline
10 & 3 &   80 $\times$ 80 $\times$ 80 &  1.4078e-04 &        - &  1.1642e-03 &          - &      8.632  \\
20 & 3 & 160 $\times$ 160 $\times$ 160&  7.4122e-07 &    7.569 &  9.2799e-06 &      6.971 &    171.818  \\
40 & 3 & 320 $\times$ 320 $\times$ 320&  2.4238e-09 &    8.256 &  2.2712e-08 &      8.675 &   3733.668  \\
80 & 3 & 640 $\times$ 640 $\times$ 640&  6.4885e-11 &    5.223 &  2.4064e-10 &       6.560 &  92752.380 \\
\hline
\multicolumn{8}{c}{ Sparse-grid, WENO interpolation for prolongation}\\ \hline
$N_r$   &$N_L$  &$N_h\times N_h\times N_h$	&$L^1$ error   &Order	    &$L^\infty$ error    &Order  &CPU(s)	     \\
\hline
10 & 3 &   80 $\times$ 80 $\times$ 80  &  1.3225e-04 &        - &  1.1997e-03 &          - &     15.106  \\
20 & 3 & 160 $\times$ 160 $\times$ 160 &  7.6655e-07 &    7.431 &  7.9147e-06 &      7.244 &    223.450  \\
40 & 3 & 320 $\times$ 320 $\times$ 320 &  2.4830e-09 &     8.270 &  2.5650e-08 &      8.269 &   4144.632  \\
80 & 3 & 640 $\times$ 640 $\times$ 640 &  6.4883e-11 &    5.258 &  2.4056e-10 &      6.736 &  96962.960  \\
\hline
\end{tabular}
\caption{\footnotesize{Example 3(b), a 3D Burgers' equation is solved by the WENO scheme.
Comparison of numerical errors and CPU times for computations on single-grid and sparse-grid. Both Lagrange and WENO interpolations for prolongation are respectively used in sparse-grid computations.
Final time $T=0.1$.
$N_r$: number of cells in each spatial direction of a root grid.
$N_L$: the finest level in a sparse-grid computation.
CPU: CPU time for a complete simulation. CPU time unit: seconds.}}
\label{tab:ex3-baweno}
\end{table}

\subsection{Nonlinear case with discontinuous solution}

In this section, we test the fifth order sparse grid WENO scheme in solving nonlinear problem which has shock waves developed in the solution.

\bigskip
\noindent{\bf Example 4 (A 3D Burgers' equation with non-smooth solution).}

\begin{figure}[p]
    \centering
    \begin{subfigure}[b]{0.48\textwidth}
        \includegraphics[width=0.9\textwidth]{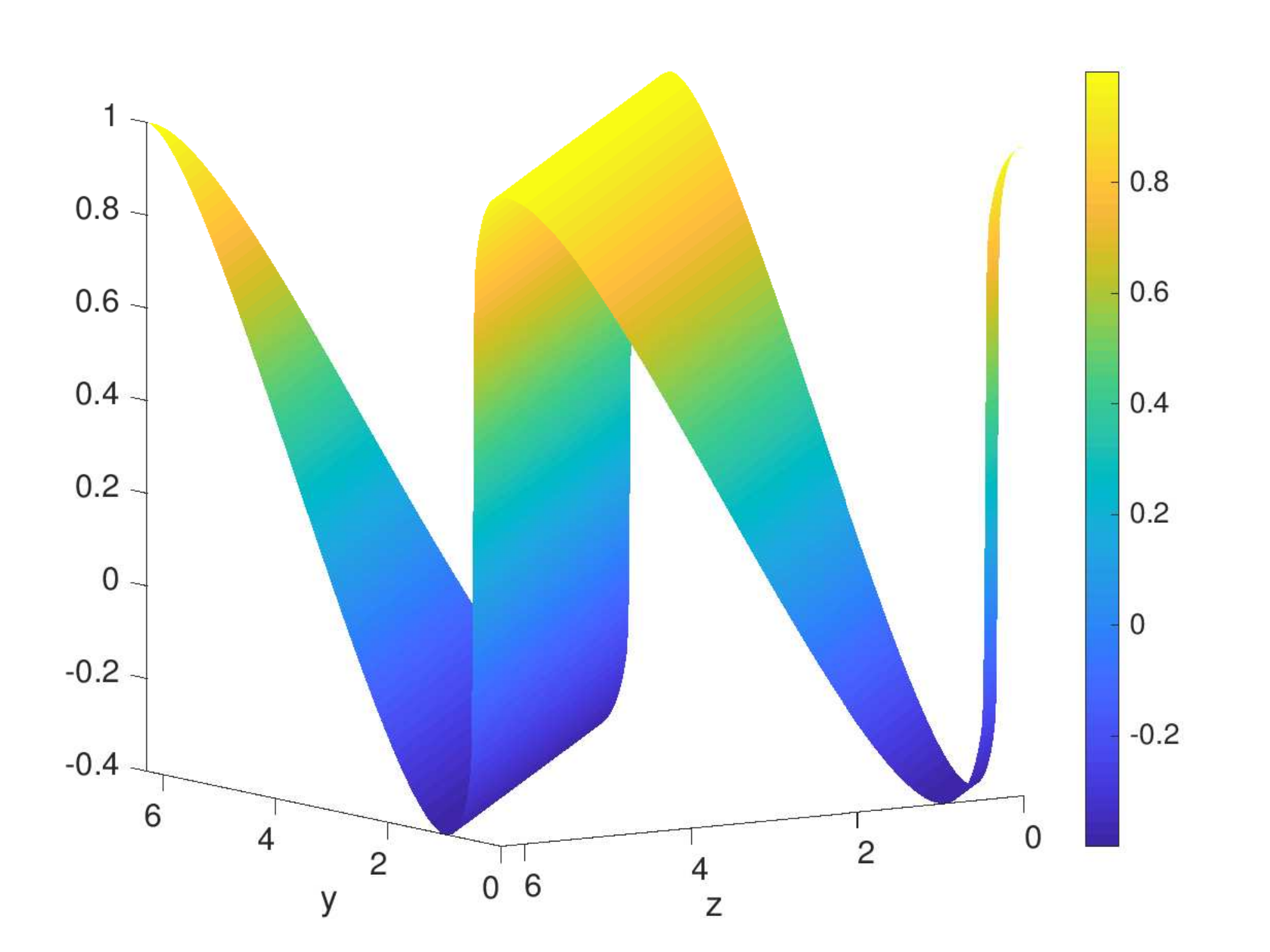}
        \caption{\footnotesize{single-grid solution on $y-z$ plane at $x = 0$}}
        \label{result on regular grid}
    \end{subfigure}
        \begin{subfigure}[b]{0.48\textwidth}
        \includegraphics[width=0.9\textwidth]{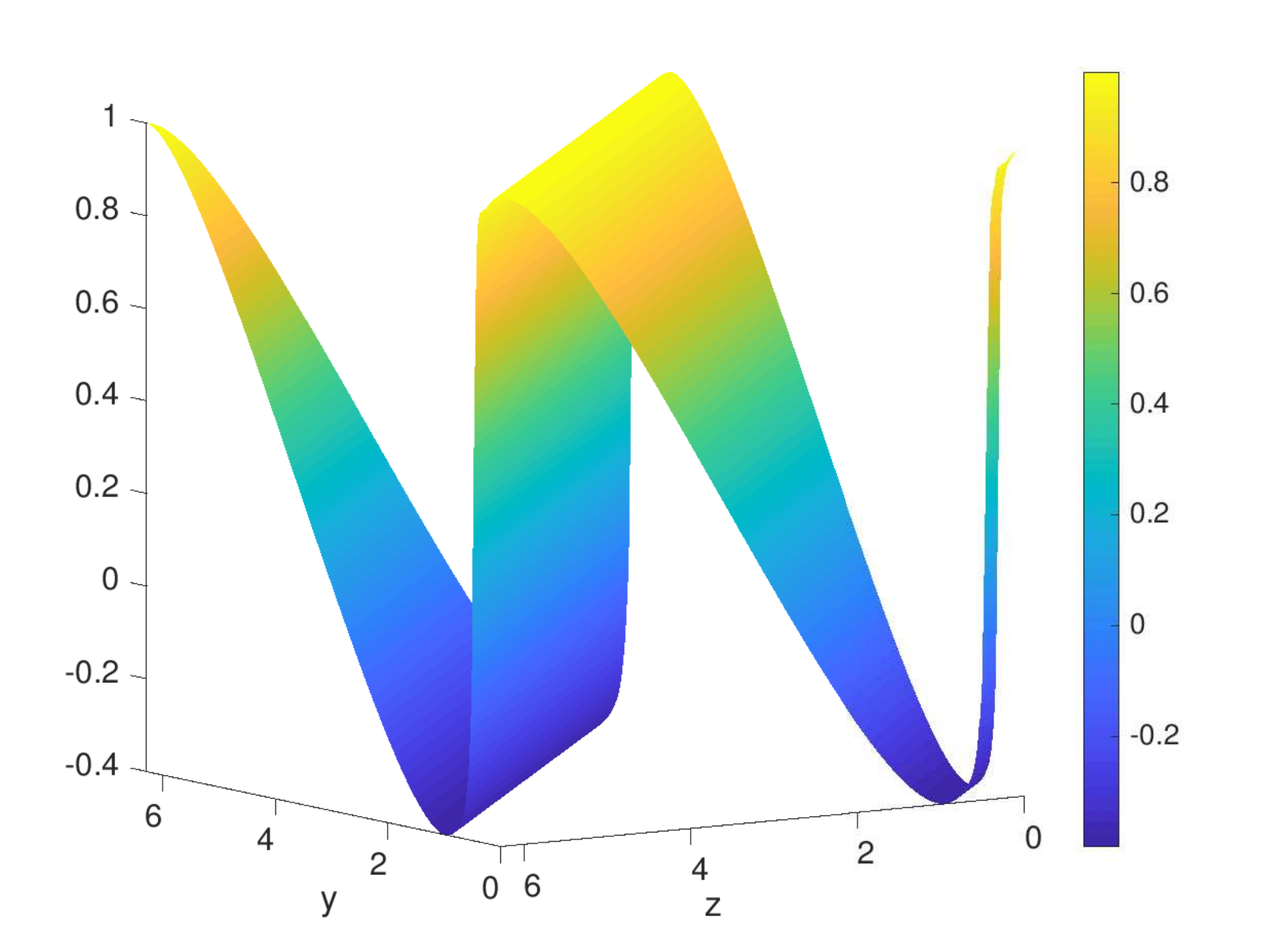}
        \caption{\footnotesize{sparse-grid solution on $y-z$ plane at $x = 0$}}
        \label{result on sparse grids}
    \end{subfigure}
    \begin{subfigure}[b]{0.48\textwidth}
        \includegraphics[width=0.9\textwidth]{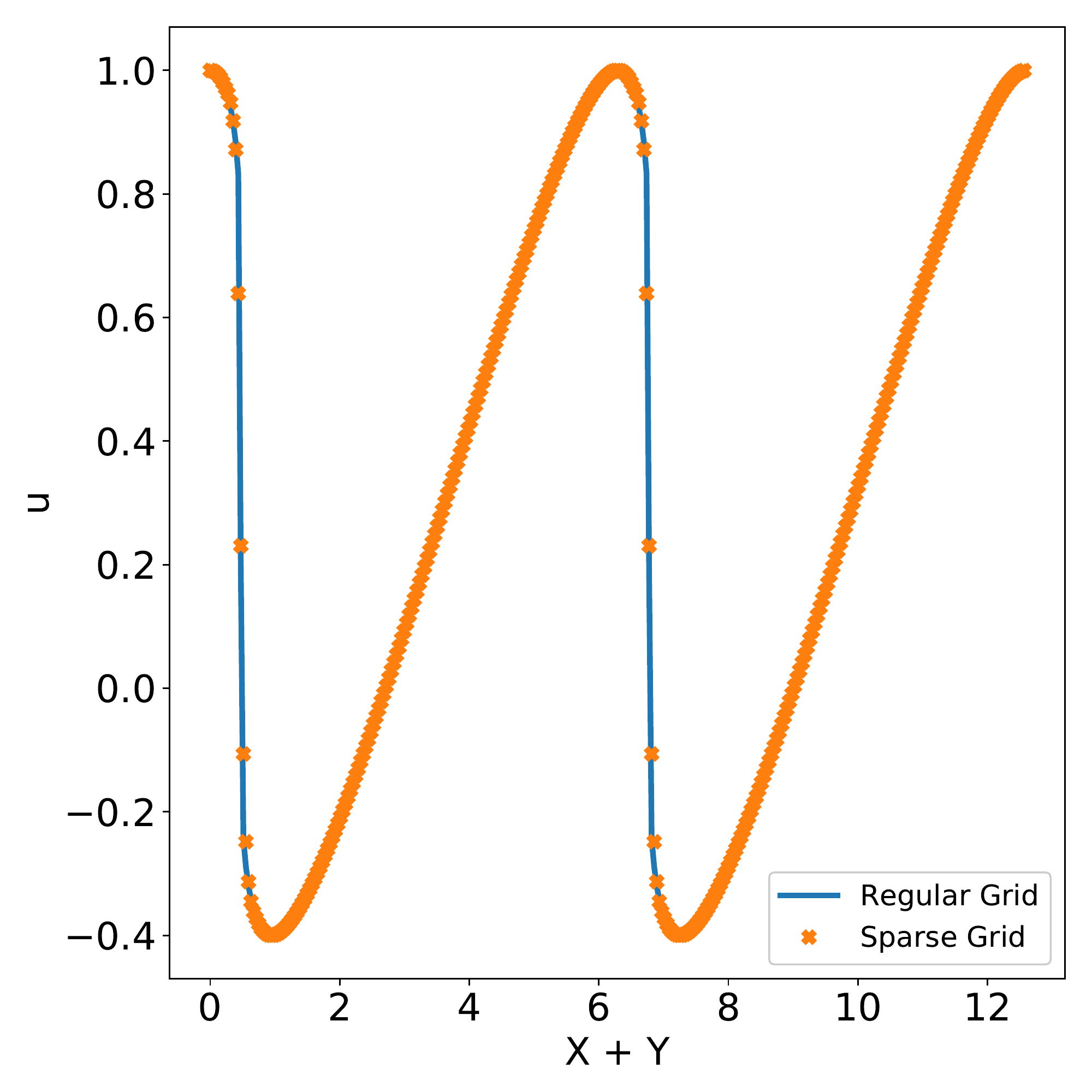}
        \caption{\footnotesize{Solutions along the line $x=y$ on plane $z=0$}}
        \label{1dcut1}
    \end{subfigure}
        \begin{subfigure}[b]{0.48\textwidth}
        \includegraphics[width=0.9\textwidth]{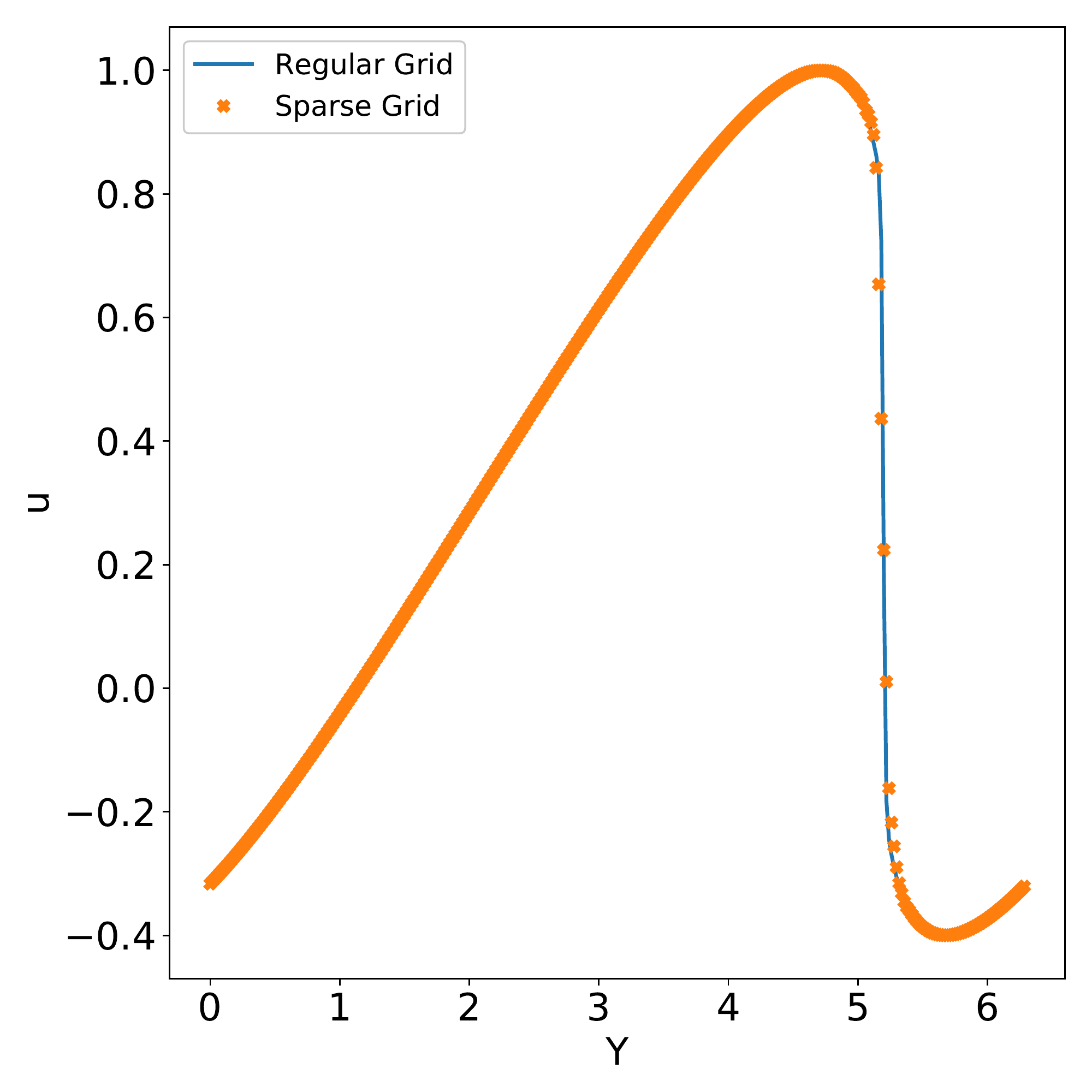}
        \caption{\footnotesize{Solutions along the line $z = 0, x = 1.5710$}}
        \label{1dcut2}
    \end{subfigure}
    \caption{\footnotesize{Example 4, solution of a 3D Burgers' equation at the time $T = 0.52$ by the fifth order WENO scheme on sparse grids ($N_r=40$ for root grid, finest level $N_L=3$ in the sparse-grid computation) and the corresponding $320 \times 320\times 320$ single grid, using the fifth order WENO interpolation for prolongation in the sparse grid combination.  $CFL=0.4$. 2D and 1D cutting-plots at different locations of the domain. Blue lines: results using the single grid; crosses: results using the sparse grids.}}
    \label{3d_burg_discon_comp}
\end{figure}

\noindent The 3D Burgers' equation is solved till a time when shock waves in the solution have formed. The equation is
\begin{equation*}
    u_t~+~\Big(\frac{u^2}{2}\Big)_x~+~\Big(\frac{u^2}{2}\Big)_y~+~\Big(\frac{u^2}{2}\Big)_z~=~0,
\end{equation*}
with the initial condition
\begin{equation}
    u_0(x,y,z) = 0.3+0.7\sin(x+y+z)
\end{equation}
and periodic boundary conditions.
The computational domain is $[0,~2\pi]\times[0,~2\pi]\times[0,~2\pi]$.
We use the fifth order WENO scheme on sparse grids and the corresponding single grids to compute the numerical solution of the problem at the final time $T = 0.52$ when the discontinuities in the solution have been developed. The WENO interpolation for prolongation is used in the sparse grid WENO scheme. The sparse grid root grid $N_r=40$, and the finest level $N_L=3$. Hence the most refined mesh in the sparse grids or the corresponding single grid has $320\times 320\times 320$ computational cells.  The CFL number is taken to be $0.4$ in the simulations. In Figure
\ref{3d_burg_discon_comp}, we show the 2D and 1D cutting-plots of the numerical solutions at different locations of the domain for both the sparse-grid computation and the corresponding single-grid computation. It is clearly seen that
the numerical solution by the fifth order sparse grid WENO scheme with the WENO prolongation is similar as that by the single grid WENO scheme. The nonlinear stability and high resolution properties of the fifth order WENO scheme for resolving shock waves are preserved well in the sparse-grid computation. However, the sparse-grid computation is much more efficient than the corresponding single-grid computation. In this example, it takes only $17893.4$ seconds of CPU time to complete the simulation of the sparse-grid computation, while $86549.8$ seconds of CPU time are needed for finishing the simulation of the corresponding single-grid computation. About $80\%$ CPU time is saved by performing the fifth order WENO simulation on the sparse grids for this 3D example with discontinuous solution.

\bigskip
\noindent{\bf Remark 1.} In this example which has a discontinuous solution such as shock wave, we show that the fifth order WENO scheme with the sparse-grid combination techniques can stably capture discontinuities in their correct locations and preserve the essentially non-oscillatory property. However, for more complex problems with very strong shock waves in applications such as compressible fluid dynamics, it is possible to have the issue of bound-preserving in our proposed sparse grid methods, for example, to satisfy the maximum principle for scalar conservation laws, and to preserve the positivity property of density and internal energy for Euler systems, etc. We will investigate such cases in the next research of this topic. A possible approach is to apply the bound-preserving techniques, e.g. \cite{Shu3}, to our proposed sparse grid methods.

\subsection{Application to kinetic simulations}
In this section, we apply the fifth order sparse grid WENO method in simulating the kinetic equations.
The solution of the distribution function in the kinetic equations
may develop violent gradients and sharp transitions. If the numerical methods are not designed carefully, a massive formation of spurious oscillations will lead to a false physical description of the phenomena. WENO methods have been a class of successful approaches to deal with this issue, e.g., see \cite{Car1, Car2, QC}.
Another major challenge in deterministic kinetic simulations is due to the high dimensionality of the related PDE systems,
which can have up to six dimensions in spatial directions of the PDEs including both space and velocity variables.
As examples to show the high efficiency of the sparse grid WENO method for simulating such systems, we solve a Vlasov-Boltzmann transport equation in both 2D and 4D cases, and a simplified 3D Vlasov-Maxwell system.

\bigskip
\noindent{\bf Example 5 (A Vlasov-Boltzmann transport equation).} We consider a collisional relaxation model described
by the Vlasov-Boltzmann transport equation (see e.g. \cite{CGP} and references therein). In \cite{GC}, a sparse grid DG
method is used to solve the model. The PDE has the following form
\begin{equation}
    f_t+\mathbf{v}\cdot\nabla_\mathbf{x}f+\mathbf{E}(\mathbf{x})\cdot\nabla_\mathbf{v}f=L(f).
\end{equation}
Here the unknown function $f$ depends on space variables $\mathbf{x}$, velocity variables $\mathbf{v}$ and the time variable $t$, i.e., $f=f(t,\mathbf{x},\mathbf{v})$. It denotes the probability distribution function of electrons.
The external electric field $\mathbf{E}(\mathbf{x})$ is given by a known electrostatic potential
\begin{equation}
    \mathbf{E}(\mathbf{x})=-\nabla_\mathbf{x}\Phi(\mathbf{x}), \qquad\Phi(\mathbf{x})=\frac{|\mathbf{x}|^2}{2}.
\end{equation}
$L(f)$ is the linear relaxation operator defined as
\begin{equation}
    L(f)\triangleq\frac{\mu_\infty(\mathbf{v})\rho(t,\mathbf{x})-f(t,\mathbf{x},\mathbf{v})}{\tau}.
\end{equation}
The absolute Maxwellian distribution $\mu_\infty$ is defined by
\begin{equation}
    \mu_\infty(\mathbf{v})\triangleq\frac{e^{-\frac{\mathbf{|v|}^2}{2\theta}}}{(2\pi\theta)^{d/2}},
\end{equation}
and the macroscopic density
\begin{equation}
\label{eq:integral_f_v}
    \rho(t,\mathbf{x})=\int_\mathbf{v}f(t,\mathbf{x},\mathbf{v})d\mathbf{v}.
\end{equation}
Here the parameters $\theta$ is the kinetic temperature, $\tau=1/k$ with $k$ being the constant transition
probability of scatters passing from one state into another one. $d$ is the space dimension.

\begin{figure}[p]
    \centering
    \begin{subfigure}[b]{0.48\textwidth}
        \includegraphics[width=0.9\textwidth]{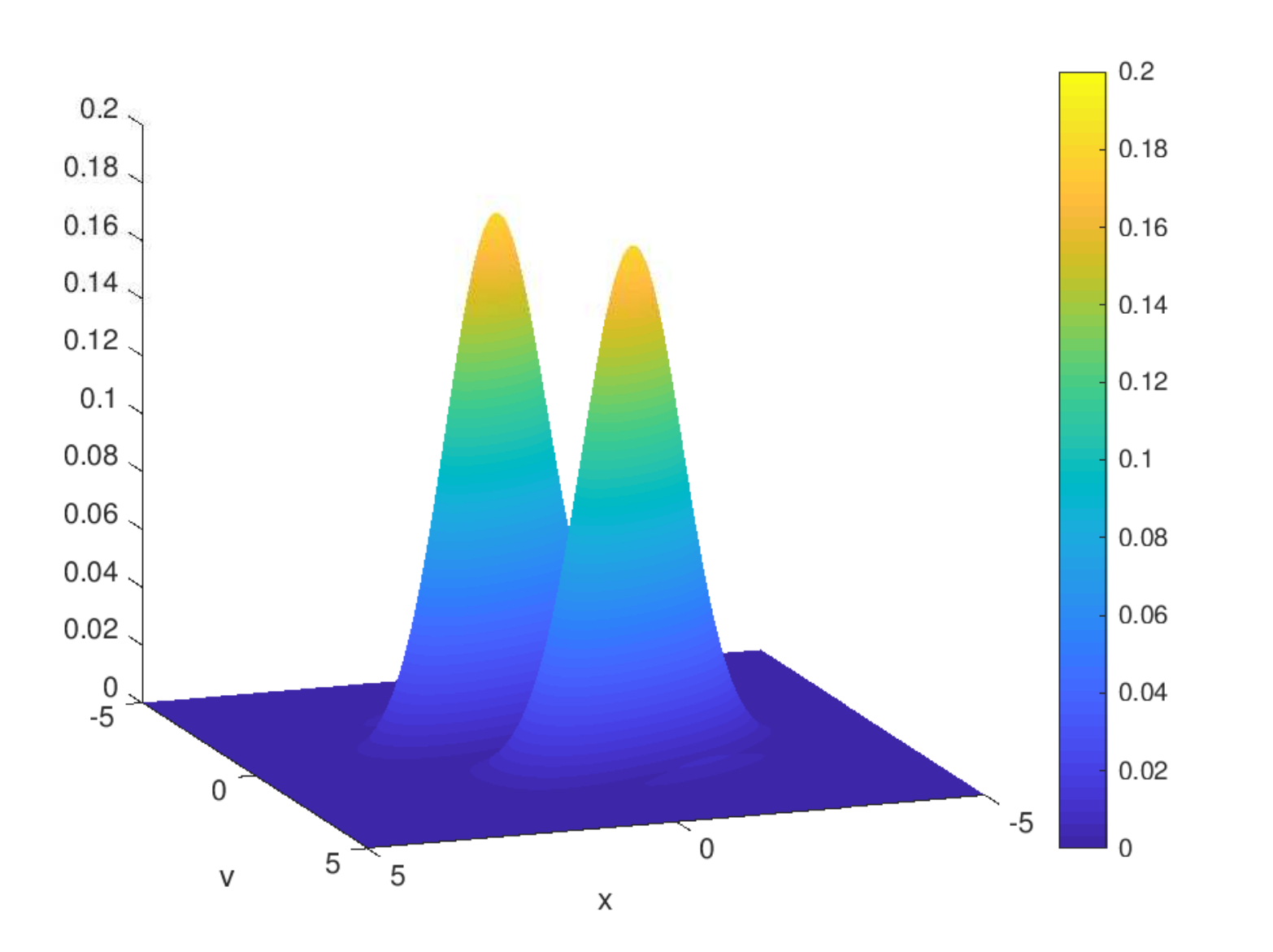}
        \caption{single-grid result at $t=0.5$}
        \label{single-t0.5}
    \end{subfigure}
        \begin{subfigure}[b]{0.48\textwidth}
        \includegraphics[width=0.9\textwidth]{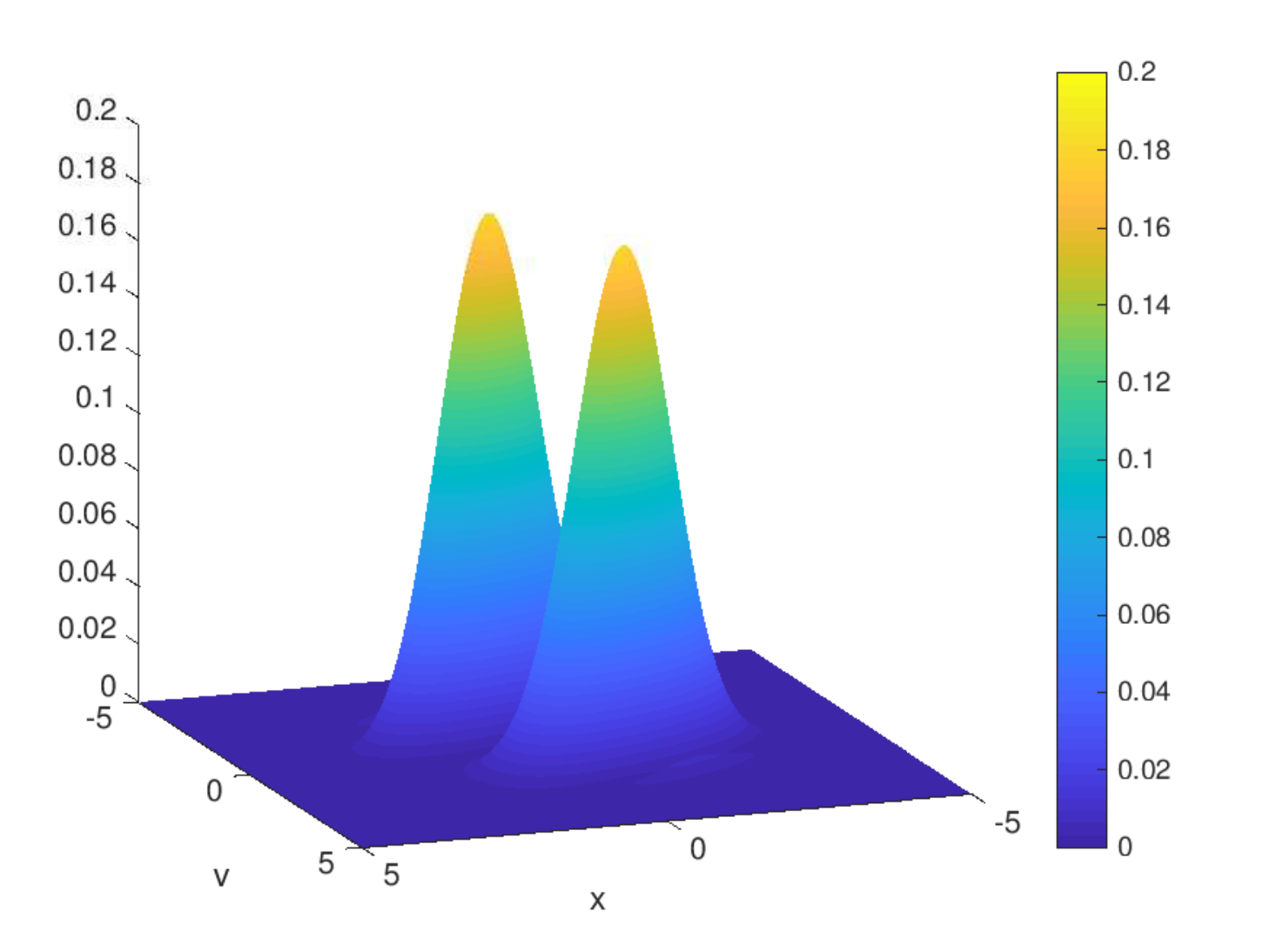}
        \caption{sparse-grid result at $t=0.5$}
        \label{spar-t0.5}
    \end{subfigure}
        \begin{subfigure}[b]{0.48\textwidth}
        \includegraphics[width=0.9\textwidth]{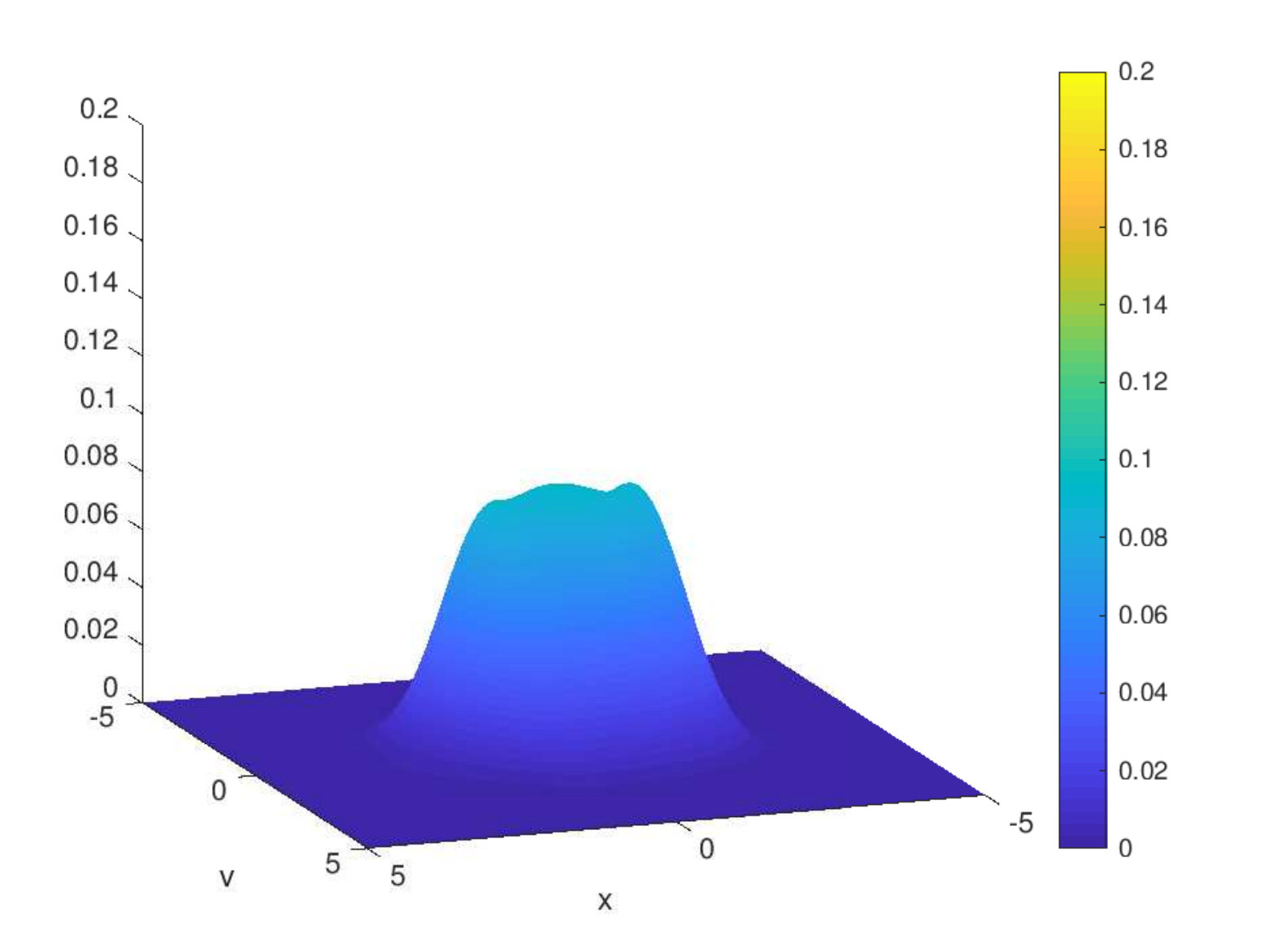}
        \caption{single-grid result at $t=2.0$}
        \label{single-t2}
    \end{subfigure}
        \begin{subfigure}[b]{0.48\textwidth}
        \includegraphics[width=0.9\textwidth]{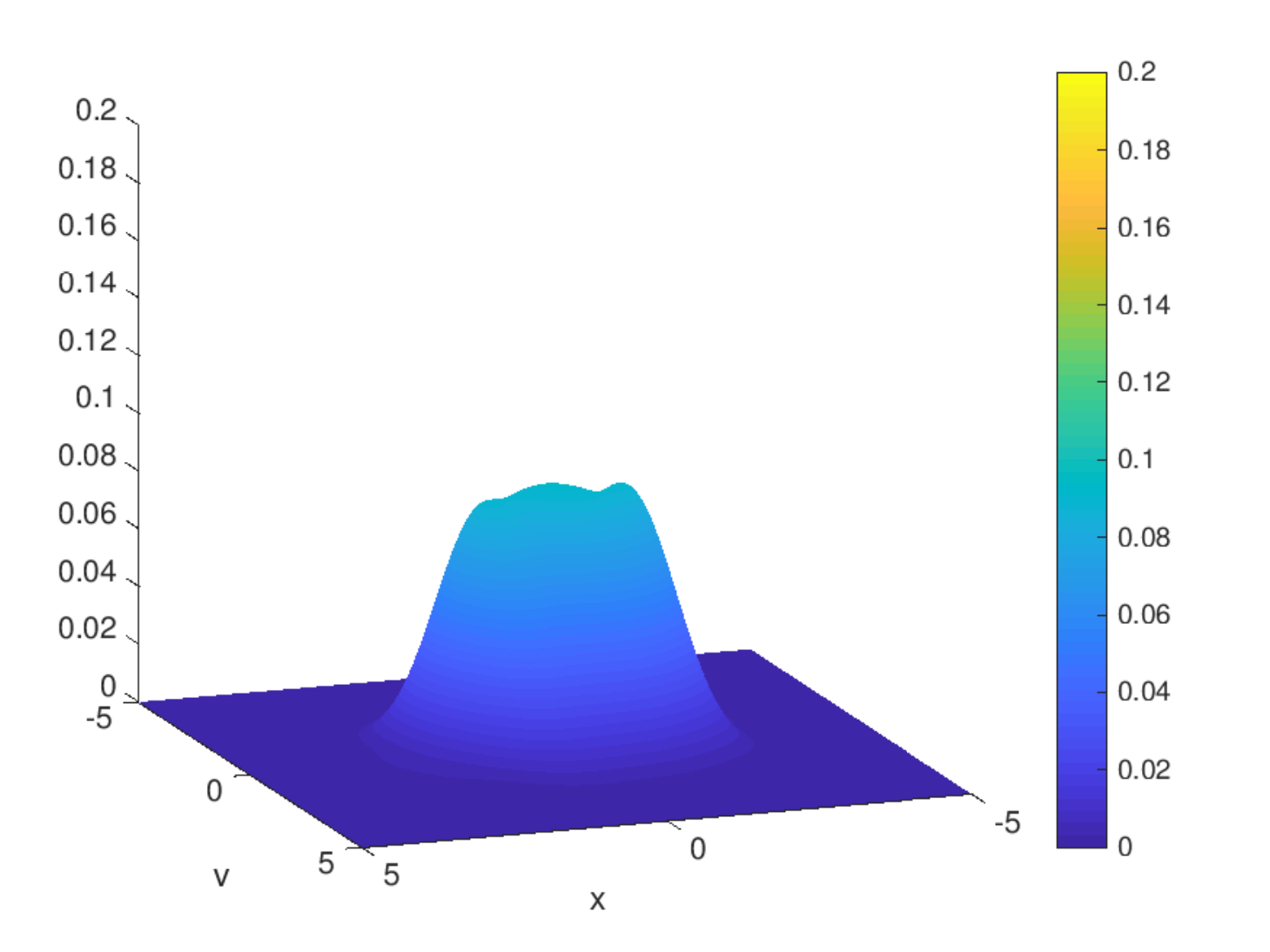}
        \caption{sparse-grid result at $t=2.0$}
        \label{sparse-t2}
    \end{subfigure}
        \begin{subfigure}[b]{0.48\textwidth}
        \includegraphics[width=0.9\textwidth]{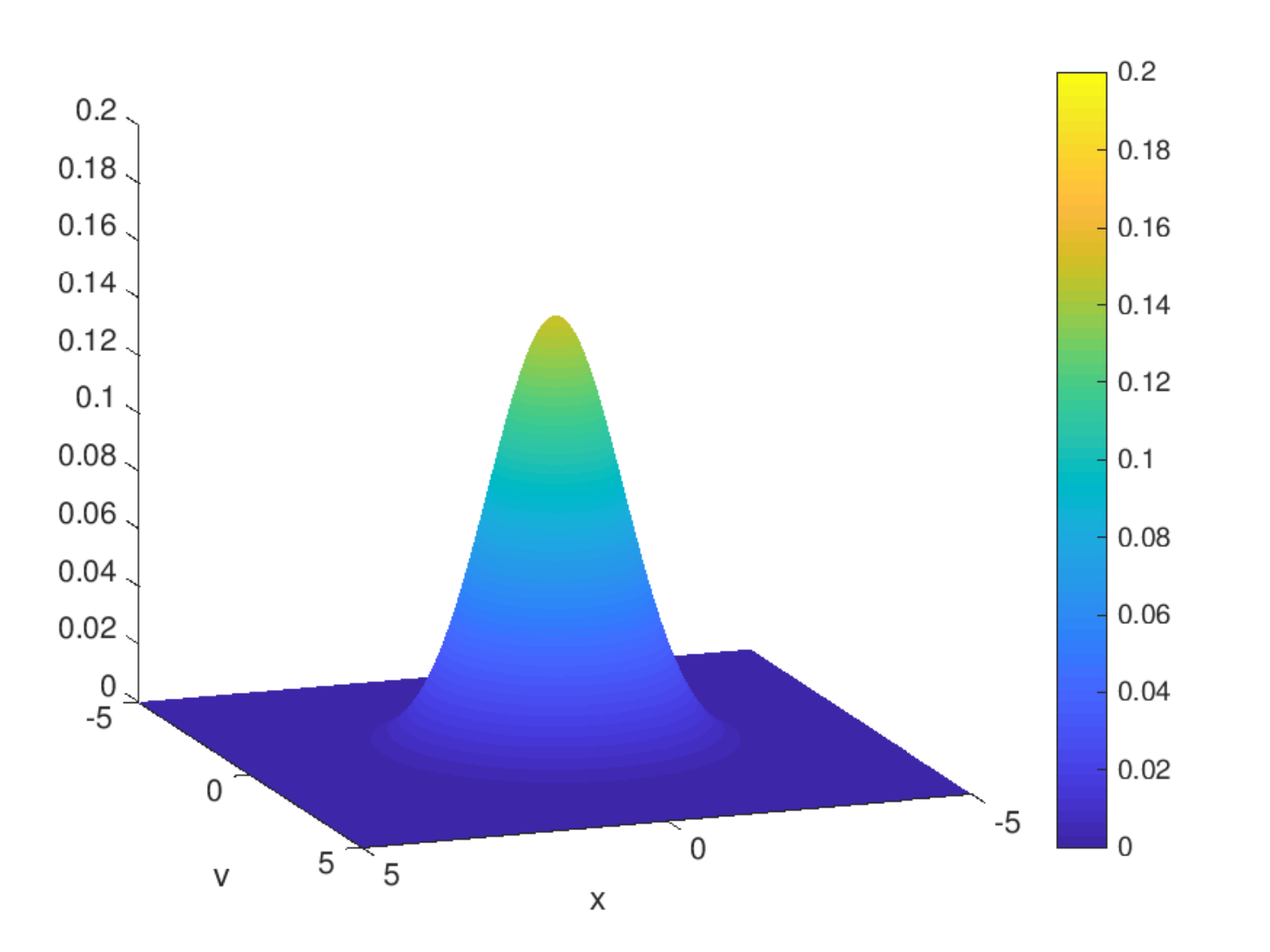}
        \caption{single-grid result at $t=6.0$}
        \label{single-t6}
    \end{subfigure}
        \begin{subfigure}[b]{0.48\textwidth}
        \includegraphics[width=0.9\textwidth]{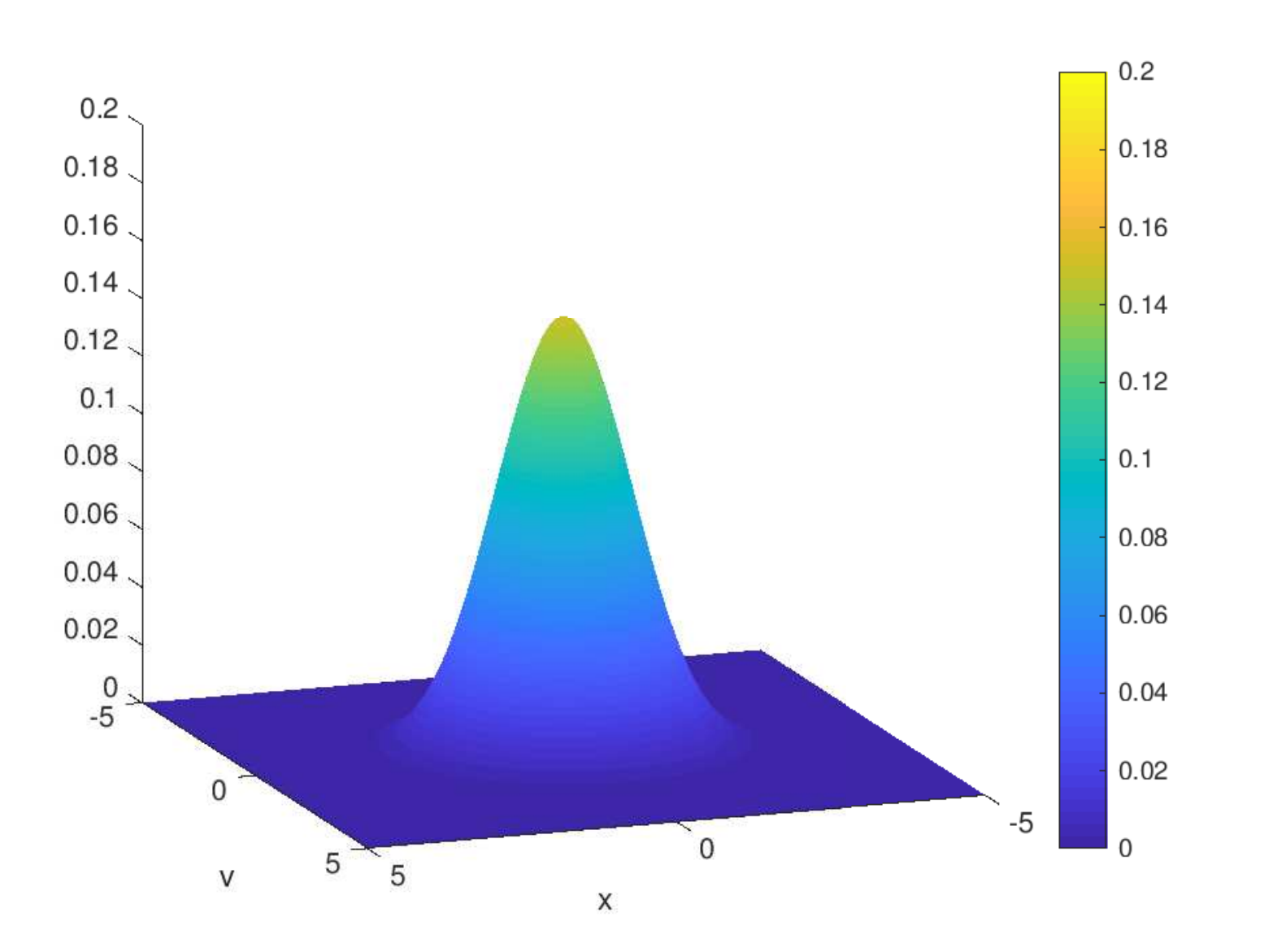}
        \caption{sparse-grid result at $t=6.0$}
        \label{sparse-t6}
    \end{subfigure}
    \caption{\footnotesize{Example 5(a), solution $f$ of the two dimensional Vlasov-Boltzmann transport equation by fifth order WENO scheme on sparse grids ($N_r=40$ for root grid, finest level $N_L=3$ in the sparse-grid computation) and the corresponding $320 \times 320$ single grid, at different time $t$. $CFL=0.4$. (a), (c), (e): single-grid results; (b), (d), (f): sparse-grid results. }}
\label{2dplot_vlasovBTE}
\end{figure}

\begin{figure}[p]
    \centering
    \begin{subfigure}[b]{0.48\textwidth}
        \includegraphics[width=0.9\textwidth]{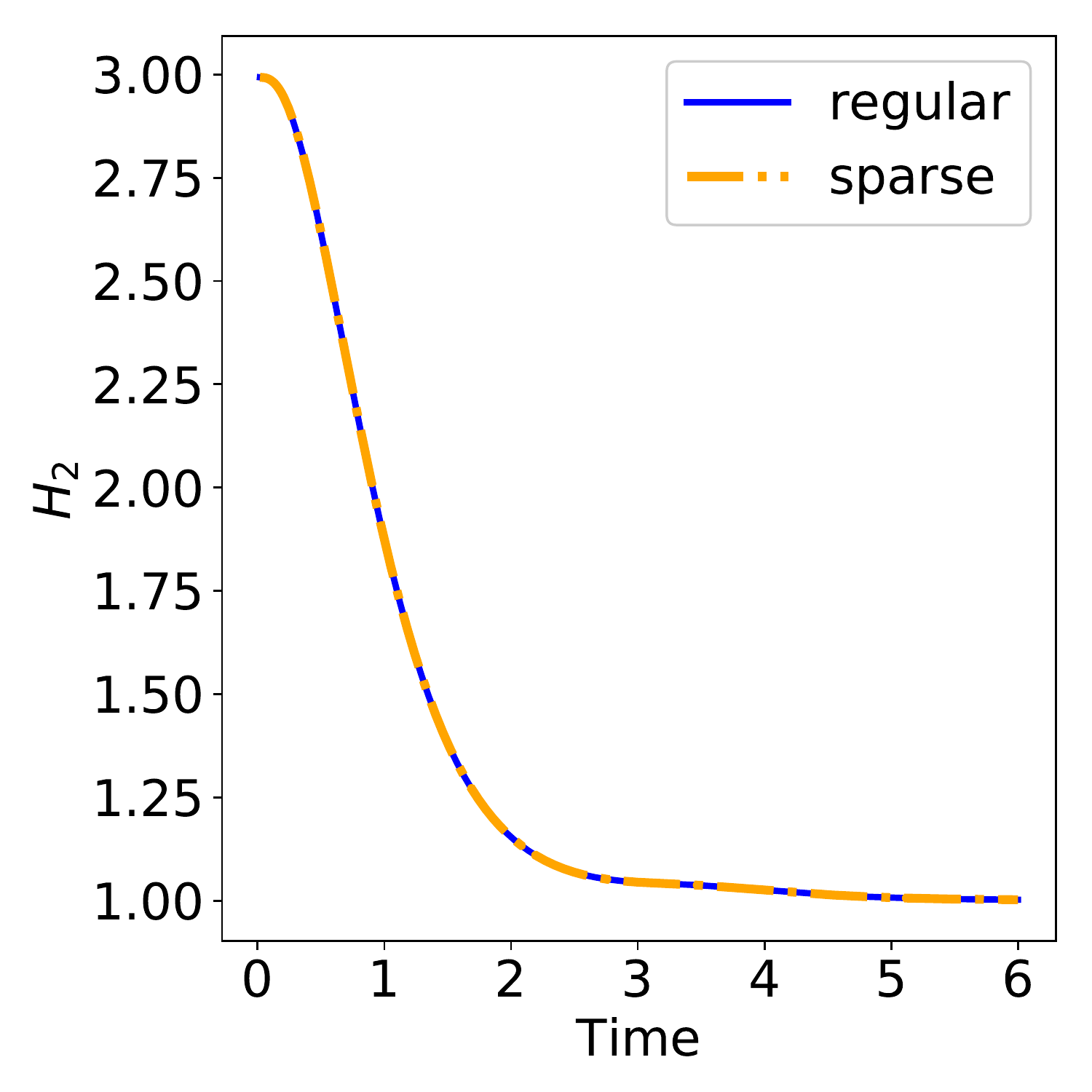}
        \caption{$\mathcal{H}_2$ entropy functional.}
        \label{h2}
    \end{subfigure}
        \begin{subfigure}[b]{0.48\textwidth}
        \includegraphics[width=0.9\textwidth]{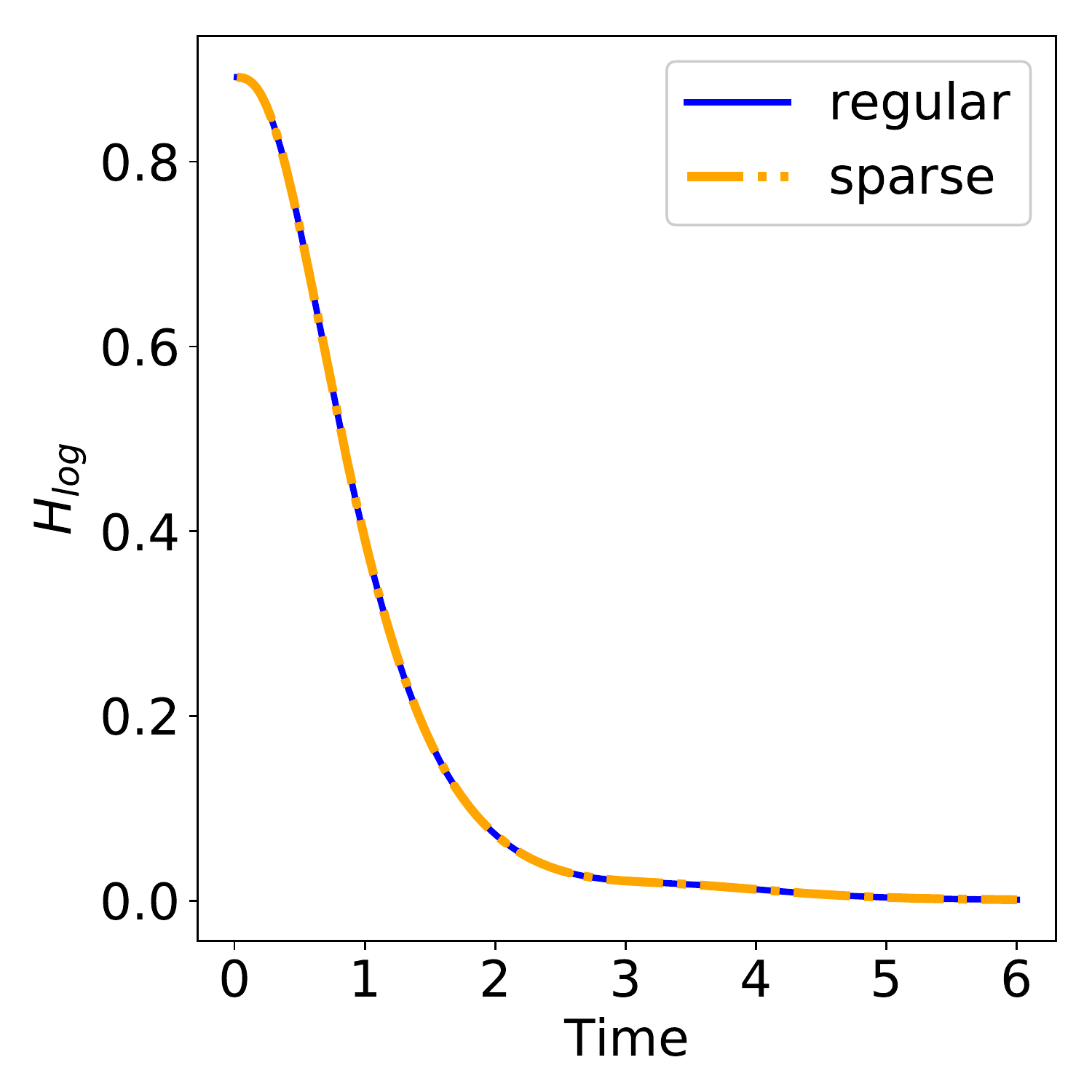}
        \caption{$\mathcal{H}_{\log}$ entropy functional.}
        \label{hlog}
    \end{subfigure}
    \caption{\footnotesize{Example 5(a), Comparison of decay rate for entropy functionals obtained by fifth order WENO scheme on sparse grids ($N_r=40$ for root grid, finest level $N_L=3$ in the sparse-grid computation) and the corresponding $320 \times 320$ single grid. $CFL=0.4$. Orange lines with dots: results by the sparse-grid computation; blue lines: results by the corresponding single-grid computation.}}
    \label{2d_vlasov_entropy}
\end{figure}

{\bf(a) 2D case.} First we solve the 2D case by the fifth order sparse grid WENO scheme. $d = 1$ and parameters $\theta=\tau=1$. The computation domain is $\Omega=[-5,5]\times[-5,5]$. The initial condition is
\begin{align}
    f(0, x, v)&=\frac{1}{s}\sin{(\frac{x^2}{2})}^2e^{(-\frac{x^2+v^2}{2})},
\end{align}
where $s$ is the normalization constant such that $\int_{\Omega}f(0, x, v)dxdv=1$.
The zero boundary conditions are prescribed at the domain boundaries.
Sparse grids with a $40 \times 40 $ root grid
and the most refined level $3$ are used. So the corresponding single grid is $320 \times 320 $.
We compare the simulation results of sparse-grid computation and
the corresponding single-grid computation. The CFL number is taken to be $0.4$ in the simulations.
From Figure \ref{2dplot_vlasovBTE}, we can observe that sparse-grid computations and
the corresponding single-grid computations generate similar results. For a further comparison,
we study the decay rate of the initial state to equilibrium of the solution. As that in \cite{CGP, GC},
time evolutions of two entropy functionals are tracked. They are defined as following:
\begin{equation}
    \mathcal{H}_2(t)~=~\int_\Omega H^2\mathcal{M}dxdv,~~~~~~~\mathcal{H}_{\log}(t)~=~\int_\Omega H\log(H)\mathcal{M}dxdv,
\end{equation}
where $H(t, x, v)=f(t,x,v)/\mathcal{M}(x,v)$ and $\mathcal{M}(x,v)$ is the unique stationary state solution of the system \cite{CGP}
\begin{equation}\label{eq:integral_M}
    \mathcal{M}(x,v)=\frac{e^{-(\frac{|v|^2}{2}+\Phi(x))/\theta}}{(2\pi\theta)^{d/2}\int_x e^{-\Phi(x)/\theta}dx}.
\end{equation}
From Figure \ref{2d_vlasov_entropy}, we can observe that the decay rates for entropy functionals obtained by
sparse-grid computation and
the corresponding single-grid computation agree to each other very well. We report the CPU times for the simulations on sparse grids and the corresponding single grids at different time $t$ in Table \ref{tb:ex5acpu}. Close to $30\%$ CPU times are saved for a longer time $t$ by performing the fifth order WENO simulation on the sparse grids for this 2D problem. As shown in the following 4D case, a much more significant CPU time saving in sparse-grid computation is obtained for a higher dimensional problem.

\begin{table}[htbp]\footnotesize
    \begin{tabular}{ c|c|c|c }
        \hline
        Time $t$ & CPU Time(s) of single-grid comp.& CPU Time(s) of sparse-grid comp. & sparse/single ratio \\\hline
        0.5 & 13111.313 & 11720.746  & 0.8939  \\
        1.0 & 25185.141  & 21371.762  & 0.8486 \\
        2.0 & 50756.395  & 39493.979  & 0.7781 \\
        3.0 & 82253.965  & 57898.835  & 0.7039\\
        6.0 & 153593.745  &110626.464  & 0.7203\\
        \hline
    \end{tabular}
    \caption{\footnotesize{Example 5(a), CPU times (unit: seconds) for the simulations on sparse grids and the corresponding single grids, at different time $t$.}}
    \label{tb:ex5acpu}
\end{table}

\begin{figure}[p]
    \centering
    \begin{subfigure}[b]{0.48\textwidth}
        \includegraphics[width=0.9\textwidth]{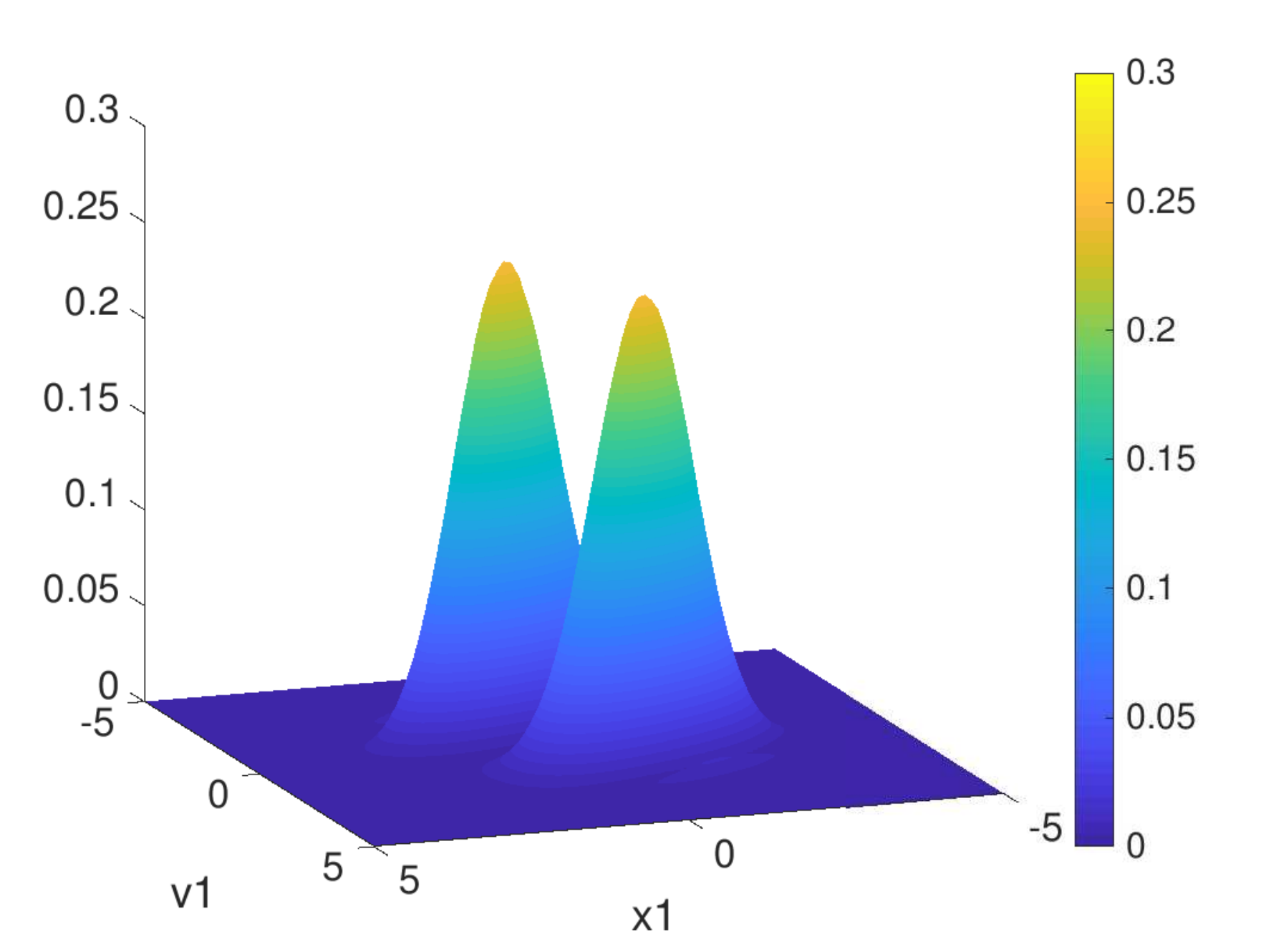}
        \caption{single-grid result at $t=0.5$}
        \label{4dsingle-t0.5}
    \end{subfigure}
        \begin{subfigure}[b]{0.48\textwidth}
        \includegraphics[width=0.9\textwidth]{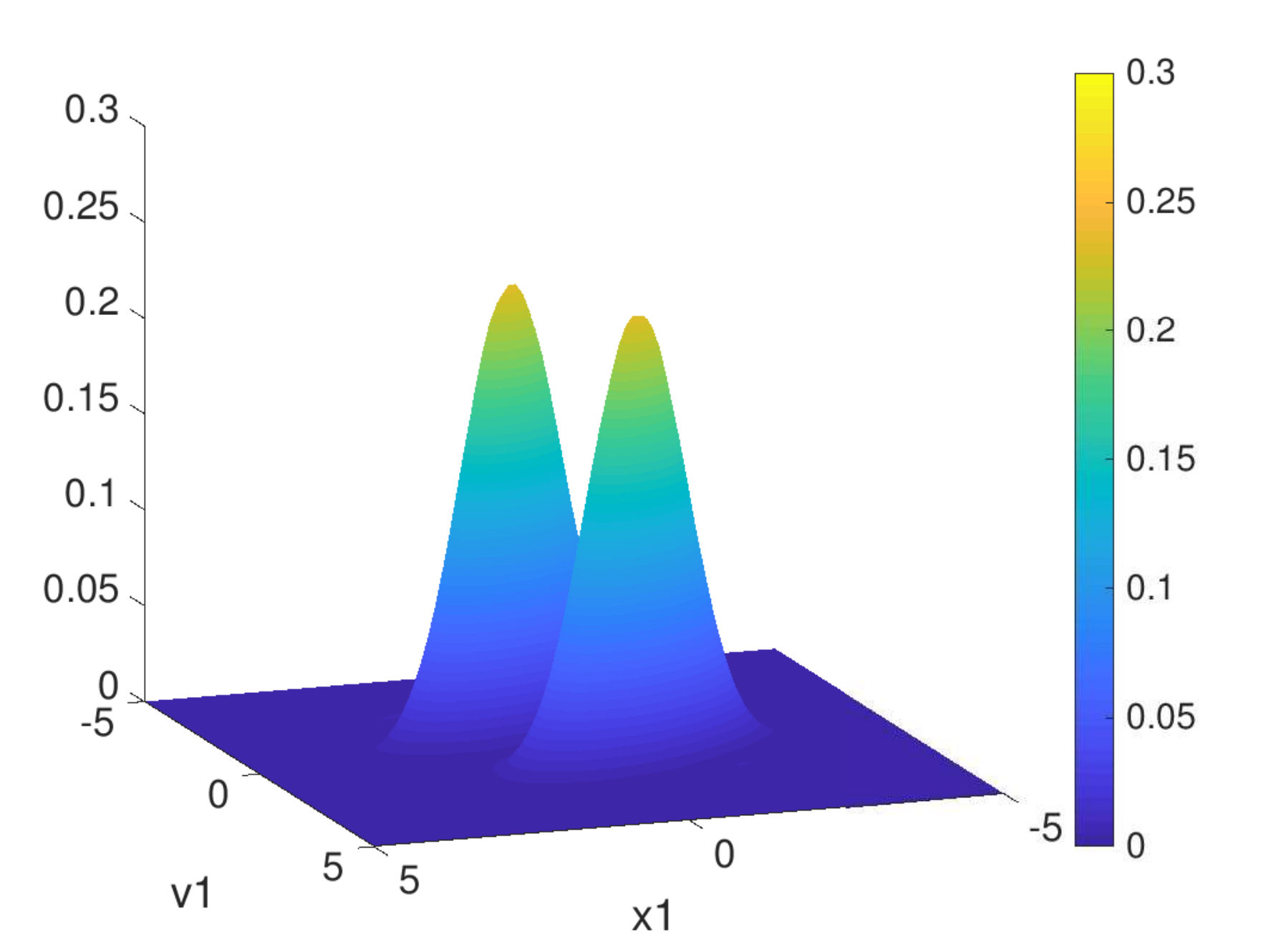}
        \caption{sparse-grid result at $t=0.5$}
        \label{4dspar-t0.5}
    \end{subfigure}
        \begin{subfigure}[b]{0.48\textwidth}
        \includegraphics[width=0.9\textwidth]{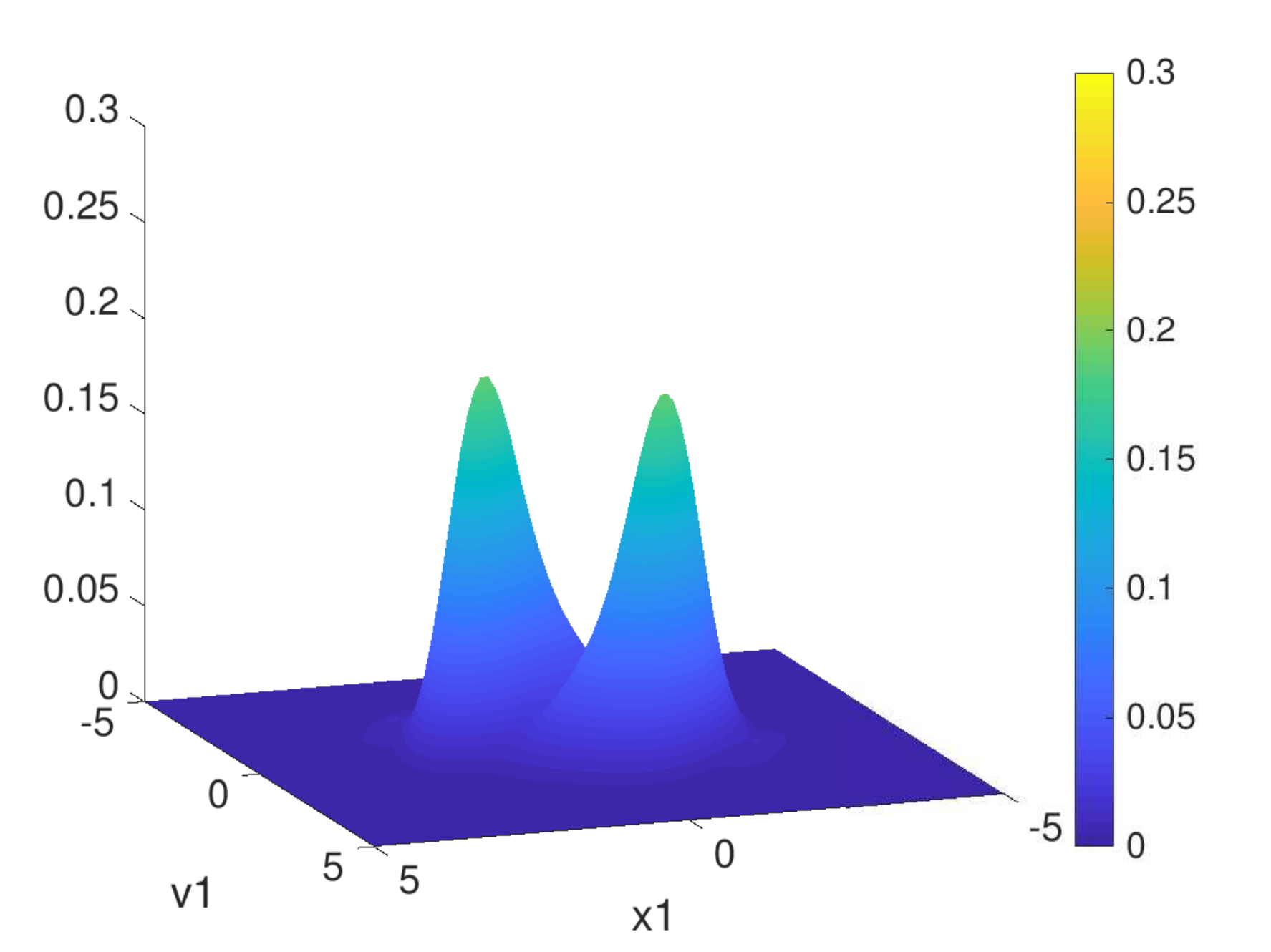}
        \caption{single-grid result at $t=1.0$}
        \label{4dsingle-t1}
    \end{subfigure}
        \begin{subfigure}[b]{0.48\textwidth}
        \includegraphics[width=0.9\textwidth]{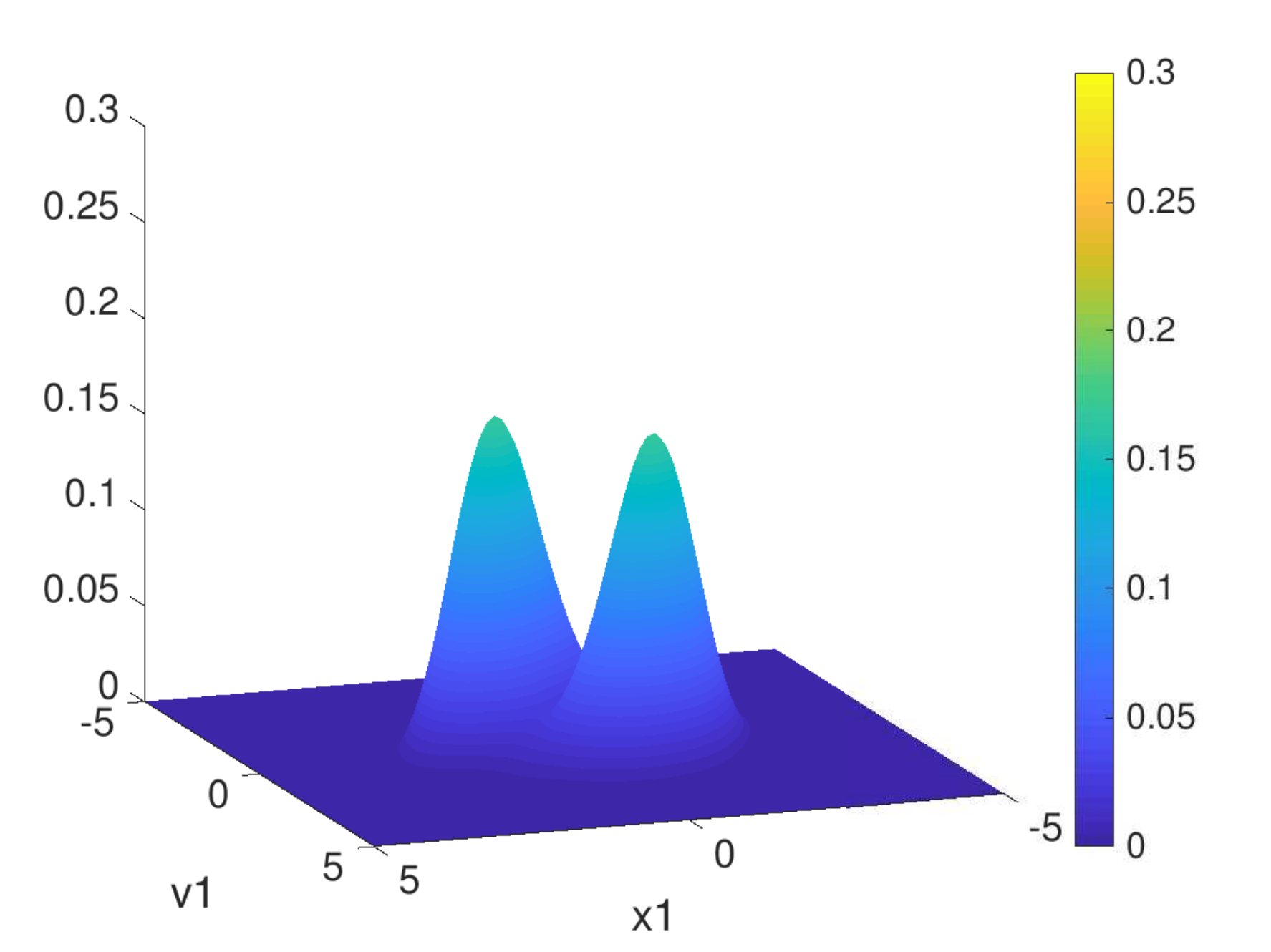}
        \caption{sparse-grid result at $t=1.0$}
        \label{4dsparse-t1}
    \end{subfigure}
        \begin{subfigure}[b]{0.48\textwidth}
        \includegraphics[width=0.9\textwidth]{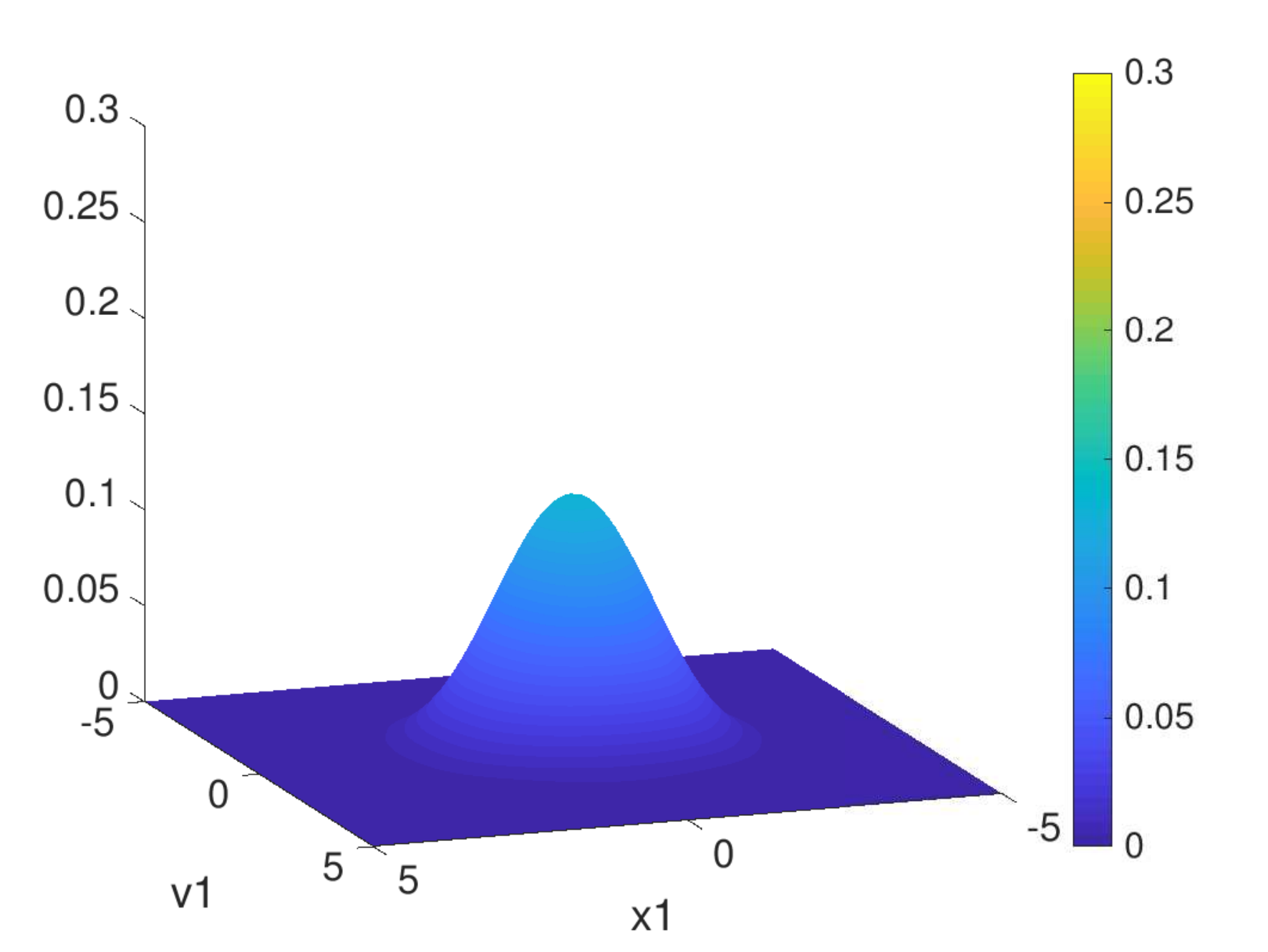}
        \caption{single-grid result at $t=3.0$}
        \label{4dsingle-t3}
    \end{subfigure}
        \begin{subfigure}[b]{0.48\textwidth}
        \includegraphics[width=0.9\textwidth]{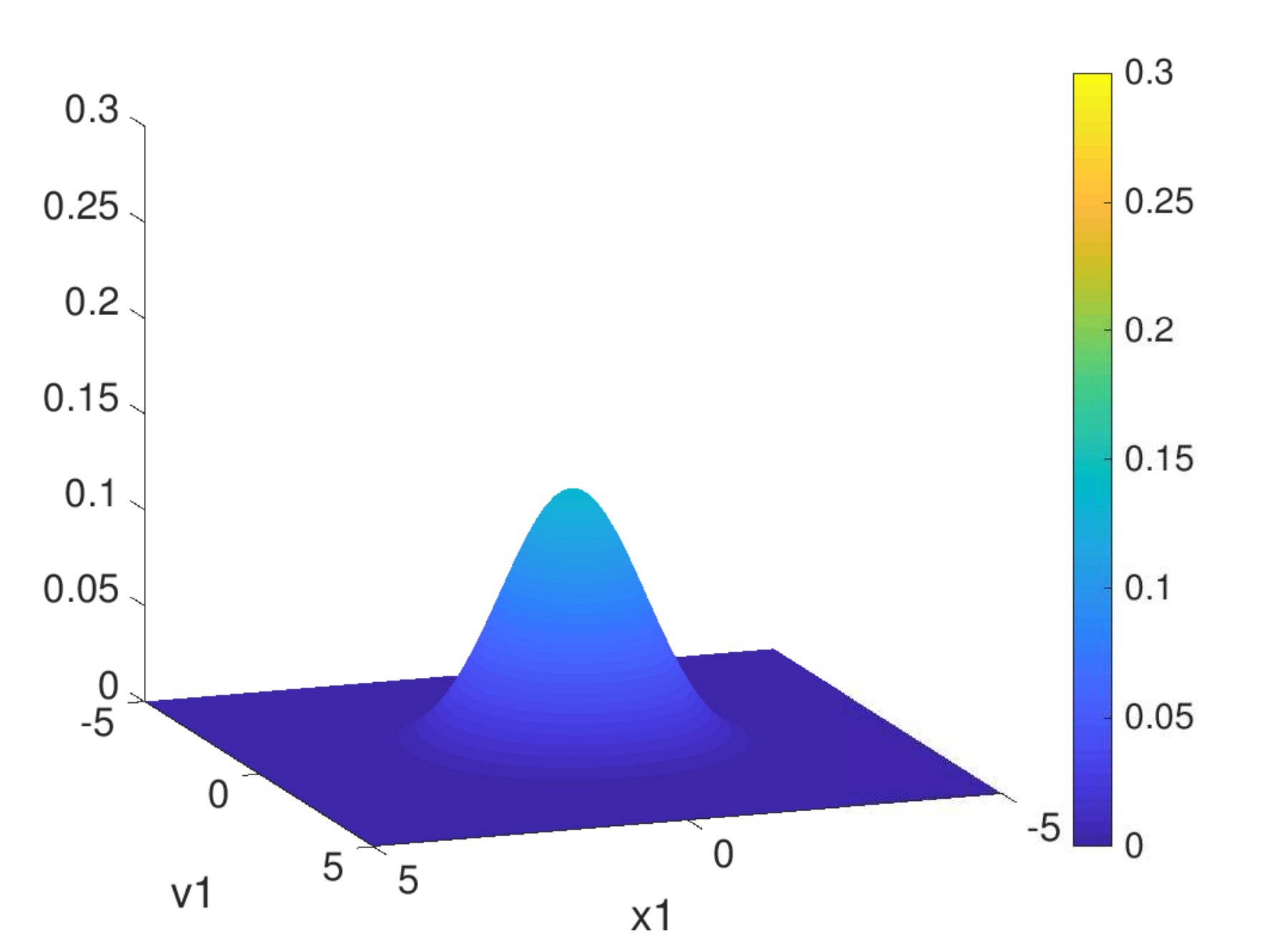}
        \caption{sparse-grid result at $t=3.0$}
        \label{4dsparse-t3}
    \end{subfigure}
    \caption{\footnotesize{Example 5(b), solution $f$ of the four dimensional Vlasov-Boltzmann transport equation by fifth order WENO scheme on sparse grids ($N_r=10$ for root grid, finest level $N_L=3$ in the sparse-grid computation) and the corresponding $80 \times 80\times 80\times 80$ single grid, at different time $t$.
    2D cuts of solutions in the $x_1-v_1$ plane at $x_2=v_2=0$.
    $CFL=0.4$. (a), (c), (e): single-grid results; (b), (d), (f): sparse-grid results. }}
\label{4dplot_vlasovBTE}
\end{figure}

\begin{figure}[p]
    \centering
    \begin{subfigure}[b]{0.48\textwidth}
        \includegraphics[width=0.9\textwidth]{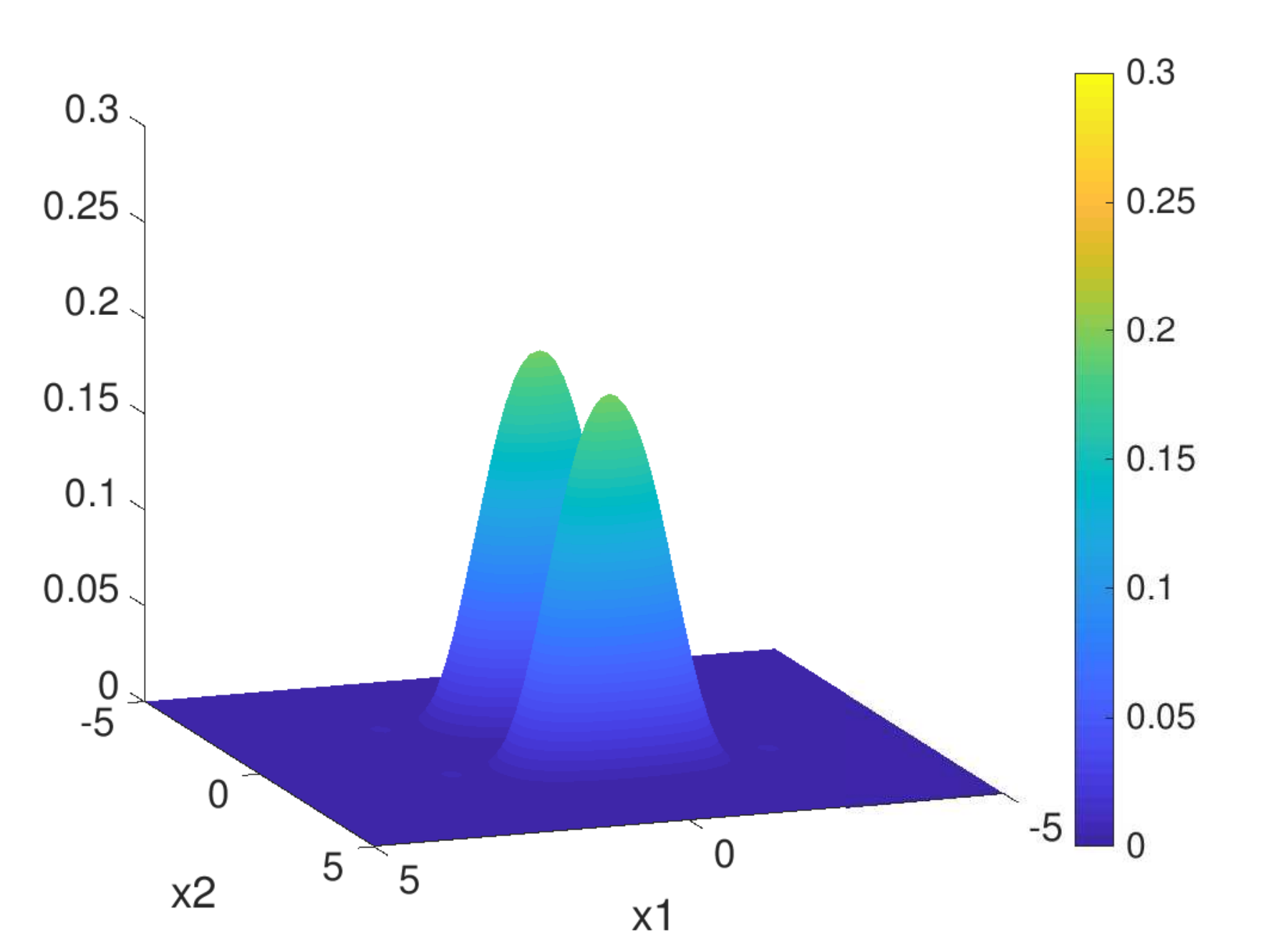}
        \caption{single-grid result at $t=0.5$}
        \label{4dsingle2-t0.5}
    \end{subfigure}
        \begin{subfigure}[b]{0.48\textwidth}
        \includegraphics[width=0.9\textwidth]{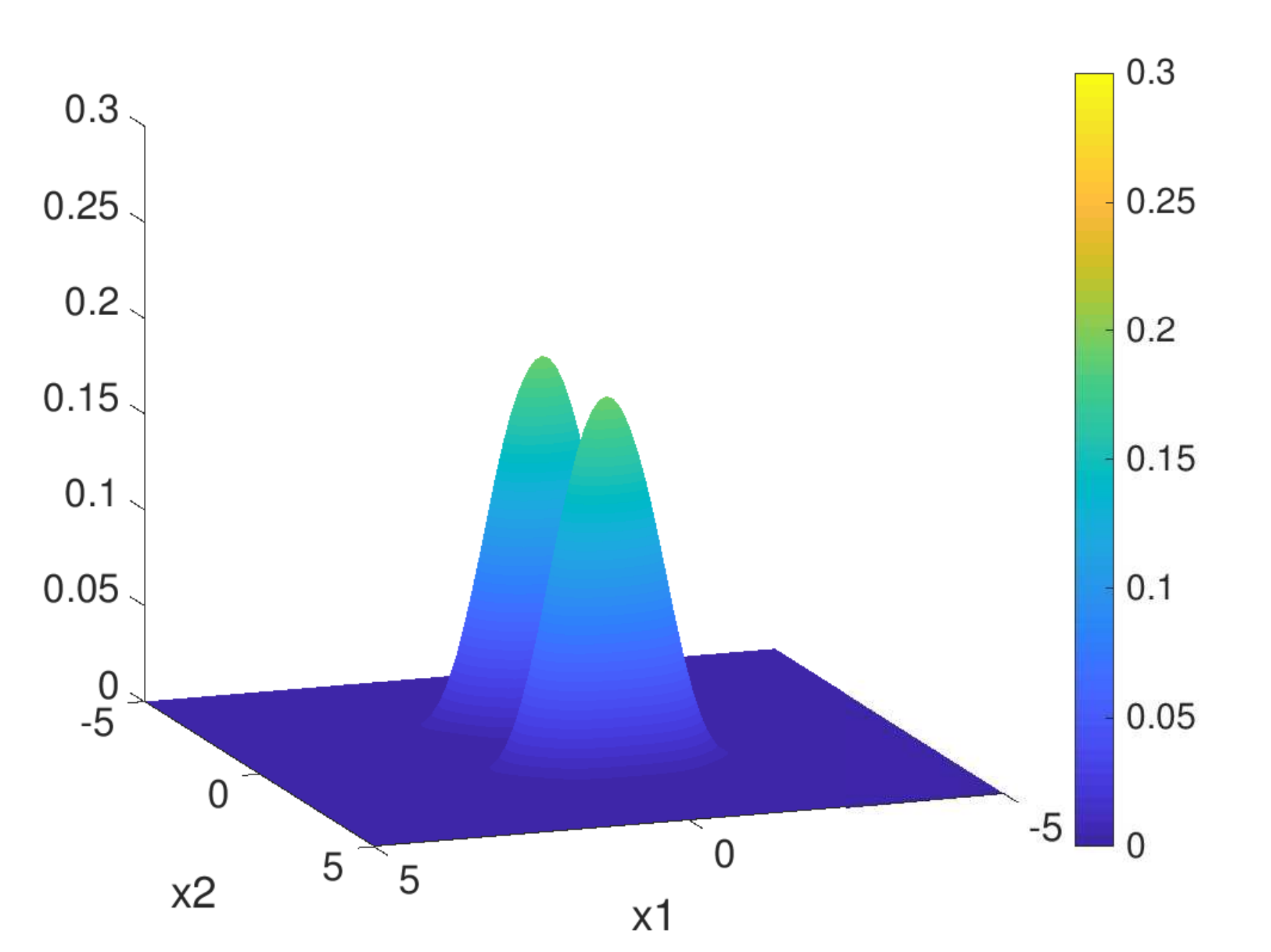}
        \caption{sparse-grid result at $t=0.5$}
        \label{4dspar2-t0.5}
    \end{subfigure}
        \begin{subfigure}[b]{0.48\textwidth}
        \includegraphics[width=0.9\textwidth]{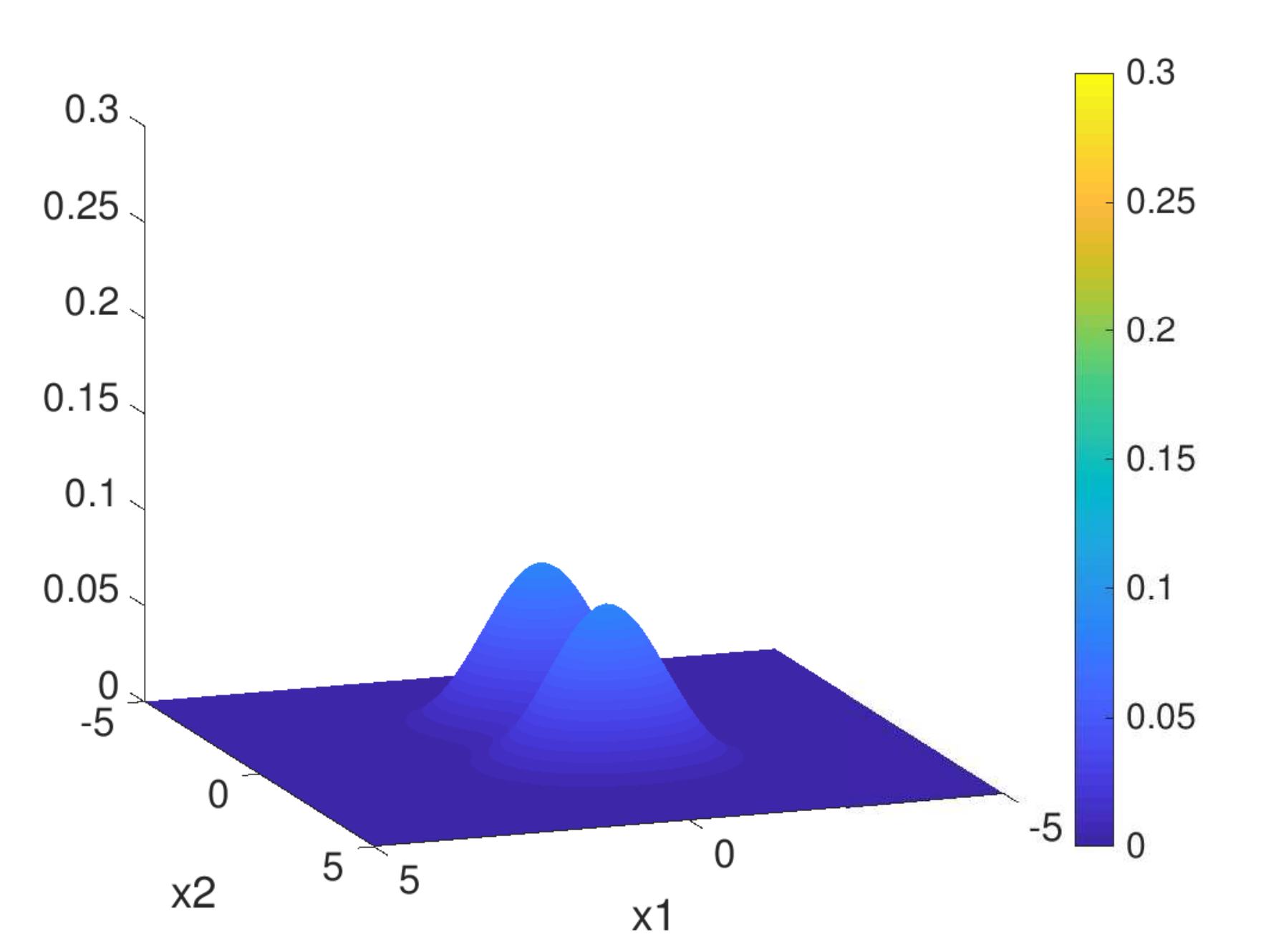}
        \caption{single-grid result at $t=1.0$}
        \label{4dsingle2-t1}
    \end{subfigure}
        \begin{subfigure}[b]{0.48\textwidth}
        \includegraphics[width=0.9\textwidth]{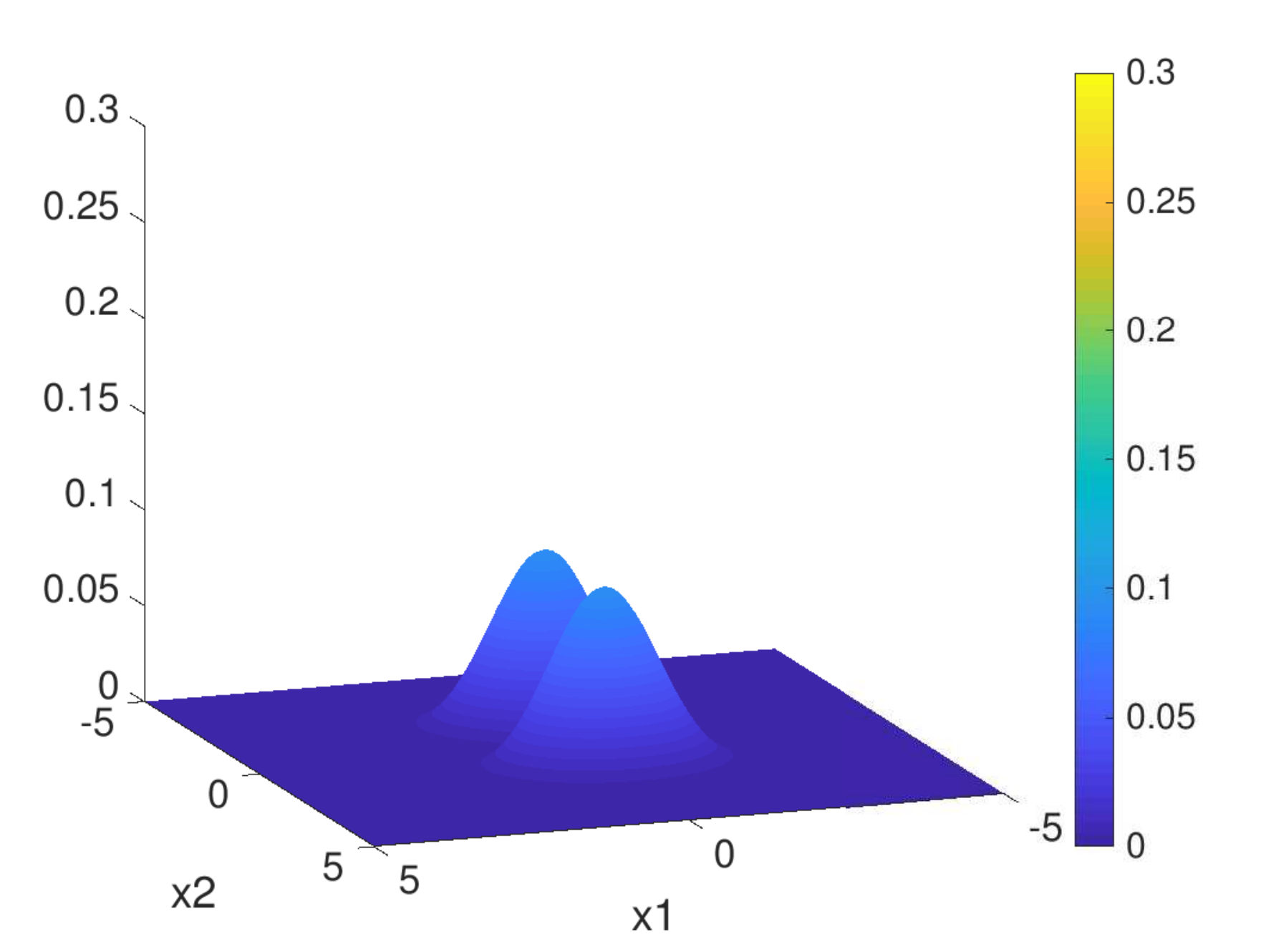}
        \caption{sparse-grid result at $t=1.0$}
        \label{4dsparse2-t1}
    \end{subfigure}
        \begin{subfigure}[b]{0.48\textwidth}
        \includegraphics[width=0.9\textwidth]{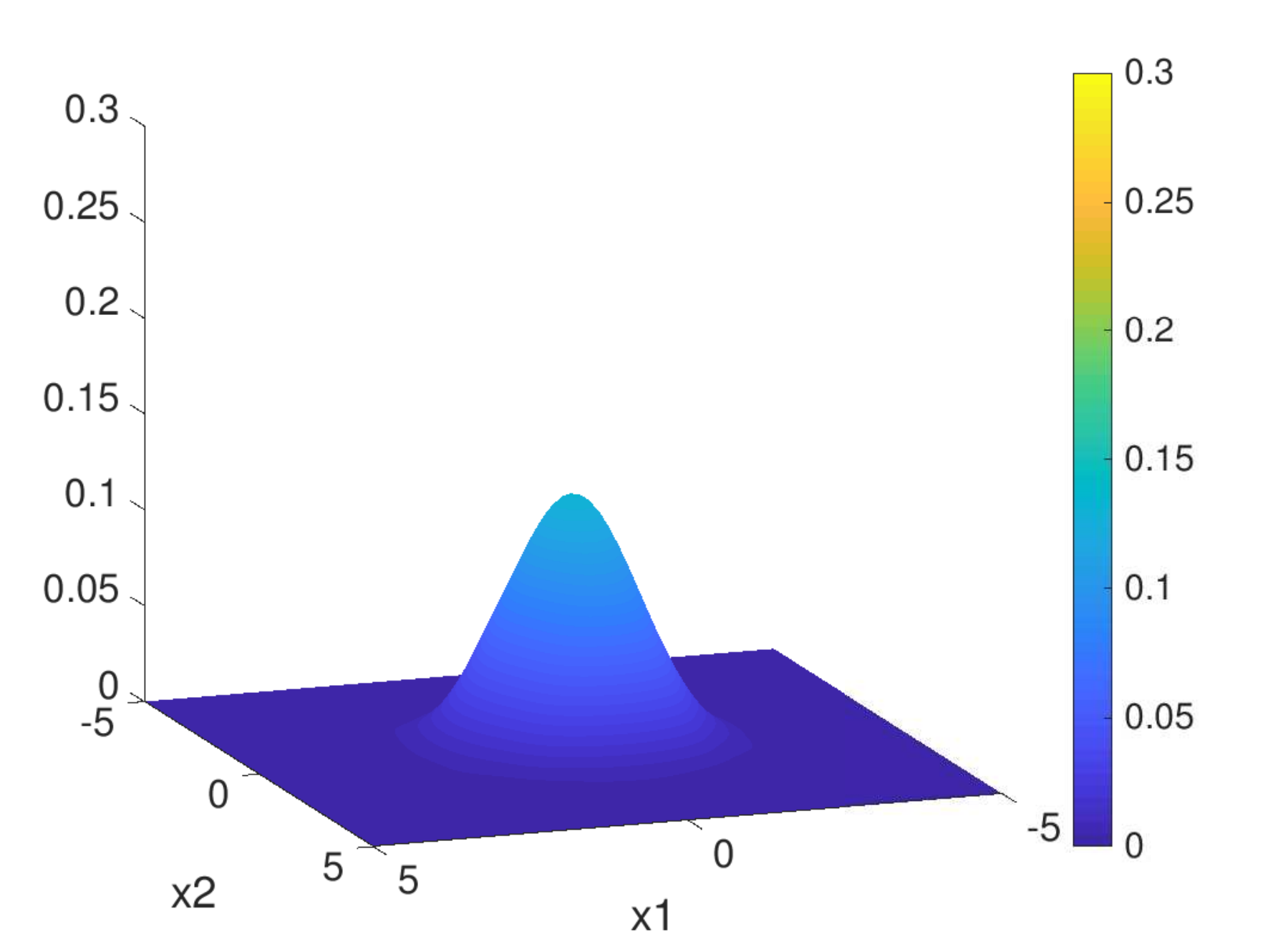}
        \caption{single-grid result at $t=3.0$}
        \label{4dsingle2-t3}
    \end{subfigure}
        \begin{subfigure}[b]{0.48\textwidth}
        \includegraphics[width=0.9\textwidth]{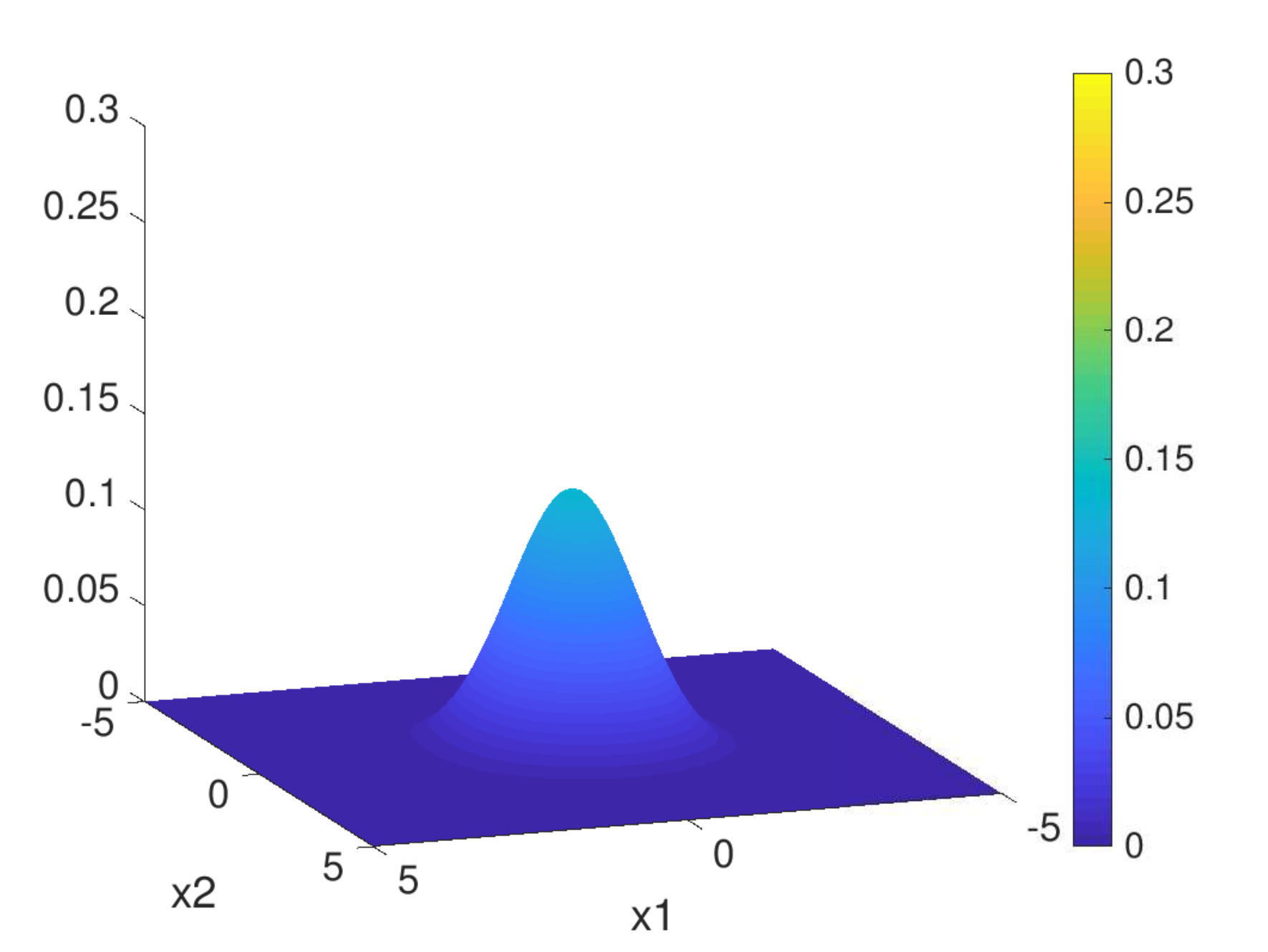}
        \caption{sparse-grid result at $t=3.0$}
        \label{4dsparse2-t3}
    \end{subfigure}
    \caption{\footnotesize{Example 5(b), solution $f$ of the four dimensional Vlasov-Boltzmann transport equation by fifth order WENO scheme on sparse grids ($N_r=10$ for root grid, finest level $N_L=3$ in the sparse-grid computation) and the corresponding $80 \times 80\times 80\times 80$ single grid, at different time $t$.
    2D cuts of solutions in the $x_1-x_2$ plane at $v_1=v_2=0$.
    $CFL=0.4$. (a), (c), (e): single-grid results; (b), (d), (f): sparse-grid results. }}
\label{4dplot2_vlasovBTE}
\end{figure}

{\bf(b) 4D case.} We solve the 4D case by the fifth order sparse grid WENO scheme. $d = 2$ and parameters $\theta=\tau=1$. The computation domain is $\Omega=[-5,5]\times[-5,5]\times[-5,5]\times[-5,5]$. The initial condition is
\begin{equation}
    f(0,x_1,x_2,v_1,v_2)=\frac{1}{s}\sin\big(\frac{x_1^2}{2}\big)^2\cos\big(\frac{x_2^2}{2}\big)^2\exp(-\frac{x_1^2+x_2^2+v_1^2+v_2^2}{2}),
\end{equation}
where $s$ is the normalization constant such that $\int_{\Omega}f(0,x_1,x_2,v_1,v_2)dx_1dx_2dv_1dv_2=1$.
Again, the zero boundary conditions are prescribed at the domain boundaries.
Sparse grids with a $10 \times 10 \times 10\times 10$ root grid
and the most refined level $3$ are used. So the corresponding single grid is $80 \times 80 \times 80\times 80$.
We compare the simulation results of sparse-grid computation and
the corresponding single-grid computation. The CFL number is taken to be $0.4$ in the simulations. In Figure \ref{4dplot_vlasovBTE} and Figure \ref{4dplot2_vlasovBTE}, we show the plots of numerical solutions of
$f$ in different 2D planes with fixed third and fourth direction coordinates. It is seen that the sparse-grid computations and
the corresponding single-grid computations generate comparable results. However, the sparse-grid computations are much
more efficient than the corresponding single-grid computations. In Table \ref{tb:ex5bcpu}, we record the simulation CPU times at different time $t$ to which the model equation is evolved. More than $93\%$ CPU times are saved by performing the fifth order WENO simulations on the sparse grids rather than the single grid for this 4D problem.

\begin{table}[htbp]\footnotesize
    \begin{tabular}{ c|c|c|c }
        \hline
        Time $t$ & CPU Time(s) of single-grid comp.  & CPU Time(s) of sparse-grid comp. &sparse/single ratio \\\hline
        0.5 & 56231.125 &  3918.468  & 0.06969 \\
        1.0 & 112539.000  & 7876.233  & 0.06999 \\
        2.0 & 225612.000  & 14468.500  & 0.06413 \\
        3.0 & 332706.250  & 22821.125  & 0.06859 \\
        \hline
    \end{tabular}
    \caption{\footnotesize{Example 5(b), CPU times (unit: seconds) for the simulations on sparse grids and the corresponding single grids, at different time $t$. }}
    \label{tb:ex5bcpu}
\end{table}

\bigskip
\noindent{\bf Example 6 (A 3D Vlasov-Maxwell system).}
In this example, we consider a 3D example of Vlasov-Maxwell system. The example is from \cite{TGC}.
It is a simplified version of the single species Vlasov-Maxwell system which has one spatial variable and two velocity variables, by assuming that the system is uniform in other variable directions of the full 6D domain. In \cite{TGC},
the system was solved by a sparse grid DG method.
The system of equations has the following form:
\begin{align*}
    &f_t+\xi_2f_{x_2}+(E_1+\xi_2B_3)f_{\xi_1}+(E_2-\xi_1B_3)f_{\xi_2}=0,\\
    &\frac{\partial B_3}{\partial t}~=~\frac{\partial E_1}{\partial x_2},\\
    &\frac{\partial E_1}{\partial t}~=~\frac{\partial B_3}{\partial x_2}-j_1,\\
    &\frac{\partial E_2}{\partial t}~=~-j_2,
\end{align*}
where $x_2$ is the spatial variable and $\xi_1,\xi_2$ are the velocity variables.
The system is defined on the domain $\Omega_x \times \Omega_\xi$. $\Omega_x$ denotes the physical space and $x_2\in \Omega_x$. $\Omega_\xi$ is the velocity space and $(\xi_1,\xi_2) \in \Omega_\xi$.
The probability distribution function of electrons $f=f(x_2,\xi_1,\xi_2,t)$. $E_1=E_1(x_2, t)$ and $E_2=E_2(x_2, t)$ are the electric field components. $B_3=B_3(x_2,t)$ is the magnetic field component.
The whole physical space has the 2D electric field $\vec{\mathbf{E}}=(E_1(x_2, t), E_2(x_2, t), 0)$ and the 1D magnetic field
$\vec{\mathbf{B}}=(0,0,B_3(x_2,t))$.
The current densities $j_1(x_2, t)$
and $j_2(x_2, t)$ are
\begin{equation}
    j_1~=~\iint\limits_{\Omega_\xi} f(x_2,\xi_1, \xi_2, t)\xi_1d\xi_1d\xi_2,~~~~j_2~=~\iint\limits_{\Omega_\xi} f(x_2,\xi_1, \xi_2, t)\xi_2d\xi_1d\xi_2.
\end{equation}
The initial condition of the system is
\begin{align*}
    f(x_2,\xi_1, \xi_2, 0)&=\frac{1}{\pi\beta}e^{-\xi_2^2/\beta}[\delta e^{-(\xi_1-v_{0,1})^2/\beta} + (1 - \delta)e^{-(\xi_1+v_{0,2})^2/\beta}],\\
    E_1(x_2, 0)&=E_2(x_2, 0)=0,~~~~~~~~B_3(x_2, 0)=b\sin(k_0x_2).
\end{align*}
The computational domain is $\Omega_x=[0, 2\pi/k_0]$ and $\Omega_\xi=[-1.2,1.2]^2$, with periodic boundary conditions applied to the system.
The parameters are taken to be $\beta=0.01, b=0.001, \delta=0.5, v_{0,1}=v_{0,2}=0.3, k_0=0.2$ as in \cite{TGC}.
Here we use this interesting 3D problem to test the efficiency of the proposed fifth order sparse grid WENO scheme in this paper. For detailed physical explanations of the system and the parameters, we refer to \cite{CPBM, TGC}.

\begin{figure}[p]
    \centering
    \begin{subfigure}[b]{0.48\textwidth}
        \includegraphics[width=0.9\textwidth]{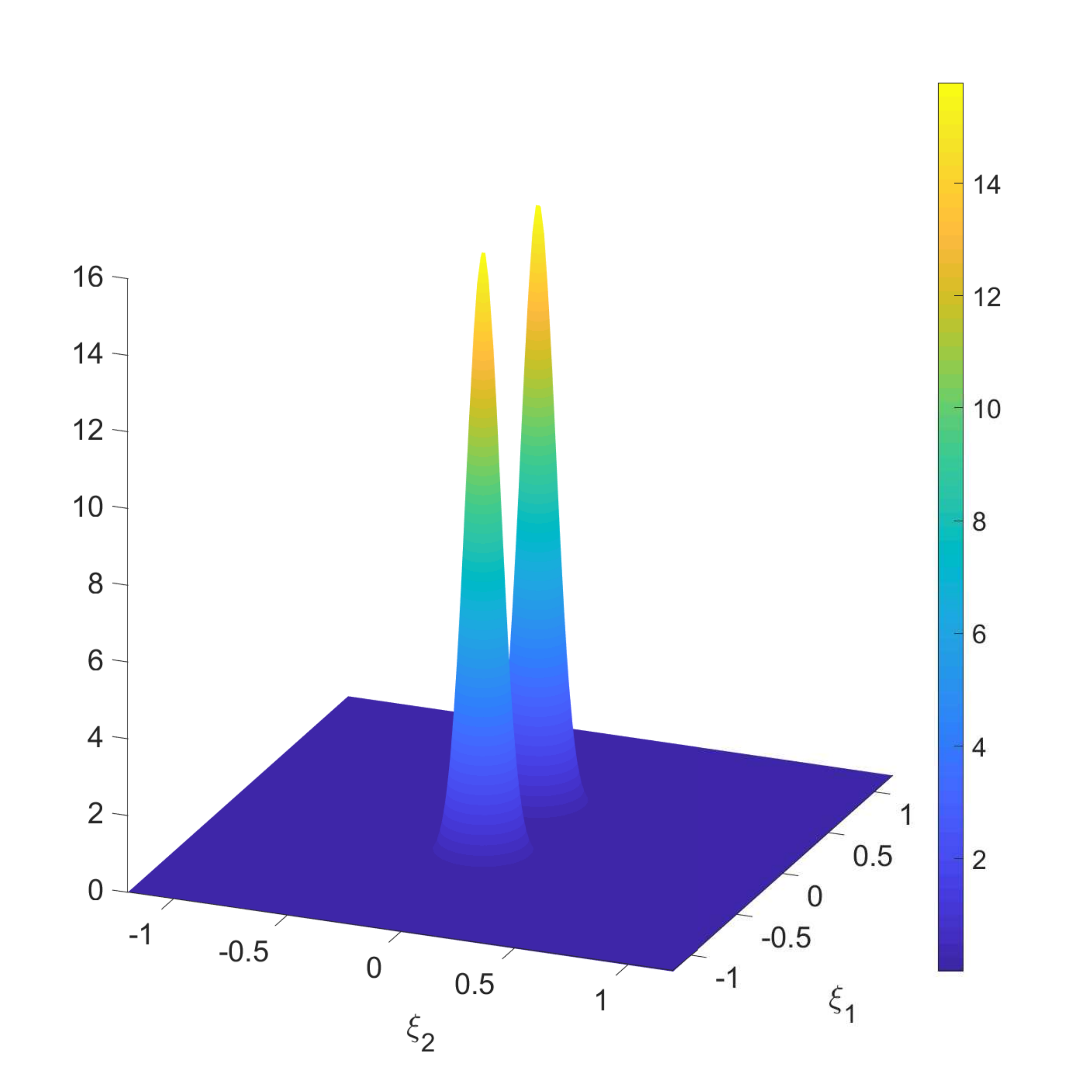}
        \caption{$x_2=5\pi$, single grid }
        \label{3dsingle1_vlamax}
    \end{subfigure}
        \begin{subfigure}[b]{0.48\textwidth}
        \includegraphics[width=0.9\textwidth]{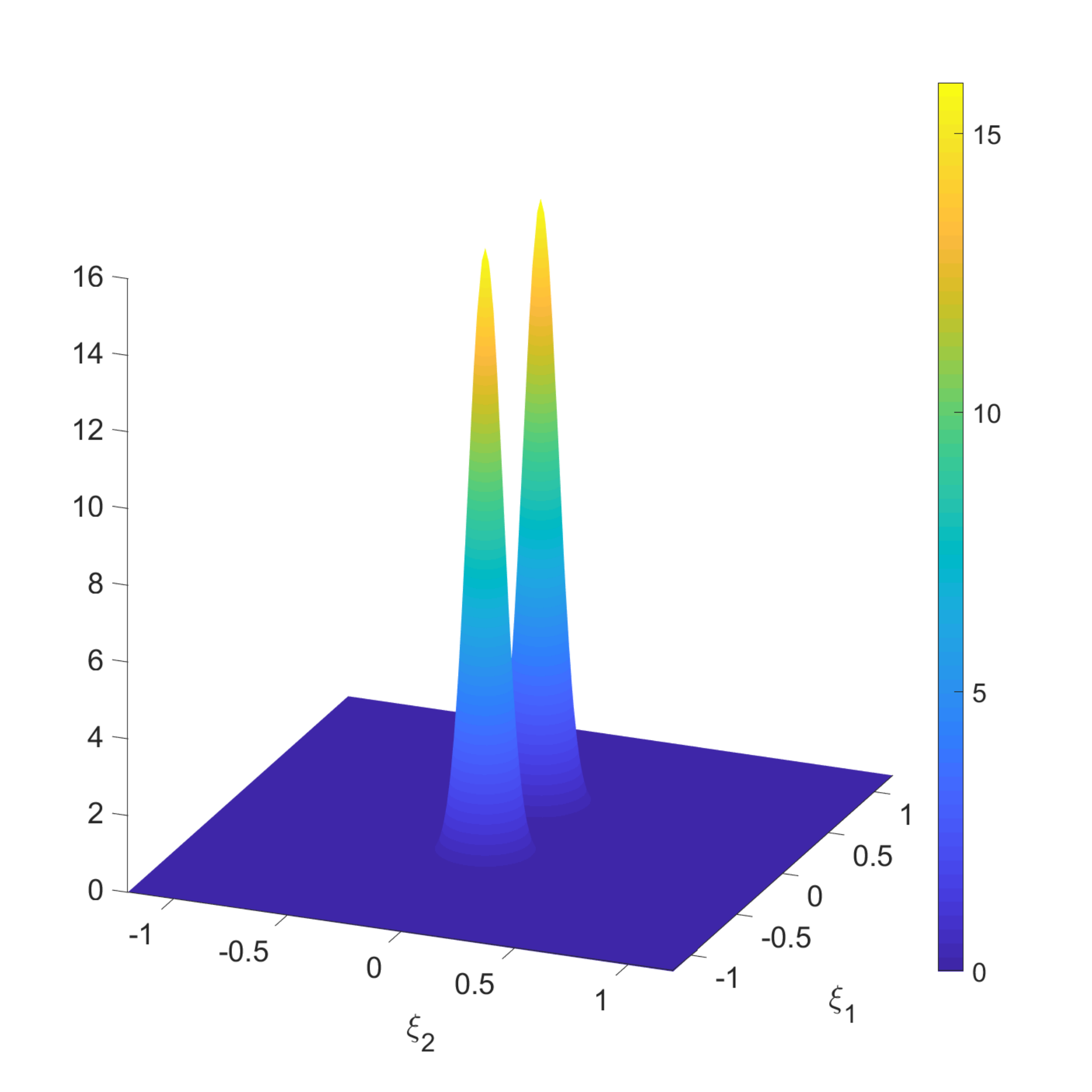}
        \caption{$x_2=5\pi$, sparse grids}
        \label{3dsparse1_vlamax}
    \end{subfigure}
        \begin{subfigure}[b]{0.48\textwidth}
        \includegraphics[width=0.9\textwidth]{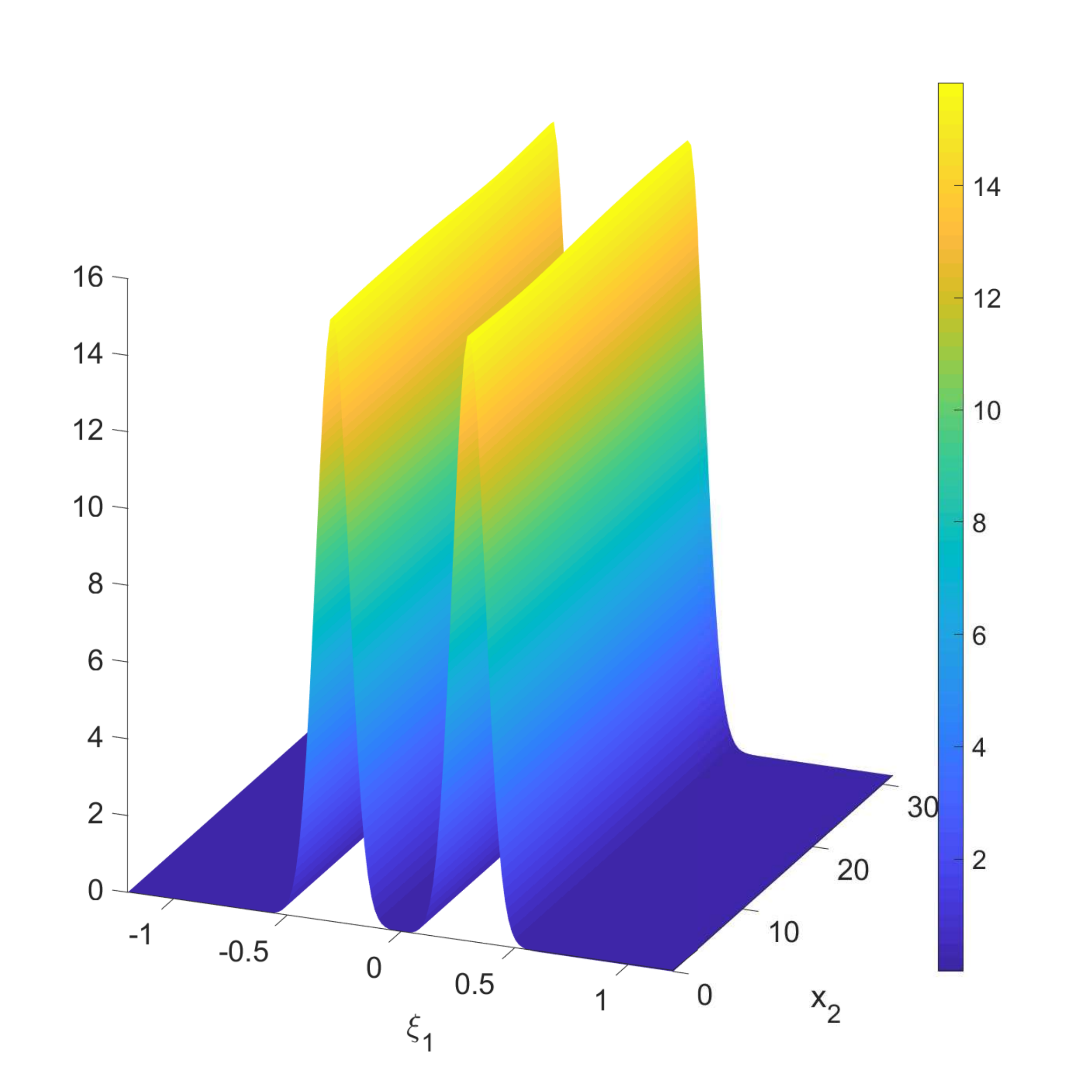}
        \caption{$\xi_2=0$, single grid}
        \label{3dsingle2_vlamax}
    \end{subfigure}
        \begin{subfigure}[b]{0.48\textwidth}
        \includegraphics[width=0.9\textwidth]{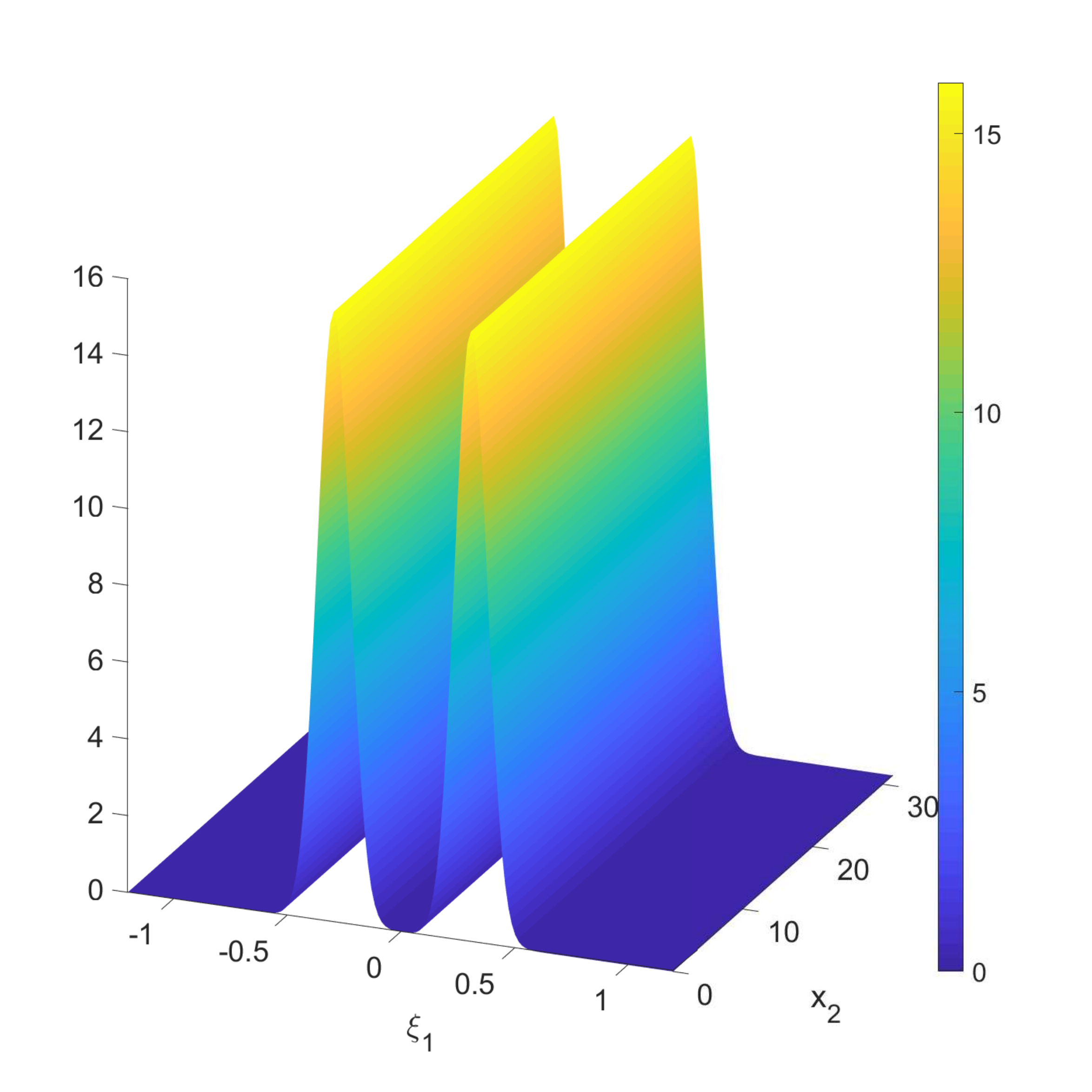}
        \caption{$\xi_2=0$, sparse grids}
        \label{3dsparse2_vlamax}
    \end{subfigure}
        \begin{subfigure}[b]{0.48\textwidth}
        \includegraphics[width=0.9\textwidth]{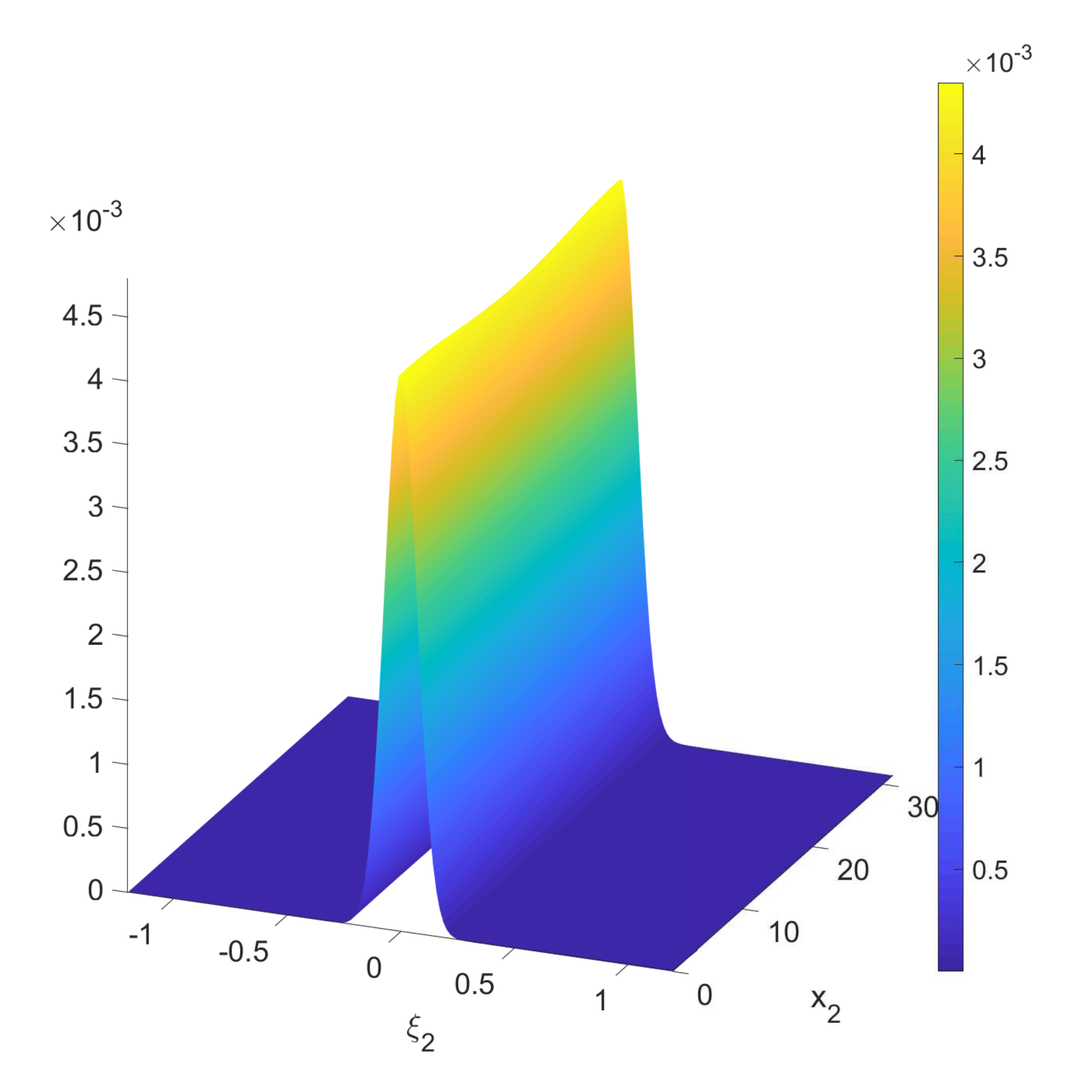}
        \caption{$\xi_1=0$, single grid}
        \label{3dsingle3_vlamax}
    \end{subfigure}
        \begin{subfigure}[b]{0.48\textwidth}
        \includegraphics[width=0.9\textwidth]{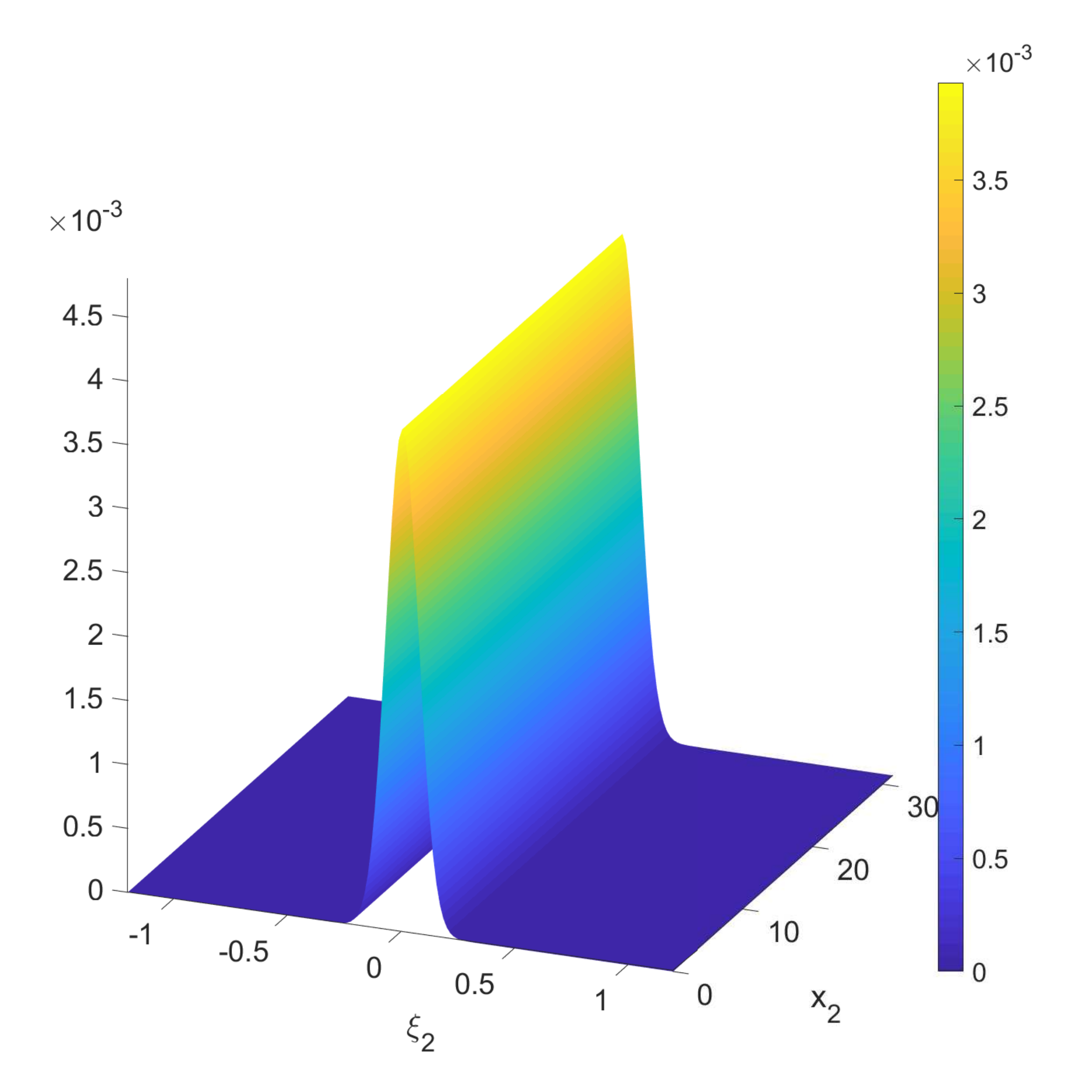}
        \caption{$\xi_1=0$, sparse grids}
        \label{3dsparse3_vlamax}
    \end{subfigure}
    \caption{Example 6, solution $f$ of the 3D Vlasov-Maxwell system at time $T=10$, by fifth order WENO scheme on sparse grids ($N_r=20$ for root grid, finest level $N_L=3$ in the sparse-grid computation) and the corresponding $160 \times 160\times 160$ single grid. Plots of solutions in 2D planes with a fixed third direction coordinate. $CFL=0.4$. }
\label{3dplot_vlasovmaxwell}
\end{figure}

We solve this 3D system till the final time $T=10$ by the fifth order sparse grid WENO scheme. Sparse grids with a $20 \times 20 \times 20$ root grid
and the most refined level $3$ are used. So the corresponding single grid is $160 \times 160 \times 160$.
The CFL number is taken to be $0.4$ in the simulations. We compare the simulation results of sparse-grid computation and
the corresponding single-grid computation. In Figure \ref{3dplot_vlasovmaxwell}, the plots of numerical solutions of
$f(x_2,\xi_1,\xi_2,10)$ in different 2D planes with a fixed third direction coordinate are presented. Similar as previous examples, we observe that the sparse-grid computation and
the corresponding single-grid computation generate comparable results. However, it only takes CPU time $3585.81$ seconds for the sparse-grid computation to finish the simulation, while CPU time $10331.04$ seconds are needed by the corresponding single-grid computation. More than $65\%$ CPU costs have been saved by performing the fifth order WENO simulations
on sparse grids in this 3D example.

\bigskip
\noindent{\bf Remark 2.} In this example we show the efficiency and CPU time saving of performing the fifth order WENO simulations
on sparse grids by solving a simplified 3D Vlasov-Maxwell system, which indicates that the sparse grid WENO scheme could be a promising method in simulating complex kinetic systems. However we would like to emphasize that the Vlasov-Maxwell system is a very complex problem and there are many challenging issues in the physics it models, its mathematical theory and numerical method developments. Although beyond the scope of this paper, more detailed and in-depth studies on the system are important and will be carried out in our future research.

\section{Conclusions and discussions}

 In this paper, we propose a general approach to perform high order accuracy WENO simulations on sparse grids for efficiently solving high dimensional hyperbolic PDEs. The broadly used fifth order finite difference WENO scheme is focused to present the method of implementing  WENO schemes on sparse grids via sparse-grid
combination technique.
Different from other sparse grid methods, in this paper novel WENO prolongations are proposed in sparse-grid
combination to obtain the robustness and high resolution properties of the fifth order WENO scheme for resolving shock waves
in solution of hyperbolic equations.
By performing the simulations of the fifth order WENO scheme
on sparse grids, we achieve a much more efficient algorithm than regular ways on single grids to solve the multidimensional hyperbolic equations.
Numerical experiments on 2D, 3D and 4D problems are carried out for the fifth order sparse grid WENO method to show its high efficiency. The savings in
computational costs are significant, especially for higher dimensional problems. For example, the fifth order sparse grid WENO simulations on a kinetic problem modeled by 4D Vlasov equation show more than $93\%$ CPU time savings, by comparing with the corresponding single-grid computations.
At the same time, comparable accuracy and resolution of the fifth order WENO scheme
to that of computations on regular single grids are obtained for sparse-grid computations on relatively refined meshes.

In this paper, we focus on the algorithm development and its numerical experiments for the nonlinear fifth order WENO scheme on sparse grids.
There are quite a few open problems to be investigated further for the method. The sparse-grid combination technique has been studied and understood well for linear schemes. For example, a linear error analysis of the sparse-grid combination technique for a linear advection equation solved by an upwind scheme is performed in
\cite{LKV1}. However for nonlinear schemes such as high order WENO schemes, it will be interesting to study
how to perform theoretical error analysis to find the relationship of numerical errors among different sparse grids. In our numerical experiments in this paper, it is found that on relatively coarse meshes, the fifth order sparse grid WENO method has larger numerical errors than the computations on corresponding single grids.
How to improve the accuracy of the sparse grid WENO scheme on coarser meshes is also an open problem.
Another topic is to extend the method to solve more complex problems with strong shock waves in applications such as compressible fluid dynamics, which needs to resolve the possible issue of bound-preserving in the proposed sparse grid methods.
These will be our future work.

\end{document}